\documentclass[11pt,letterpaper]{article}

\usepackage[T1]{fontenc}
\usepackage{lmodern}
\usepackage[margin=1in]{geometry}
\usepackage{microtype}
\usepackage{amsmath,amssymb,amsthm}
\usepackage{array,booktabs,multirow}
\usepackage{graphicx,subcaption,pdflscape}
\usepackage{textcomp}
\usepackage{algorithm}
\usepackage[noend]{algpseudocode}
\usepackage[authoryear,round]{natbib}
\usepackage[hidelinks]{hyperref}
\newcommand{\PaperTitle}{Doing Well by Doing Good: Alignment of Revenue and Welfare in Robust Product Line Design through Data-Driven Optimization}
\hypersetup{
  pdftitle={\PaperTitle},
  pdfauthor={Maoqi Liu, Liangjie Li, Yun Hui Lin, Xiuli Chao, Jiajie Zhang},
  pdfkeywords={Random Utility Maximization; Distributionally Robust Optimization; Welfare Maximization}
}

\theoremstyle{plain}
\newtheorem{theorem}{Theorem}
\newtheorem{lemma}{Lemma}
\newtheorem{proposition}{Proposition}
\newtheorem{corollary}{Corollary}
\newtheorem{aplemma}{Lemma}

\theoremstyle{definition}
\newtheorem{definition}{Definition}
\theoremstyle{remark}
\newtheorem{remark}{Remark}

\newcommand{\vhat}{\widehat v}
\newcommand{\Shat}{\widehat {\mathcal S}_{N}}

\newcolumntype{L}[1]{>{\raggedright\arraybackslash}p{#1}}
\newcommand{\epsdag}{\epsilon^{\dagger}}
\newcommand{\Aplus}{\mathcal{A}^{+}}
\newcommand{\IAPS}{\mathcal{I}_{\lambda}}
\newcommand{\Mail}[1]{\href{mailto:#1}{\nolinkurl{#1}}}

\begin{document}
\begin{center}
  {\LARGE\bfseries\PaperTitle\par}
  \vspace{1em}
  {\large
    Maoqi Liu\textsuperscript{1}, Liangjie Li\textsuperscript{1},
    Yun Hui Lin\textsuperscript{2}, Xiuli Chao\textsuperscript{3},
    Jiajie Zhang\textsuperscript{4,5,*}\par}
  \vspace{0.6em}
  {\small
    \textsuperscript{1}School of Management, Shandong University\\
    \textsuperscript{2}School of Economics and Management, Dalian University of Technology\\
    \textsuperscript{3}Department of Industrial and Operations Engineering, University of Michigan, Ann Arbor\\
    \textsuperscript{4}School of Statistics and Data Science, Southeast University\\
    \textsuperscript{5}Institute of Operations Research and Analytics, National University of Singapore\par}
  \vspace{0.6em}
  {\footnotesize
    \Mail{liumq@sdu.edu.cn}\enspace
    \Mail{mogoooojay@gmail.com}\enspace
    \Mail{linyh@dlut.edu.cn}\enspace
    \Mail{xchao@umich.edu}\\[2pt]
    \textsuperscript{*}Corresponding author: Jiajie Zhang, \Mail{jiajiez@u.nus.edu}\par}
\end{center}
\medskip

\noindent\textbf{Abstract}\par\nobreak\smallskip

\textbf{Problem definition.}
Product line design (PLD) determines which variants a firm should offer in order to maximize expected revenue. A candidate line is usually evaluated with a random utility maximization (RUM) model, in which choice probabilities are determined by the joint distribution of latent customer utilities. Because that distribution is not directly observable, it must be inferred from choice data, which leads to inevitable distributional uncertainty. We study how to protect the PLD decision against that uncertainty.

\textbf{Methodology/results.}
We propose a distributionally robust optimization (DRO) model that maximizes worst-case expected revenue over an $\ell_1$-Wasserstein ball centered at a nominal utility distribution. We first derive a dual reformulation for a general nominal and then specialize the nominal in two ways. When utility shocks are independent Gumbel, the dual reduces to a closed-form robust multinomial logit. We bound the revenue gap when the parameters are estimated from finitely many full-assortment choices. When the  model specifics are not known a priori, 
we adopt a data-driven approach by applying the empirical distribution of sampled utility vectors as an approximation and extend the optimality-gap analysis to a general RUM oracle.
We develop Adaptive Parametric Search (APS) for the NP-hard
empirical problem, with a certified relative optimality guarantee.
Numerical experiments show faster computation and improved
out-of-sample performance.

\textbf{Managerial implications.}
Analysis of the dual shows that protecting revenue against preference uncertainty creates an endogenous incentive to raise buyers' welfare. The incentive takes a different form in each implementation: robustness favors larger assortments in the unconstrained logit model and, under the empirical nominal, products with lower prices and greater utility relative to the outside option. Numerical results further show a ``doing well by doing good'' effect: a product line with higher in-sample buyers' welfare delivers higher out-of-sample revenue and buyers' welfare.

\smallskip
{\small\noindent\textbf{Keywords:} Random Utility Maximization, Distributionally Robust Optimization, Welfare Maximization.\par}

\section{Introduction}\label{sec:intro}

The product line design (PLD) problem concerns the selection of a limited subset of related product variants from a large candidate set to serve a heterogeneous market \citep{zufryden1977conjoint,10.2307/183706}. 

A key factor in the success of PLD is understanding heterogeneous customer preferences and predicting choice behavior across alternative product lines. Random utility maximization (RUM) models provide a standard framework to analyze this choice behavior. Specifically, it assumes that a random customer assigns a utility for every available alternative and selects the alternative with the highest realized utility. Therefore, the utility distribution over the population determines the choice probability of each alternative---the probability that it realizes the highest utility---which in turn allows the firm to evaluate the expected revenue of a product line. 

A central challenge in constructing an RUM model is the inevitable uncertainty in customer utilities. This is because neither realized utilities nor their underlying distribution is directly observable and therefore both must be inferred from data such as transactions or conjoint surveys. More specifically, the error jointly reflects the inconsistency in respondents' choices and the limited number of observations for each customer, leading to substantial estimation error in the inferred utilities \citep{Toubia2007Probabilistic,bertsimas2016robust}.

As such, it is essential to design product lines that remain robust to uncertainty in the utility distribution. In this paper, we study the robust revenue-maximizing PLD problem under a general RUM framework. Specifically, we propose a distributionally robust optimization (DRO) model that maximizes the worst-case expected revenue over an $\ell_1$-Wasserstein ball centered at the nominal utility distribution. Analysis of the dual reformulation reveals a general mechanism for robust PLD: protecting revenue against customer preference uncertainty creates an endogenous incentive to increase buyers' welfare, defined as the expected utility customers obtain from their chosen alternatives. This finding is surprising because revenue and buyers' welfare are typically viewed as competing objectives: revenue-oriented designs may favor higher-priced or less-preferred variants, thereby reducing customer welfare \citep{mussa1978monopoly}. In contrast, our results show that, under preference uncertainty, buyers' welfare serves not as a competing objective but as an insurance of revenue against customer preference uncertainty.

Two representaive implementations are then considered. 
We first study the problem when the random utilities  follows  
independent Gumbel shocks, yielding a closed-form robust multinomial logit (MNL) model. We show that in the unconstrained MNL setting, this structure further preserves the revenue-order property: an optimal robust product line can be found within the class of revenue-ordered assortments, substantially simplifying the resulting optimization problem. We quantify the optimality gap when the parameters are estimated from finitely many full-assortment choices.
To accommodate more general choice behaviors, we then present a nonparametric implementation that represents the distribution of the complete utility vector directly by its empirical distribution. Its objective can be evaluated from sampled utility vectors and therefore can approximate any general RUM models, including those for which choice probabilities do not admit a closed-form expression. We establish a finite-sample bound linking this empirical first-choice formulation to the product line design problem under oracle distribution following arbitrary RUM. 

The two implementations translate the general welfare mechanism into distinct operational principles. In the unconstrained MNL setting, the DRO model admits an optimal revenue-ordered product line that expands as the extent of distributional uncertainty increases. Under the empirical first-choice representation, robustness favors a more cost-effective product line that delivers greater customer utility relative to the revenue extracted.

\subsection{Main Contributions}\label{sec:intro:contrib}
Our main contributions are summarized as follows.
\begin{itemize}
\item \textbf{A general Wasserstein DRO framework for revenue-maximizing PLD.} We formulate a robust PLD model under a general RUM by maximizing worst-case expected revenue over an $\ell_1$-Wasserstein neighborhood of the nominal joint utility distribution. Unlike approaches that perturb only a prespecified set of parameters or scenarios, our framework operates directly on the utility distribution of general RUM models. Performance guarantee of two specific implementations are established. Extensive numerical experiments demonstrate the superior performance of the proposed DRO model.

\item \textbf{An endogenous alignment between revenue robustness and buyers' welfare.} Analysis of the dual reformulation reveals that protecting revenue against adverse utility shifts creates an endogenous incentive to increase buyers' welfare. Specifically, greater aversion to distributional uncertainty raises the guaranteed lower bound on the buyers' welfare delivered by the selected product line, thereby making higher-welfare product lines increasingly attractive as the extent of robustness increases. Under specific RUM specifications, this general mechanism translates into distinct operational principles that benefits customers. In the unconstrained MNL setting, robustness favors larger assortments; with an empirical nominal distribution, it favors more cost-effective product lines composed of products with lower prices and greater marginal utility relative to the outside option. Our numerical results on a semi-synthetic case study shows that the ``goodness'' to benefit more to the customer eventually pays back in the way of better out-of-sample revenue and buyers' welfare, and therefore, yields a ``doing well by doing good'' effect in the robust PLD problem.

\item \textbf{A certified solution method for empirical DRO.}
When the joint utility distribution is unknown, we use the
empirical distribution of sampled utility vectors as the nominal.
For the resulting cardinality-constrained problem, we develop
Adaptive Parametric Search (APS). Given a prescribed tolerance
$\delta\in(0,1)$, APS certifies that the returned product line
attains at least a $(1-\delta)$ fraction of the empirical robust
optimum. Across 360 synthetic instances, APS meets its prescribed
tolerance in every run and, at $\delta=0.05$, is $70.7\times$
faster on average and up to $662.4\times$ faster than a direct
mixed-integer linear reformulation.

\end{itemize}

\subsection{Related Literature}\label{sec:LR}

Our work is most closely related to the literature on robust product line design and robust assortment planning. Because these two problems share closely related mathematical structures in our setting, we discuss them jointly and use the term product line design throughout for consistency. 

Existing robust PLD models differ primarily in their assumptions about customer choice and the structure of the uncertainty set. Although a growing literature studies robust PLD under emerging choice models, such as ranking-based choice \citep{sturt2021value} and Markov chain choice \citep{desir2021robust}, we focus on RUM-based models. RUMs offer a clear economic interpretation of customer choice and allow the outcomes of product-line decisions to be evaluated from both the seller’s perspective, through revenue, and the customer’s perspective, through buyers’ welfare.

Within the RUM framework, most studies focus on either the first-choice model or the multinomial logit (MNL) model. For the first-choice model, \citet{Wang2012A} allow each respondent's estimated partworth utilities to vary within a box uncertainty set and maximize the resulting worst-case market share. \citet{liu2023explainable} propose a related formulation in which each respondent's utility vector varies within a more general polyhedral uncertainty set constructed from conjoint-analysis responses. \citet{Wang2021Robust} adopt a scenario-based framework in which each scenario represents a possible first-choice model and optimize worst-case market share or expected revenue across scenarios.  \citet{liu2023from} introduce a Wasserstein DRO framework for the first-choice model, but focus on single-product design and share-of-choice maximization. For the MNL model, \citet{rusmevichientong2012robust} study robust assortment optimization when the attraction parameters are allowed to vary within a prescribed uncertainty set. \citet{jin2022distributionally} instead adopt a distributionally robust approach by placing a Kullback--Leibler divergence ball around a nominal MNL model.
Moreover, \citet{bertsimas2016robust} jointly address parameter and model uncertainty. They construct multiple uncertainty sets around different nominal choice models, including first-choice and mixture-logit specifications, and maximize worst-case revenue over their union.

Most of these studies focus on improving performance under uncertainty but provide limited insight into the economic mechanisms through which robustness reshapes product-line decisions. Although \citet{rusmevichientong2012robust,jin2022distributionally} show that revenue-ordered assortments remain optimal and that greater uncertainty aversion leads to larger assortments, these results are specific to MNL-based models. Appendix~\ref{appendix:capacity_separation} shows that these properties need not survive a cardinality constraint, and that they do not extend to a general RUM. In contrast, our analysis shows that the underlying driver is an endogenous emphasis on buyers' welfare induced by revenue robustness, which not only explain why a larger assortment is required for a more robust product line but also persists under a general RUM framework and cardinality constraints.

\citet{liu2023from} is the most closely related study. From the perspective of share-of-choice maximization, they show that a Wasserstein DRO model under the first-choice framework also places an endogenous emphasis on buyers' welfare.
Our work differs from theirs in four respects.
First, they maximize share-of-choice, the fraction of customers who select any of the offered products rather than the outside option. Product prices are therefore ignored or treated as identical, and the objective only requires knowing whether a customer prefers some offered product to the outside option. We maximize worst-case expected revenue. Because products earn different revenues, the firm must identify the product each customer chooses and balance higher-priced products chosen by fewer customers against lower-priced products chosen by more customers. Revenue maximization therefore creates a sharper conflict with buyers' welfare: a larger share of customers who buy is only sometimes consistent with higher buyers' welfare. The resulting price--utility trade-off yields a further implication of the welfare mechanism. Under the empirical model, greater robustness favors products that deliver more utility per unit of revenue.
Second, \citet{liu2023from} design a single product, whereas we select a product line from a candidate set. The number of feasible lines grows exponentially with the cardinality of the line, a computational difficulty that their analysis leaves unaddressed. We develop a set-branching procedure with modular bounds and an adaptive parametric search that returns a feasible line whose value is at least \((1-\delta)\) times the empirical optimum. Across 360 instances with up to several hundred products and customers, the algorithm certifies every solution and, at \(\delta=0.05\), is \(70.7\times\) faster on average than the mixed-integer formulation.
Third, their nominal distribution is the empirical distribution of sampled utilities. We formulate the DRO model for an arbitrary RUM nominal, which covers standard parametric specifications and broadens the scope of the structural results. Under an MNL nominal, the robust objective has a closed form, an optimal product line is revenue-ordered, and, without a cardinality constraint, that line expands as the robustness radius increases.
Fourth, we provide finite-sample guarantees for the revenue loss when the nominal distribution is estimated from data, which is not investigated by \cite{liu2023from}.

The literature on semiparametric choice models also addresses uncertainty in the utility distribution \citep{natarajan2009persistency,mishra2014theoretical,feng2017relation}. These studies share our motivation that the joint distribution of customer utilities may be imperfectly known. However, they address this uncertainty by selecting the distribution that maximizes the expected utility of the customer’s chosen alternative. This optimization is embedded in the construction of the choice model and is therefore independent of the firm’s downstream objective. In contrast, our model selects the worst-case utility distribution according to the firm’s revenue objective, so the effect of distributional uncertainty depends explicitly on both product revenues and the product-line decision.

Finally, another relevant stream of literature incorporates buyers' welfare into product line design. Buyers' welfare has long been considered as an explicit objective, dating back at least to \citet{KohliHeuristics1990}. For a recent review of this literature, we refer readers to \citet{sumida2021revenue}. In these studies, revenue and buyers' welfare are typically viewed as conflicting objectives, and the resulting models are designed to balance the trade-off between them. In contrast, our model begins from a purely revenue-maximization perspective and shows that robustness can endogenously align the firm's revenue objective with buyers' welfare.

\subsection{Organization}\label{sec:intro:org}

Section~\ref{sec:MF} establishes the classical product line design problem under a general random utility model. Section~\ref{sec:analysis} derives the Wasserstein robust model for a general nominal and specializes it to the multinomial logit, including the expected revenue loss when the logit parameters are estimated from choice data. Section~\ref{sec:DRO_FC} develop a tractable implementation of any RUM based on approximating the nominal utility distribution with its empirical counterpart and derives the corresponding expected revenue guarantee. Section~\ref{sec:delta_optimal} presents the certified algorithm for the empirical problem. Section~\ref{sec:experiments} evaluates the algorithm against the mixed-integer benchmark, and Section~\ref{sec:Case} reports the conjoint case study. Section~\ref{sec:Con} concludes the paper and proposes further research directions. Proofs and supporting formulations appear in the appendices.

\section{The Classical Model}\label{sec:MF}

In this section, we introduce the classical non-robust PLD model as a benchmark for the analysis that follows. Throughout the paper, sets are denoted by calligraphic letters and their cardinalities by the corresponding uppercase letters (e.g., \(|\mathcal{A}|=A\)). The symbol \(:=\) indicates a definition. Vectors and matrices are written in boldface, and \(\mathbb{I}\{\cdot\}\) denotes the indicator function.

Following a large body of prior work (e.g., \citealp{10.2307/183706}; \citealp{Chen2000Mathematical}; \citealp{bertsimas2019exact}), we formulate the PLD problem as the selection of a product line \(\mathcal{S}\subseteq\mathcal{A}\) from a set of candidate products \(\mathcal{A}=\{1,2,\ldots,A\}\). Each product \(i\in\mathcal{A}\) generates revenue \(r_i>0\), and we write \(\bar r=\max_{i\in\mathcal A}r_i\). The outside option is indexed by \(0\) and yields zero revenue, i.e., \(r_0=0\). Let \(\mathcal{L}\) denote the family of feasible product lines that satisfy the firm's structural or managerial constraints. For example, if the firm may offer at most \(p\) products, then \(\mathcal{L} = \{ \mathcal{S} \subseteq \mathcal{A} : |\mathcal{S}| \le p \}\).

The choice model describes how customers select among the offered products and the outside option and therefore forms the basis for evaluating a product line. We adopt the random utility maximization framework. Under this framework, each customer is associated with a random utility vector \(\boldsymbol{U}=(U_0,U_1,\ldots,U_A)\), where the randomness captures preference heterogeneity across customers. Given a product line \(\mathcal{S}\), the customer chooses the alternative with the highest utility from \(\mathcal{S}\cup\{0\}\). 
Let \(Q\) denote the population distribution of \(\boldsymbol{U}\). For any \(i\in\mathcal{S}\cup\{0\}\), the choice probability of alternative \(i\) under product line \(\mathcal{S}\) is \(q_i(\mathcal{S};Q)=\mathbb{P}_{Q}\!\left(U_i \ge U_j,\ \forall j\in\mathcal{S}\cup\{0\}\right)\).
For \(i\in\mathcal{S}\), this probability represents the expected market share of product \(i\)\footnote{To define realized revenue on the entire utility space, any residual tie is resolved in favor of the utility-maximizing alternative with the lowest revenue.}.
Therefore, the expected revenue of product line \(\mathcal{S}\) under the population distribution \(Q\) is \(R(\mathcal{S};Q)=\sum_{i\in\mathcal{S}} r_i q_i(\mathcal{S};Q)\), and the population PLD problem is
\begin{equation}
    \max_{\mathcal{S}\in\mathcal{L}} R(\mathcal{S};Q).\label{eq:nominal model}
\end{equation}

The choice probabilities and the resulting expected revenue depend critically on the underlying distribution of the utility vector \(\boldsymbol{U}\). Different distributional assumptions give rise to different random utility choice models. For example, if the utility shocks are independently and identically distributed according to a type-I extreme value distribution, the RUM reduces to the multinomial logit (MNL) model; if they follow a multivariate normal distribution, the resulting model is the multinomial probit model.

However, the population distribution \(Q\) is not directly observable and must be inferred from limited data, yielding a proxy distribution \(\widehat{Q}\) for optimization. A parametric approach, such as the MNL model, imposes a distributional structure and estimates its parameters, gaining parsimony at the risk of model misspecification. In contrast, a nonparametric approach, such as an empirical approximation based on i.i.d. samples, avoids imposing a specific distributional form but may require a large sample to accurately approximate the population distribution and is therefore susceptible to sampling error when data are limited. In either case, the proxy \(\widehat{Q}\) used for optimization may differ from the true population distribution \(Q\). This discrepancy motivates a robust formulation that accommodates both representations while protecting product-line decisions against ambiguity in the underlying customer preference distribution.

\section{The Distributionally Robust Optimization Model and Analysis}
\label{sec:analysis}

We now replace the population objective in Problem~\eqref{eq:nominal model} with a distributionally robust counterpart centered at a nominal utility distribution. Section~\ref{robust reformulation} develops the general model and its dual reformulation without imposing a parametric structure on the nominal. Section~\ref{sec:DRO_MNL} specializes the nominal to the multinomial logit.

\subsection{The DRO Model with a General RUM Nominal Distribution}\label{robust reformulation}

As discussed above, the proxy distribution \(\widehat{Q}\) may differ from the population distribution \(Q\) because of sampling and estimation errors. We account for this discrepancy by evaluating each product line over a Wasserstein neighborhood centered at a nominal distribution \(Q\). In the implementations below, this nominal is the proxy \(\widehat{Q}\) or an empirical distribution constructed from it. The Wasserstein distance between two distributions measures the minimum transportation cost required to move probability mass from one distribution to the other.

\begin{definition}[Wasserstein Metric, \citealp{esfahani2018data}] \label{WasMetric}
Let \(\mathcal M(u)\) be the set of probability measures supported on \(u\), and let \(\Pi(P,Q)\) collect their joint distributions with marginals \(P\) and \(Q\). The Wasserstein distance is
\[
d_W(P,Q)=\inf_{\Pi\in\Pi(P,Q)}
\int_{u^2}\|\boldsymbol{U}-\tilde{\boldsymbol{U}}\|
\,\Pi(d\boldsymbol{U},d\tilde{\boldsymbol{U}}),
\]
where \(\|\cdot\|\) is the chosen ground norm.
\end{definition}

Unless stated otherwise, the ground norm is \(\ell_1\). For a radius \(\epsilon\ge0\), define the ambiguity set centered at \(Q\) as \(\mathcal{U}_\epsilon(Q)=\left\{P\in\mathcal{M}(u):d_W(P,Q)\le\epsilon\right\}\).
The radius \(\epsilon\) controls the extent of distributional protection. Increasing \(\epsilon\) enlarges the ambiguity set, whereas \(\epsilon=0\) gives \(\mathcal U_0(Q)=\{Q\}\) and recovers the nominal model.

For a product line \(\mathcal S\), define its worst-case expected revenue over this ambiguity set as \(R_\epsilon(\mathcal S;Q)=\inf_{P\in\mathcal U_\epsilon(Q)}R(\mathcal S;P)\).
The corresponding Wasserstein DRO problem is \(\max_{\mathcal S\in\mathcal L}R_\epsilon(\mathcal S;Q)\).

Direct evaluation of the inner problem requires optimization over probability distributions. For a fixed product line and \(\epsilon>0\), this infinite-dimensional problem admits the following one-dimensional dual reformulation.
\begin{proposition}[Dual reformulation of the worst-case revenue]
\label{prop:dual_worst_case_revenue}
Fix a product line \(\mathcal{S}\in\mathcal{L}\), and suppose that the Wasserstein distance is induced by the \(\ell_1\)-norm. If the nominal distribution \(Q\) has a finite first moment, then, for any \(\epsilon>0\), the worst-case expected revenue satisfies
\begin{equation}
  \label{e4}
  \begin{aligned}
  R_\epsilon(\mathcal S;Q)
  &=
  \max_{\lambda\ge0}
  \left\{
  -\lambda\epsilon
  +\lambda\mathbb E_{U\sim Q}
  \left[
  \max_{i\in\mathcal S\cup\{0\}}U_i
  \right]
  -\lambda\mathbb E_{U\sim Q}
  \left[
  \max_{i\in\mathcal S\cup\{0\}}
  \left\{U_i-\frac{r_i}{\lambda}\right\}
  \right]
  \right\}
  \\
  &=
  \max_{\lambda\ge0}
  \left\{
  -\lambda\epsilon
  +\mathbb E_{U\sim Q}
  \left[
  \min_{i\in\mathcal S\cup\{0\}}
  \left\{
  r_i+\lambda
  \left(
  \max_{j\in\mathcal S\cup\{0\}}U_j-U_i
  \right)
  \right\}
  \right]
  \right\}.
  \end{aligned}
  \end{equation}
\end{proposition}

Proposition~\ref{prop:dual_worst_case_revenue} replaces the search over distributions with an optimization over the scalar dual variable \(\lambda\). Both expectations are taken under \(Q\), so the least favorable distribution need not be constructed explicitly. The reformulation can also be combined with the outer problem, yielding a joint optimization over \((\mathcal S,\lambda)\). %
The DRO problem for PLD is
\begin{equation}
    \max_{\mathcal{S}\in\mathcal{L}} R_\epsilon(\mathcal{S};Q)
    \label{dro}
\end{equation}

The dual reveals that the first expectation in~\eqref{e4} is the utility customers obtain from their chosen alternatives, which is the standard buyers' welfare measure in product-line design \citep{Gupta1990Designing,KohliHeuristics1990,liu2023from}. We denote this expectation by \(W(\mathcal S;Q)\). Because \(\lambda\ge0\), welfare enters the dual objective through the nonnegative term \(\lambda W(\mathcal S;Q)\), which creates an endogenous incentive to offer product lines that deliver greater customer value.

Moreover, observe that, 
the final term in~\eqref{e4} satisfies
\[
\lambda\max_{i\in\mathcal{S}\cup\{0\}}
\left\{U_i-{r_i}/{\lambda}\right\}
=
\max_{i\in\mathcal{S}\cup\{0\}}
\left\{\lambda U_i-r_i\right\}.
\]
Therefore, the right-hand side is piecewise linear and convex in \(\lambda\), as it is the pointwise maximum of affine functions. Since expectations preserve convexity,
$
\lambda\mathbb{E}_{\boldsymbol{U}\sim Q}
\left[
\max_{i\in\mathcal{S}\cup\{0\}}
\left\{U_i-{r_i}/{\lambda}\right\}
\right]
$
is also convex in \(\lambda\). It follows that the objective in~\eqref{e4} is concave in \(\lambda\).
This concavity further reveals the welfare incentive induced by the DRO model. Because \(r_0=0\) and \(r_i>0\) for every \(i\in\mathcal S\), the outside option becomes the maximizer of \(\lambda U_i-r_i\) as \(\lambda\to0^+\). Consequently, the right derivative of the dual objective at \(\lambda=0\) is $W(\mathcal S;Q)-\mathbb E_Q[U_0]-\epsilon.$
By concavity, \(\lambda=0\) is optimal if and only if
$
W(\mathcal S;Q)-\mathbb E_Q[U_0]\le\epsilon.
$
When \(\lambda^*=0\), the robust revenue equals zero, the smallest attainable revenue level. Hence, such a product line cannot be robust optimal whenever another feasible product line attains strictly positive robust revenue. This observation imposes a welfare requirement on any robust-optimal product line.

\begin{corollary}[Lower Bound on Buyers' Welfare]
    \label{corollary:welfare_lb}
    Let \(W^*(Q)=\max_{\mathcal S\in\mathcal L}W(\mathcal S;Q)\) be the maximum buyers' welfare over all feasible product lines under \(Q\), and define the critical radius as $\epsilon^{\dagger}(Q)=W^*(Q)-\mathbb E_Q[U_0]$. Then, for any \(0\le \epsilon < \epsilon^{\dagger}(Q)\), every optimal solution \(\mathcal S^*_{\epsilon}(Q)\) satisfies
    \[
    W(\mathcal S^*_{\epsilon}(Q);Q)-\mathbb E_Q[U_0]\ge \epsilon.
    \]
\end{corollary}

In the corollary above, \(W(\mathcal S;Q)-\mathbb E_Q[U_0]\) measures the incremental welfare generated by offering \(\mathcal S\) relative to offering only the outside option. Corollary~\ref{corollary:welfare_lb} therefore shows that robustness imposes an implicit welfare requirement, and that this requirement becomes more stringent as \(\epsilon\) increases. 

Without further structure on the nominal distribution, the expectations in the general reformulation may not be directly computable, and the form of the welfare requirement depends on both the nominal and the feasible family. The next subsection specializes the nominal to the multinomial logit, for which these expectations are available in closed form. 

\subsection{DRO Model with MNL Nominal Distribution}\label{sec:DRO_MNL}
We first specialize the DRO model to an MNL nominal distribution. Let \(Q^{\mathrm{MNL}}\) denote the distribution induced by utilities \(U_i=\mu_i+\varepsilon_i\), \(i\in\mathcal{A}\cup\{0\}\), where the shocks \(\varepsilon_i\) are i.i.d. standard Gumbel random variables. Writing \(v_i=\exp(\mu_i)\), the two expectations in Proposition~\ref{prop:dual_worst_case_revenue} are \(\theta+\log\sum_{j\in\mathcal{S}\cup\{0\}}v_j\) and \(\theta+\log\sum_{i\in\mathcal{S}\cup\{0\}}v_i e^{-r_i/\lambda}\), respectively \citep{sumida2021revenue}, where \(\theta\) is the Euler--Mascheroni constant. The robust model is therefore fully determined by the nominal attraction vector $\boldsymbol{v}=(v_0,v_1,\ldots,v_A)$ and radius $\epsilon>0$. The robust MNL revenue of $\mathcal{S}\in\mathcal{L}$ is given by
\begin{equation}
\label{eq:rho-v}
R_{\epsilon}(\mathcal{S};Q^{\mathrm{MNL}})
  =\sup_{\lambda\ge0}
  \lambda\left[
    -\epsilon+
    \log\!\left(
      \frac{v_0+\sum_{j\in\mathcal{S}}v_j}
           {v_0+\sum_{i\in\mathcal{S}}v_i e^{-r_i/\lambda}}
    \right)
  \right],
\end{equation}
with the value at $\lambda=0$ defined by continuity. Therefore, the product-line design problem is \(\text{[DRO-MNL]}\ \max_{\mathcal{S}\in\mathcal{L}}R_{\epsilon}(\mathcal{S};Q^{\mathrm{MNL}})\).

Without a cardinality constraint, the next proposition establishes that [DRO-MNL] preserves the revenue-ordered structure of its nominal counterpart.
\begin{proposition}[Revenue-Ordered Structure under MNL]\label{prop:revenue_ordered_mnl}
Suppose the nominal utility distribution follows the MNL model, $\epsilon>0$, and there is no cardinality constraint on the offered product line. Relabel the products such that \(r_1\ge r_2\ge \cdots \ge r_A>0\). Then there exists an optimal product line of the form \(\mathcal S_k=\{1,\ldots,k\}\)
for some \(k\in\{1,\ldots,A\}\). Consequently, the robust product line design problem reduces to
\[
\max_{k=1,\ldots,A}
\max_{\lambda\ge 0}
\lambda
\left[
-\epsilon
+
\log
\frac{
v_0+\sum_{i=1}^k v_i
}{
v_0+\sum_{i=1}^k v_i e^{-r_i/\lambda}
}
\right].
\]
\end{proposition}

We next characterize how the optimal prefix changes with the robustness radius. The following proposition implies that the optimal product line expands as the robustness radius increases, thereby raising buyers' welfare.
\begin{proposition}[Expansion of the Optimal Product Line]
    \label{prop:monotone_assortment_epsilon}
    Suppose the nominal utility distribution follows the MNL model and there is no cardinality constraint. Relabel the products such that \(r_1\ge r_2\ge \cdots \ge r_A>0\). For each robustness radius \(\epsilon>0\), let \((\mathcal S^*(\epsilon),\lambda_\epsilon^*)\) be an optimal pair for [DRO-MNL], where \(\mathcal S^*(\epsilon)=\{1,\ldots,k^*(\epsilon)\}\), with ties among optimal prefixes resolved by selecting the largest prefix. Then \(k^*(\epsilon)\) is weakly increasing in \(\epsilon\). Equivalently, if \(0<\epsilon_1\le \epsilon_2\), then
    \[
    |\mathcal S^*(\epsilon_1)|
    \le
    |\mathcal S^*(\epsilon_2)|.
    \]
\end{proposition}

Because buyers' welfare under the MNL model is monotone with respect to set inclusion, the product-line expansion established in Proposition~\ref{prop:monotone_assortment_epsilon} directly translates into a welfare improvement. Thus, a greater degree of robustness induces the firm to offer a product line delivering weakly greater buyers' welfare.

\subsubsection{Data-driven Implementation with Finite Consumer Choice Data.}
When the parameters of the nominal MNL model are not known {\it a priori}, we should first estimate them using the observed choice data.
In this case, we consider a case where the nominal MNL distribution using $N$ choices observed when the offered product line is the candidate set $\mathcal{A}$. 
Specifically, the observed choice data $Y_1,\ldots,Y_N$ are i.i.d. samples denoting the customer choice when offering the full assortment, where $Y_t=0$ denotes no purchase which is assumed to be also observed. Under the MNL assumption and observation of the full product line, the maximum-likelihood estimator of the attraction parameters admits the closed form in the following proposition.

\begin{proposition}[MNL maximum-likelihood estimator]
\label{prop:mle}
Denote the choice counts and empirical shares as
${\hat n}_i := \sum_{t=1}^N \mathbb{I}\{Y_t = i\}$, $\hat{p}_i :={{\hat n}_i}/{N}$, $i = 0, 1, \ldots, A$,
and write \(\hat{\boldsymbol{p}}=(\hat{p}_0,\hat{p}_1,\ldots,\hat{p}_A)\).
If ${\hat n}_0>0$, the extended MLE of $v_i$ %
is
\begin{equation}
\label{eq:mle}
  \vhat_i={{\hat n}_i}/{{\hat n}_0}
  ={\hat{p}_i}/{\hat{p}_0},
  \qquad i=0, 1,\ldots,A.
\end{equation}
It is an interior finite MLE if and only if every ${\hat n}_i>0$, including ${\hat n}_0>0$.
\end{proposition}

We estimate the choice probability of alternative $i$ on a product line $\mathcal{S}$ by
\begin{equation}
\label{eq:assortment-choice-probs}
\hat q_i(\mathcal{S})
:=
\frac{\vhat_i}{1+\sum_{j\in\mathcal{S}}\vhat_j}
=
\frac{\hat{p}_i}{\hat{p}_0+\sum_{j\in\mathcal{S}}\hat{p}_j},
\end{equation}
where the second equality uses $\vhat_i=\hat{p}_i/\hat{p}_0$ for $i\ge 1$.

\begin{lemma}[Uniform choice-probability error]
\label{lem:assortment-l1}
Suppose the true model is MNL with attractions $\boldsymbol{v}=(v_0,v_1,\ldots,v_A)$, normalized by $v_0=1$, and write $Q$ for the induced distribution. Denote the choice probabilities on $\mathcal{A}$ by
\begin{equation}
\label{eq:full-choice-probs}
  p_0=Q(\mathcal{A};Q)=\frac{1}{1+\sum_{j=1}^A v_j},
  \qquad
  p_i=q_i(\mathcal{A};Q)=\frac{v_i}{1+\sum_{j=1}^A v_j},
  \quad i=1,\ldots,A,
\end{equation}
and %
\[
q_i(\mathcal{S};Q)
=
\frac{p_i}{p_0+\sum_{j\in\mathcal{S}}p_j}.
\]
If $p_0\ge\kappa>0$, and $\hat{p}_0>0$, we have
\begin{equation}
\label{eq:assortment-l1-expectation}
\mathbb{E}\left(\sup_{\mathcal{S}\subseteq\mathcal{A}}
\sum_{i\in\mathcal{S}\cup\{0\}}
\bigl|\hat q_i(\mathcal{S})-q_i(\mathcal{S};Q)\bigr|\right) 
\le
\frac{2}{\kappa}\sqrt{\frac{A}{N}}.
\end{equation}
\end{lemma}

The estimated nominal yields a guarantee for the robust objective under the true distribution.
Let $$\widehat{\mathcal{S}}_{N}\in\arg\max_{\mathcal{S}\in\mathcal{L}}\hat R_{\epsilon}(\mathcal{S})$$ denote an optimal product line of [DRO-MNL] under the estimated DRO objective function. Then, we have the following result.

\begin{proposition}[Expected robust revenue loss]
\label{prop:mnl_nominal_revenue}
Under the assumptions of Lemma~\ref{lem:assortment-l1},
suppose $\hat{p}_0>0$.
Then, the expected loss of from the data-driven approach is 
\begin{equation}
\label{eq:mnl_expected_loss}
\begin{aligned}
\mathbb{E}\!\left[
\max_{\mathcal{S}\in\mathcal{L}}
R_\epsilon(\mathcal{S};Q)-
R_\epsilon(\widehat{\mathcal{S}}_{N};Q)
\right]
&\le
\frac{4\bar r}{\kappa}\sqrt{\frac{A}{N}}
+\bar r\exp\!\left(-\frac{N\kappa}{8}\right) =O(N^{-1/2}).
\end{aligned}
\end{equation}
\end{proposition}

The loss bound in Lemma \ref{prop:mnl_nominal_revenue} 
immediately gives the sample size needed to achieve certain prescribed tolerance level. 

\begin{corollary}[Sample size for a 
prescribed expected robust loss]
\label{cor:mnl_sample_size}
Fix $\epsilon\ge 0$ and $\eta>0$. Under the assumptions of Proposition~\ref{prop:mnl_nominal_revenue}, if
\begin{equation}
\label{eq:mnl_sample_size}
N
\ge
\max\left\{
\frac{64A\bar r^2}{\kappa^2\eta^2},\;
\frac{8}{\kappa}\max\!\left\{0,\log\frac{2\bar r}{\eta}\right\}
\right\},
\end{equation}
then
\[
\mathbb{E}\!\left[
\max_{\mathcal{S}\in\mathcal{L}}
R_\epsilon(\mathcal{S};Q)-R_\epsilon(\widehat{\mathcal{S}}_{N};Q)
\right]
\le
\eta.
\]
\end{corollary}

\section{DRO with Empirical Approximation 
}\label{sec:DRO_FC}

The MNL specification yields sharp structural results but imposes restrictive substitution model, and the preceding assortment-expansion result relies on an unconstrained product line. In many applications, the type of choice model may also be unknown.
We therefore develop a nonparametric data-driven formulation based on an empirical nominal distribution, which is also known as the first-choice model \citep{bertsimas2016robust,liu2023from}.

To construct a sample average for the DRO,  we define for convenience,  for $\boldsymbol{U}=({U}_{0},{U}_{1},\ldots,{U}_{A})$,
\begin{align}
\label{eq:phi}
  \phi_{ \mathcal{S},\lambda}(\boldsymbol{U})
  =\min_{i\in \mathcal{S}\cup\{0\}}
     \left\{r_i+\lambda\left(
       \max_{j\in  \mathcal{S}\cup\{0\}}U_j-U_i
     \right)\right\}.
\end{align}
Thus, the outside option is included in \emph{both} the maximum and the minimum, with $r_0=0$.  This is essential: the first term represents the
utility of the customer's preferred available alternative, which may be the
outside option, while the second term allows the adversarially induced choice
to be the outside option.  Restricting either index set to $ \mathcal{S}$ alone would
define a different objective and would incorrectly rule out no purchase.
In terms of this function, we have
\begin{eqnarray}
    \label{008}
R_\epsilon(\mathcal S;Q)=\max_{\lambda\ge 0} \Big\{-\lambda \epsilon+
\mathbb{E}_{\boldsymbol{U}\sim Q}\big[\phi_{ \mathcal{S}, \lambda}(\boldsymbol{U})\big]\Big\}.
\end{eqnarray}

Collect $N$ independent full utility samples  $\hat{\boldsymbol{u}}_n=(\hat{u}_{n0},\hat{u}_{n1},\ldots,\hat{u}_{nA})$ that are i.i.d. with distribution $Q$.
These utility-vector samples can be obtained through, for example, 
    repeated choice experiments used to survey the customer preference such as those used in conjoint-analysis settings, individual-level utility vectors can be estimated from respondents' stated preferences and treated as sampled realizations of heterogeneous customer utilities. See e.g., the discussions in \cite{toubia2003fast,evgeniou2005generalized}.

The sample-average approximation (SAA) of the DRO objective function is
\begin{equation}
\label{eq:empirical-objective}
 \hat R_{\epsilon} ( \mathcal{S})
  :=\sup_{\lambda\ge0}
  \left\{-\lambda\varepsilon+
  \frac1N\sum_{n=1}^N\phi_{ \mathcal{S},\lambda}(\hat{\boldsymbol{u}}_n)\right\}.
\end{equation}
The resulting first-choice product-line design problem, denoted by [DRO-FC], is
\begin{equation}\label{DRO-FC}
    \text{[DRO-FC]}\qquad
    \Shat\in\arg\max_{ \mathcal{S}\in\mathcal{L}}\; 
    \hat R_{\epsilon} ( \mathcal{S}).
\end{equation}

Beyond its implementability, the empirical formulation permits a sample-level characterization of worst-case selection.
\begin{proposition}[Worst-case Selection of Customer \(n\)]\label{lemma:2}
Given a fixed product line \(\mathcal{S}\) and dual multiplier \(\lambda\), let \(k_n \in \arg\min_{i \in \mathcal{S} \cup \{0\}} \{ r_i - \lambda \hat{u}_{ni} \}\).
Then, \(k_n\) satisfies the following conditions:
\begin{itemize}
    \item \emph{Dominance exclusion.} There is no \(i \in \mathcal{S}\cup\{0\}\) such that \(r_i < r_{k_n}\) and \(\hat{u}_{ni} \ge \hat{u}_{nk_n}\).

    \item \emph{Upper-margin condition.} For every \(i \in \mathcal{S} \cup \{0\}\) such that \(r_i > r_{k_n}\) and \(\hat{u}_{ni} > \hat{u}_{nk_n}\),
    \[
    \frac{r_i - r_{k_n}}{\hat{u}_{ni} - \hat{u}_{nk_n}} \ge \lambda.
    \]

    \item \emph{Lower-margin condition.} For every \(i \in \mathcal{S}\cup\{0\}\) such that \(r_i < r_{k_n}\) and \(\hat{u}_{ni} < \hat{u}_{nk_n}\),
    \[
    \frac{r_{k_n} - r_i}{\hat{u}_{nk_n} - \hat{u}_{ni}} \le \lambda.
    \]
\end{itemize}
\end{proposition}

The dominance condition states that the worst-case choice \(k_n\) is Pareto-undominated under the adversary's revenue--utility ordering: no available alternative can generate strictly lower revenue while delivering weakly higher utility. 

The two margin conditions provide a consumer-value interpretation of
\(\lambda\). Relative to the worst-case choice \(k_n\), a higher-revenue
alternative cannot offer enough extra utility per unit of additional revenue to
justify moving upward; equivalently, its utility gain per unit of revenue
premium is at most \(1/\lambda\). Conversely, a lower-revenue alternative must
sacrifice sufficiently much utility per unit of revenue saved, so moving
downward is not attractive from the consumer's perspective. Thus, \(1/\lambda\)
acts as a utility-per-revenue benchmark for comparing products around
\(k_n\).

If \(k_n\neq0\), dominance exclusion implies \(\hat u_{nk_n}>\hat u_{n0}\). Applying the lower-margin condition to the outside option then gives \(r_{k_n}/(\hat u_{nk_n}-\hat u_{n0}) \le \lambda\). Summing over customers, and noting that both sides contribute zero when \(k_n=0\), gives
\[
\sum_{n\in \mathcal{N}} r_{k_n} \le \lambda \sum_{n\in \mathcal{N}} (\hat u_{nk_n}-\hat u_{n0}).
\]
The left-hand side is the aggregate revenue generated by the worst-case choices, and the right-hand sum is their aggregate utility gain relative to the outside option. Thus, the inequality indicates that there must exists an alternative for each sample such that these alternatives must deliver at least \(1/\lambda\) units of incremental utility per unit of revenue on average. 

A smaller \(\lambda\) tightens this requirement and forces the worst-case demand realization to place weight only on products with greater consumer value per unit of revenue. The next result links this interpretation to the robustness radius.

\begin{proposition}[Monotonicity of \(\lambda^*_\epsilon\)]\label{pro:mon_lam}
For [DRO-FC], let \((\mathcal S_j,\lambda_j)\) be any joint optimizer at robustness level \(\epsilon_j\), \(j=1,2\). If \(0<\epsilon_1<\epsilon_2\), then \(\lambda_1\ge\lambda_2\).
\end{proposition}

As \(\epsilon\) increases, Proposition~\ref{pro:mon_lam} implies that the upper bound on revenue per unit of incremental utility becomes tighter for products selected under the worst-case demand realization. Accordingly, under a greater degree of uncertainty aversion, the firm should offer products that deliver higher consumer value per unit of revenue.

There are two questions that need to be addressed: 1) How good is the data-driven solution of (\ref{DRO-FC}) compared with that of the origin PLD problem? and 2) how to compute the optimal solution for (\ref{DRO-FC})? In the rest of this section we address the first question, and in the next section, we discuss the second.

\subsection{Statistical Guarantees of the Empirical DRO Model}
\label{subsubsec:overall_empirical_guarantee}

We quantify the statistical loss from replacing the population
DRO problem with its sample-average approximation. Let $\Shat$
be an exact optimizer of the empirical problem~\eqref{DRO-FC}.
We evaluate this product line under the population distribution
$Q$ and compare it with the optimum under the same distribution:
\[
L_N
=
\max_{\mathcal S\in\mathcal L}R_\epsilon(\mathcal S;Q)
-
R_\epsilon(\Shat;Q).
\]

The finite-sample analysis uses the following bound on the dual integrand, where $\bar r=\max_{i\in\mathcal{A}}r_i$ is the largest product revenue.

\begin{lemma}
\label{lem:envelope}
For every product line $\mathcal{S}$, every $\lambda\ge0$, and every utility vector $\boldsymbol{U}$,
\[
0\le \phi_{\mathcal{S},\lambda}(\boldsymbol{U})\le\bar r.
\]
Moreover, the supremum over $\lambda$ in both the population objective and the empirical objective can be restricted to $[0,\bar r/\epsilon]$.
\end{lemma}

For $\mathcal{S}\in\mathcal{L}$ and $\lambda\ge0$, write $\mu_{\mathcal{S}}(\lambda)=\mathbb{E}_{Q}[\phi_{\mathcal{S},\lambda}(\boldsymbol{U})]$ and let $\hat\mu_{\mathcal{S}}(\lambda)$ denote the corresponding sample average. Define the uniform dual deviation
\begin{equation}
\label{eq:ZN}
  Z_N:=\sup_{\mathcal{S}\in \mathcal{L}}
  \sup_{\lambda\ge0}
  \bigl|\hat\mu_{\mathcal{S}}(\lambda)-\mu_{\mathcal{S}}(\lambda)\bigr|.
\end{equation}
Lemma~\ref{lem:envelope} implies $Z_N\le\bar r$. For every $\epsilon>0$, the penalty $-\lambda\epsilon$ cancels in the comparison of the two dual objectives, so
\[
\sup_{\mathcal{S}\in\mathcal{L}}
\bigl|R_{\epsilon}(\mathcal{S};Q)-\hat R_{\epsilon}(\mathcal{S})\bigr|
\le Z_N.
\]
For $\alpha\in(0,1)$, define
\begin{equation}
\label{eq:b_N}
b_N(\alpha)
=
\bar r
\min\left\{
1,\,
\frac{\sqrt{\pi}+\sqrt{\log(|\mathcal{L}|/\alpha)}}{\sqrt{2N}}
\right\}.
\end{equation}

\begin{theorem}[Finite-sample accuracy]
\label{thm:loss}
Suppose the $N$ utility vectors are drawn independently from $Q$, and let $\alpha\in(0,1)$. With probability at least $1-\alpha$,
\begin{equation}
\label{eq:uniform_statistical_error}
Z_N\le b_N(\alpha).
\end{equation}
The same event yields, for every $\epsilon>0$,
\begin{align}
\bigl|V^{\star}_{N,\epsilon}-V_\epsilon^\star(Q)\bigr|
&\le b_N(\alpha), \label{eq:value_statistical_error}\\
0\le L_N
&\le 2b_N(\alpha), \label{eq:decision_statistical_error}
\end{align}
where $V_\epsilon^\star(Q)=\max_{\mathcal{S}\in\mathcal{L}}R_\epsilon(\mathcal{S};Q)$ and $V^{\star}_{N,\epsilon}=\max_{\mathcal{S}\in\mathcal{L}}\hat R_\epsilon(\mathcal{S})$. The expected loss satisfies
\begin{equation}
\label{eq:expected-loss}
\mathbb{E}[L_N]
\le
\frac{2\bar r}{\sqrt{2N}}
\left(
\sqrt{\pi}+\frac12+\sqrt{1+\log|\mathcal{L}|}
\right).
\end{equation}
Since $|\mathcal{L}|\le 2^{A}$, both bounds remain valid when $|\mathcal{L}|$ is replaced by $2^{A}$.
\end{theorem}

The bound in~\eqref{eq:expected-loss} gives a sufficient sample
size for a target expected loss due to sampling. The next section
accounts for the additional optimization error when the empirical
problem is solved to a prescribed tolerance.

\section{Computational Algorithm} %
\label{sec:delta_optimal}

Although the empirical approximation error deminishes as $N$ increases, [DRO-FC] with a large sample size is computationally challenging. Even its nominal counterpart, obtained by setting $\epsilon=0$, is NP-hard, and robustness couples the continuous multiplier with the combinatorial product-line decision. McCormick envelopes and big-$M$ bounds yield the exact formulation [DRO-MILP] in Section~\ref{appendix:milp} of the online supplement. Its weak linear programming relaxation, however, makes [DRO-MILP] impractical beyond relatively small instances. We therefore develop a certified approximation framework for [DRO-FC] that returns a $(1-\delta)$-optimal solution. Throughout this section, we consider the cardinality-constrained feasible family $\mathcal{L}=\left\{\mathcal{S}:|\mathcal{S}|\le p\right\}$.
We write $\bar{\mathcal{S}}=\mathcal{S}\cup\{0\}$ for the offered alternatives including the outside option.
For the empirical distribution, we write $\epsdag=\epsilon^\dagger(\widehat Q_N)=\max_{\mathcal S\in\mathcal L}W_N(\mathcal S)-N^{-1}\sum_{n\in\mathcal N}\hat u_{n0}$ for the critical radius in Corollary~\ref{corollary:welfare_lb}.
We maintain $0\le\epsilon<\epsdag$, so [DRO-FC] has a positive optimal value.

To expose the two-level structure of the problem, define the functions
$C(\cdot,\cdot),G(\cdot,\cdot):\mathcal{L}\times\mathbb{R}_+\to\mathbb{R}$ by
\[
C(\mathcal{S},\lambda)
=
\frac{1}{N}
\sum_{n\in\mathcal{N}}
\min_{i\in\bar{\mathcal{S}}}
\left(r_i-\lambda\hat{u}_{ni}\right)
\]
and $G(\mathcal{S},\lambda)=\lambda W_N(\mathcal{S})+C(\mathcal{S},\lambda)$.
Here, $C(\mathcal S,\lambda)$ is the average of the customer-level minima $\min_{i\in\bar{\mathcal S}}(r_i-\lambda\hat u_{ni})$ obtained by separating $\lambda\max_{j\in\bar{\mathcal S}}\hat u_{nj}$ from the minimum in~\eqref{eq:phi}, and $G(\mathcal S,\lambda)$ is the sample average of $\phi_{\mathcal S,\lambda}$ in~\eqref{eq:empirical-objective} before the ambiguity penalty $-\lambda\epsilon$.
For a fixed $\lambda$, let
\begin{equation}
\label{eq:flambda}
f(\lambda)
=
\max_{\mathcal{S}\in\mathcal{L}}
G(\mathcal{S},\lambda)
\end{equation}
denote the optimal value of the product-line subproblem at multiplier $\lambda$. For a fixed product line $\mathcal S$, write $h_{\mathcal S}(\lambda)=-\lambda\epsilon+G(\mathcal S,\lambda)$. By~\eqref{eq:empirical-objective}, $\hat R_{\epsilon}(\mathcal S)=\sup_{\lambda\ge0}h_{\mathcal S}(\lambda)$. Write $R^{\mathrm{FC}}_{N,\epsilon}(\mathcal S):=\hat R_{\epsilon}(\mathcal S)$ and $V^{\star}_{N,\epsilon}:=\max_{\mathcal S\in\mathcal L}\hat R_{\epsilon}(\mathcal S)$. It follows that $h(\lambda)=-\lambda\epsilon+f(\lambda)=\max_{\mathcal S\in\mathcal L}h_{\mathcal S}(\lambda)$.

Before constructing either level of the algorithm, we restrict the multiplier to a finite, data-dependent interval. For customer $n$, consider an ordered pair $(i,j)$ for which product $i$ has both higher revenue and higher utility than product $j$. This pair can affect the right slope of $h$ at $\lambda$ only when its revenue-to-utility ratio is strictly greater than $\lambda$. Accordingly, define
\begin{equation}
\label{eq:D_n}
D_n(\lambda)
=
\max\left\{
\hat u_{ni}-\hat u_{nj}:
i\in\mathcal A,\ j\in\mathcal A\cup\{0\},\
r_i>r_j,\ \hat u_{ni}>\hat u_{nj},\
    \frac{r_i-r_j}{\hat u_{ni}-\hat u_{nj}}>\lambda
\right\}.
\end{equation}
If the set in~\eqref{eq:D_n} is empty, we set $D_n(\lambda)=0$. Averaging these customer-level bounds gives
\begin{equation}
\label{eq:Phi}
\Phi(\lambda)
=
\frac{1}{N}
\sum_{n\in\mathcal{N}}
D_n(\lambda).
\end{equation}
\begin{lemma}[Valid Upper Bound on $\lambda$]
\label{prop:lambda_bounds}
The envelope $\Phi(\lambda)$ is nonincreasing, right-continuous, and piecewise constant in $\lambda$. Moreover, for every $\mathcal S\in\mathcal L$,
\begin{equation}
\label{eq:lambda_bounds}
R^{\mathrm{FC}}_{N,\epsilon}(\mathcal S)
=\max_{0\le\lambda\le\lambda^u}h_{\mathcal S}(\lambda),
\qquad
\lambda^u := \min \left\{ \lambda \ge 0 : \Phi(\lambda) \le \epsilon \right\}.
\end{equation}
Consequently, [DRO-FC] admits an optimal multiplier $\lambda^*\in[0,\lambda^u]$.
\end{lemma}
The only breakpoints of $\Phi$ are the $O(N|\mathcal A|^2)$ ratios in~\eqref{eq:D_n}. Thus, $\lambda^u$ is obtained by sorting these ratios and scanning the value immediately after each breakpoint until $\Phi(\lambda)\le\epsilon$; the implementation uses the next representable floating-point value to perform this right-hand scan.

The [DRO-FC] problem is therefore equivalent to
\begin{equation}
\label{eq:hlambda}
h^\star
=
\max_{\lambda\in[0,\lambda^u]}
h(\lambda),
\end{equation}
where $h^\star=V^{\star}_{N,\epsilon}$. Our method separates the two decisions in~\eqref{eq:hlambda}. Section~\ref{sec:SBMB} first constructs an exact combinatorial solver for $f(\lambda)$ at a fixed multiplier. Section~\ref{sec:APS} then uses this solver as an oracle in a certified search over $[0,\lambda^u]$.

\subsection{Set-Branching with Modular Bounds (SBMB)}\label{sec:SBMB}
Fix $\lambda\in[0,\lambda^u]$. The inner problem~\eqref{eq:flambda} remains difficult even though the multiplier is no longer a decision variable. A direct MILP formulation of this fixed-$\lambda$ problem retains $N(A+1)$ assignment variables, $N$ epigraph variables, and $O(NA)$ big-$M$ and linking constraints. SBMB avoids this overhead by branching directly on product sets and bounding each node in closed form, without solving a continuous relaxation. The fixed-$\lambda$ problem is nevertheless NP-hard: for sufficiently large $\lambda$, it reduces to nominal first-choice product-line design~\citep{KohliHeuristics1990}.

SBMB keeps a feasible incumbent $(\mathcal S_\lambda,G_{\mathrm{best}})$ and a collection $\mathcal Q$ of open nodes. A node is a pair $(\mathcal I,\mathcal F)$ of disjoint subsets of $\mathcal A$. Products in $\mathcal I$ have been included, and those in $\mathcal F$ have not yet been decided. The remaining products are excluded, so $\mathcal E=\mathcal A\setminus(\mathcal I\cup\mathcal F)$ need not be stored. The root is $(\mathcal I_0,\mathcal F_0)=(\emptyset,\mathcal A)$. At any node, the residual capacity is $k=p-|\mathcal I|$, and its feasible descendants are $\mathcal I\cup\mathcal T$, where $\mathcal T\subseteq\mathcal F$ and $|\mathcal T|\le k$.

The outside option is always available, but it is not a member of $\mathcal I$, $\mathcal F$, or $\mathcal E$. In the code, index~0 is retained in the internal included-set array, which therefore has $p+1$ slots. Our notation counts candidate products only and hence has $|\mathcal I|\le p$. For any product line,
\[
G(\mathcal S,\lambda)
=\frac{1}{N}\sum_{n\in\mathcal N}
\min_{i\in\bar{\mathcal S}}
\left\{r_i+\lambda\left(\max_{j\in\bar{\mathcal S}}\hat u_{nj}-\hat u_{ni}\right)\right\}
\ge0.
\]
Thus, the best singleton weakly dominates the empty line, whose value is zero. The code omits $t=0$ whenever $\mathcal I=\emptyset$; we retain it in~\eqref{eq:node_ub} to state one bound for every node.

After a node is bounded, it is stored in $\mathcal Q$ as $(\mathcal I,\mathcal F,U)$, where $U=\mathrm{UB}_{\mathrm{node}}(\mathcal I,\mathcal F)$. The pair $(\mathcal I,\mathcal F)$ identifies the subproblem, and $U$ determines its priority in the queue. The multiplier is fixed during an SBMB call and is therefore suppressed in the notation for $\mathrm{UB}_{\mathrm{node}}$.

The bound rests on two diminishing-returns properties. Adding a product increases $W_N(\cdot)$, with smaller marginal gains for larger product lines. It decreases $C(\cdot,\lambda)$, and the magnitude of this decrease also becomes smaller as the product line grows. Lemma~\ref{lem:modularity} records these properties.
\begin{lemma}[Submodularity and Supermodularity]\label{lem:modularity}
Given a fixed $\lambda \ge 0$, $W_N(\cdot)$ is monotone
nondecreasing and submodular, and $C(\cdot,\lambda)$ is monotone nonincreasing and
supermodular. That is, for every $\mathcal{S}_1 \subseteq \mathcal{S}_2 \subseteq \mathcal{A}$
and $j \notin \mathcal{S}_2$,
\begin{align*}
W_N(\mathcal{S}_1) &\le W_N(\mathcal{S}_2),
& W_N(\mathcal{S}_1 \cup \{j\}) - W_N(\mathcal{S}_1) &\ge W_N(\mathcal{S}_2 \cup \{j\}) - W_N(\mathcal{S}_2),\\
C(\mathcal{S}_1,\lambda) &\ge C(\mathcal{S}_2,\lambda),
& C(\mathcal{S}_1 \cup \{j\},\lambda) - C(\mathcal{S}_1,\lambda) &\le C(\mathcal{S}_2 \cup \{j\},\lambda) - C(\mathcal{S}_2,\lambda).
\end{align*}
\end{lemma}

At a node $(\mathcal I,\mathcal F)$, the included products determine the customer-level values
\begin{align}
    \bar{u}_n(\mathcal I) &= \max_{i \in \mathcal I \cup \{0\}} \hat{u}_{ni},
    &
    \bar{c}_n(\mathcal I,\lambda) &= \min_{i \in \mathcal I \cup \{0\}} (r_i - \lambda \hat{u}_{ni}). \label{eq:node_vals}
\end{align}
Thus, $W_N(\mathcal I)=N^{-1}\sum_n\bar u_n(\mathcal I)$ and $C(\mathcal I,\lambda)=N^{-1}\sum_n\bar c_n(\mathcal I,\lambda)$. For each free product $j\in\mathcal F$, let $m^W_j(\mathcal I)$ be its welfare gain and $m^C_j(\mathcal I,\lambda)$ its decrease in $C(\cdot,\lambda)$ when added alone:
\begin{equation}\label{eq:marginals}
\begin{aligned}
    m^W_j(\mathcal I)
    &=W_N(\mathcal I\cup\{j\})-W_N(\mathcal I)\\
    &= \tfrac{1}{N}\sum_{n \in \mathcal N} [\hat{u}_{nj} - \bar{u}_n(\mathcal I)]_+,
    \\
    m^C_j(\mathcal I,\lambda)
    &=C(\mathcal I,\lambda)-C(\mathcal I\cup\{j\},\lambda)\\
    &= \tfrac{1}{N}\sum_{n \in \mathcal N} [\bar{c}_n(\mathcal I,\lambda) - (r_j - \lambda\hat{u}_{nj})]_+.
\end{aligned}
\end{equation}
All marginals can be computed in $O(N|\mathcal F|)$ operations. The modular component of the node bound uses only the $k$ largest welfare gains and the $k$ smallest decreases in $C(\cdot,\lambda)$. The code maintains these order statistics in $O(k|\mathcal F|)$ time instead of sorting all of $\mathcal F$.

Consider a completion $\mathcal T\subseteq\mathcal F$ with $|\mathcal T|=t$. By submodularity, its total welfare gain is at most the sum of the $t$ largest values of $m^W_j(\mathcal I)$. For the $C(\cdot,\lambda)$ term, any $t$ chosen products must include one whose one-product decrease is at least the $t$-th smallest value of $m^C_j(\mathcal I,\lambda)$. Hence, $C(\mathcal I\cup\mathcal T,\lambda)$ is at most $C(\mathcal I,\lambda)$ minus that order statistic. Proposition~\ref{prop:modular_ub} combines the two bounds.

\begin{proposition}[Node Upper Bound]\label{prop:modular_ub}
At a node $(\mathcal I,\mathcal F)$ with residual capacity $k=p-|\mathcal I|$,
let $m^W_{[1]}\ge m^W_{[2]}\ge\cdots$ denote the values $m^W_j(\mathcal I)$, $j\in\mathcal F$, in decreasing order and
$m^C_{(1)}\le m^C_{(2)}\le\cdots$ the values $m^C_j(\mathcal I,\lambda)$, $j\in\mathcal F$, in increasing order, with
$m^C_{(0)}=0$. For $t=0,\ldots,\min\{k,|\mathcal{F}|\}$, define
\begin{equation}\label{eq:node_ub}
\begin{aligned}
\mathrm{UB}(t\mid\mathcal I,\mathcal F)
&=\lambda\left(W_N(\mathcal{I})+\sum_{q=1}^{t}m^W_{[q]}\right)
+C(\mathcal{I},\lambda)-m^C_{(t)},\\
\mathrm{UB}_{\mathrm{mod}}(\mathcal I,\mathcal F)
&=\max_{0\le t\le\min\{k,|\mathcal{F}|\}}\mathrm{UB}(t\mid\mathcal I,\mathcal F).
\end{aligned}
\end{equation}
For $n\in\mathcal N$ and $j\in\mathcal F$, define the customer-level net gain
\[
\eta_{nj}(\mathcal I,\lambda)
=\lambda[\hat u_{nj}-\bar u_n(\mathcal I)]_+
-[\bar c_n(\mathcal I,\lambda)-(r_j-\lambda\hat u_{nj})]_+,
\]
and let
\begin{equation}\label{eq:node_sep_ub}
\begin{aligned}
\mathrm{UB}_{\mathrm{sep}}(\mathcal I,\mathcal F)
&=G(\mathcal I,\lambda)
+\frac{1}{N}\sum_{n\in\mathcal N}
\left[\max_{j\in\mathcal F}\eta_{nj}(\mathcal I,\lambda)\right]_+,\\
\mathrm{UB}_{\mathrm{node}}(\mathcal I,\mathcal F)
&=\min\{\mathrm{UB}_{\mathrm{mod}}(\mathcal I,\mathcal F),
\mathrm{UB}_{\mathrm{sep}}(\mathcal I,\mathcal F)\}.
\end{aligned}
\end{equation}
The inner maximum in~\eqref{eq:node_sep_ub} is defined as zero when $\mathcal F=\emptyset$.
Then $\max_{\mathcal{T}\subseteq\mathcal{F},\,|\mathcal{T}|\le k}G(\mathcal{I}\cup\mathcal{T},\lambda)\le\mathrm{UB}_{\mathrm{node}}(\mathcal I,\mathcal F)$.
\end{proposition}
The separable component couples the utility improvement and cost decrease at the customer--product level before maximizing over free products. It requires the same $O(N|\mathcal F|)$ pass used to compute the marginals and is particularly effective when the two terms in $G(\cdot,\lambda)$ are individually large but offset one another.

SBMB begins with a greedy incumbent. Starting from the best singleton, it repeatedly adds the product with the largest positive increment in $G(\cdot,\lambda)$ until no positive increment remains or the line contains $p$ products. The root triple $(\mathcal I_0,\mathcal F_0,U_0)$ is then inserted into $\mathcal Q$.

The search processes open nodes in decreasing order of $U$. Suppose $(\mathcal I,\mathcal F,U)$ has the largest bound in $\mathcal Q$. If $U\le G_{\mathrm{best}}$, neither this node nor any other open node can improve the incumbent, so the algorithm terminates. Otherwise, $\mathcal I$ is itself a feasible product line. It replaces the incumbent if $\mathcal I\neq\emptyset$ and $G(\mathcal I,\lambda)>G_{\mathrm{best}}$.

If $|\mathcal I|<p$ and $\mathcal F\neq\emptyset$, the node is not a leaf. SBMB chooses $j^*\in\arg\max_{j\in\mathcal F}m^W_j(\mathcal I)$ and forms two children:
\[
(\mathcal I^{\mathrm{in}},\mathcal F^{\mathrm{in}})
=(\mathcal I\cup\{j^*\},\mathcal F\setminus\{j^*\}),
\qquad
(\mathcal I^{\mathrm{out}},\mathcal F^{\mathrm{out}})
=(\mathcal I,\mathcal F\setminus\{j^*\}).
\]
The include child fixes $j^*$ in the line, while the exclude child removes it from further consideration. For each child $b\in\{\mathrm{in},\mathrm{out}\}$, SBMB recomputes $U^b=\mathrm{UB}_{\mathrm{node}}(\mathcal I^b,\mathcal F^b)$ and inserts $(\mathcal I^b,\mathcal F^b,U^b)$ into $\mathcal Q$ only if $U^b>G_{\mathrm{best}}$. A leaf has no children. The greedy start and the choice of branching product affect runtime, but not the validity of the bounds.

\begin{algorithm}[H]
\caption{Set-Branching with Modular Bounds (SBMB)}\label{alg:SBMB}
{\small
\begin{algorithmic}[1]
\Require Instance parameters $(\boldsymbol{r}, \widehat{\boldsymbol{u}}, N, \mathcal A, p)$ with $1\le p\le|\mathcal A|$ and fixed $\lambda\ge0$
\Ensure An optimal product line $\mathcal S_\lambda$ and $f(\lambda)=G(\mathcal S_\lambda,\lambda)$
\State Construct the greedy feasible line described above; set $(\mathcal S_\lambda,G_{\mathrm{best}})$ to this line and its value
\State Set $(\mathcal I_0,\mathcal F_0)\gets(\emptyset,\mathcal A)$ and $U_0\gets\mathrm{UB}_{\mathrm{node}}(\mathcal I_0,\mathcal F_0)$; initialize $\mathcal Q\gets\{(\mathcal I_0,\mathcal F_0,U_0)\}$
\While{$\mathcal{Q} \neq \emptyset$}
    \State Select $(\mathcal I,\mathcal F,U)\in\mathcal Q$ with maximal $U$; set $\mathcal Q\gets\mathcal Q\setminus\{(\mathcal I,\mathcal F,U)\}$
    \If{$U\le G_{\mathrm{best}}$}
        \State \Return $(\mathcal S_\lambda,G_{\mathrm{best}})$
    \EndIf
    \State If $\mathcal I\neq\emptyset$ and $G(\mathcal I,\lambda)>G_{\mathrm{best}}$, set $(\mathcal S_\lambda,G_{\mathrm{best}})\gets(\mathcal I,G(\mathcal I,\lambda))$
    \If{$|\mathcal I|<p$ and $\mathcal F\neq\emptyset$} \Comment{a nonleaf node}
        \State Compute $m^W_j(\mathcal I),m^C_j(\mathcal I,\lambda)$ for $j\in\mathcal F$ via~\eqref{eq:marginals}; select $j^*\gets\arg\max_{j\in\mathcal F}m^W_j(\mathcal I)$
        \State Set the include child $(\mathcal I^{\mathrm{in}},\mathcal F^{\mathrm{in}})\gets(\mathcal I\cup\{j^*\},\mathcal F\setminus\{j^*\})$
        \State $U^{\mathrm{in}}\gets\mathrm{UB}_{\mathrm{node}}(\mathcal I^{\mathrm{in}},\mathcal F^{\mathrm{in}})$
        \State \textbf{if} $U^{\mathrm{in}}>G_{\mathrm{best}}$ \textbf{then} $\mathcal Q\gets\mathcal Q\cup\{(\mathcal I^{\mathrm{in}},\mathcal F^{\mathrm{in}},U^{\mathrm{in}})\}$
        \State Set the exclude child $(\mathcal I^{\mathrm{out}},\mathcal F^{\mathrm{out}})\gets(\mathcal I,\mathcal F\setminus\{j^*\})$
        \State $U^{\mathrm{out}}\gets\mathrm{UB}_{\mathrm{node}}(\mathcal I^{\mathrm{out}},\mathcal F^{\mathrm{out}})$
        \State \textbf{if} $U^{\mathrm{out}}>G_{\mathrm{best}}$ \textbf{then} $\mathcal Q\gets\mathcal Q\cup\{(\mathcal I^{\mathrm{out}},\mathcal F^{\mathrm{out}},U^{\mathrm{out}})\}$
    \EndIf
\EndWhile
\State \Return $(\mathcal S_\lambda,G_{\mathrm{best}})$
\end{algorithmic}}
\end{algorithm}
Proposition~\ref{prop:modular_ub} makes every pruning step valid, and each branch removes one product from $\mathcal F$, so the search tree is finite. When SBMB terminates, every unexplored node has an upper bound no larger than $G_{\mathrm{best}}$; Algorithm~\ref{alg:SBMB} therefore returns the exact value $f(\lambda)$.

\subsection{Adaptive Parametric Search (APS)}
\label{sec:APS}

With the exact fixed-$\lambda$ oracle in place, it remains to maximize $h$ over $[0,\lambda^u]$. For a fixed product line $\mathcal S$, $G(\mathcal S,\cdot)$ is concave and piecewise linear. Its pointwise maximum over $\mathcal L$, however, need not be concave or convex. APS therefore does not rely on the global shape of $h$. Instead, it combines exact endpoint values with the uniform slope bounds in Lemma~\ref{lem:bilipschitz}.

\begin{lemma}[Slope Bounds]
\label{lem:bilipschitz}\label{prop:interval_ub}
The function $h(\lambda) = -\lambda\epsilon + f(\lambda)$ satisfies, for all
$0 \le \lambda_1 < \lambda_2 \le \lambda^u$,
\begin{equation}\label{eq:bilipschitz}
    -L^- \;\le\; \frac{h(\lambda_2) - h(\lambda_1)}{\lambda_2 - \lambda_1} \;\le\; L^+,
\end{equation}
where $L^- = \epsilon$, $L^+ = N^{-1} \sum_{n \in \mathcal{N}} \Delta_n - \epsilon$, and $\Delta_n = \max_{i \in \mathcal{A}\cup\{0\}} \hat{u}_{ni} - \hat{u}_{n0}$.
\begin{itemize}
  \item If $L^+>0$, fix any interval $[a,b]\subseteq[0,\lambda^u]$ and define $x_{ab}=\bigl[h(b)-h(a)+L^-(b-a)\bigr]/(L^++L^-)$.
  Then $x_{ab}\in[0,b-a]$ and
  \begin{equation}\label{eq:ub_star}
  \begin{aligned}
      \max_{\lambda \in [a,b]} h(\lambda)
      &\;\le\; \mathrm{UB}^*([a,b])
      :=h(a)+L^+x_{ab}\\
      &=h(b)+L^-(b-a-x_{ab})
      =\frac{L^-h(a)+L^+h(b)+L^+L^-(b-a)}{L^++L^-}.
  \end{aligned}
  \end{equation}
  \item If $L^+\le0$, then $h(\cdot)$ is nonincreasing and $\lambda=0$ is an optimal multiplier.
\end{itemize}

\end{lemma}

The maintained condition $\epsilon<\epsdag$ implies $L^+>0$, so the first case of Lemma~\ref{lem:bilipschitz} applies throughout Algorithm~\ref{alg:APS}; the second case is included for completeness. On an interval $[a,b]$, the slope bounds give $h(\lambda)\le h(a)+L^+(\lambda-a)$ and $h(\lambda)\le h(b)+L^-(b-\lambda)$.
The pointwise minimum of these affine functions is a tent-shaped upper bound on $h$, whose maximum occurs at their intersection $a+x_{ab}$.

APS initially calls SBMB at $0$ and $\lambda^u$ and sets $\IAPS=\{[0,\lambda^u]\}$. Every interval in $\IAPS$ stores its exact endpoint values and carries the bound in~\eqref{eq:ub_star}. Let $\mathrm{ub}_{\max}=\max_{[a,b]\in\IAPS}\mathrm{UB}^*([a,b])$ denote the largest active bound. The incumbent $(\widehat{\mathcal S},\widehat\lambda)$ provides the lower bound $h_{\mathrm{best}}=h_{\widehat{\mathcal S}}(\widehat\lambda)$. An interval is pruned only when its upper bound does not exceed the incumbent; hence, whenever $\IAPS\neq\emptyset$, $h_{\mathrm{best}}\le h^\star\le\mathrm{ub}_{\max}$. The certificate gap is therefore $\Gamma:=\mathrm{ub}_{\max}-h_{\mathrm{best}}\ge0$,
and APS stops once $\Gamma\le\delta\,\mathrm{ub}_{\max}$. If the certificate is not yet satisfied, APS selects
$[a,b]\in\arg\max_{[c,d]\in\IAPS}\mathrm{UB}^*([c,d])$ and calls SBMB at $\lambda_{\mathrm{new}}=(a+b)/2$. The exact value $h(\lambda_{\mathrm{new}})$ becomes the common endpoint of the two child intervals, while the returned product line is processed through the solution pool described below. APS then replaces $[a,b]$ by its two children and prunes every interval whose upper bound does not exceed $h_{\mathrm{best}}$:
\[
\IAPS\gets
\bigl(\IAPS\setminus\{[a,b]\}\bigr)
\cup\{[a,\lambda_{\mathrm{new}}],[\lambda_{\mathrm{new}},b]\},
\qquad
\IAPS\gets
\{[c,d]\in\IAPS:\mathrm{UB}^*([c,d])>h_{\mathrm{best}}\}.
\]
Both children are created and assessed independently; APS does not impose unimodality on $h$.

The endpoint values support the certificate; the product lines returned by SBMB strengthen the incumbent. Once a product line is fixed, its multiplier can be optimized exactly without another combinatorial search. Although this fixed-line problem can be written as a linear program, APS does not call an LP solver. Proposition~\ref{lemma:2} identifies every value of $\lambda$ at which a customer's worst-case selection can change. APS therefore enumerates these breakpoints and evaluates the objective directly. For $\mathcal S\in\mathcal L$, define
\begin{equation}\label{eq:pool_bp}
\begin{aligned}
\mathcal B_{\mathcal S}
=\{0,\lambda^u\}\ \cup\ 
\left\{
\frac{r_i-r_j}{\hat u_{ni}-\hat u_{nj}}:
\begin{array}{l}
n\in\mathcal N,\ \{i,j\}\subseteq\bar{\mathcal S},\ \hat u_{ni}\ne\hat u_{nj},\\[-1mm]
0\le (r_i-r_j)/(\hat u_{ni}-\hat u_{nj})\le\lambda^u
\end{array}
\right\},\\
\text{[Pool-BP($\mathcal S$)]}\qquad
H_{\mathcal S}
=\max_{\lambda\in\mathcal B_{\mathcal S}}h_{\mathcal S}(\lambda),
\qquad
\lambda_{\mathcal S}\in\arg\max_{\lambda\in\mathcal B_{\mathcal S}}h_{\mathcal S}(\lambda).
\end{aligned}
\end{equation}
For fixed $\mathcal S$, the worst-case alternative of customer $n$ can change only when two terms $r_i-\lambda\hat u_{ni}$ intersect. Because $h_{\mathcal S}(\cdot)$ is concave and piecewise linear, an optimum on $[0,\lambda^u]$ is attained at one of these intersections or at an endpoint. Lemma~\ref{prop:lambda_bounds} therefore gives $H_{\mathcal S}=h_{\mathcal S}(\lambda_{\mathcal S})=R^{\mathrm{FC}}_{N,\epsilon}(\mathcal S)$.
Thus, [Pool-BP($\mathcal S$)] is an exact finite enumeration, not a solver-based subproblem. The candidate set has at most $2+N\binom{|\bar{\mathcal S}|}{2}$ elements. The implementation evaluates these candidates in vectorized, bounded-memory batches and caches the exact pair $(\lambda_{\mathcal S},H_{\mathcal S})$ when a product line first enters the solution pool $\mathcal P$. The incumbent is $\widehat{\mathcal S}\in\arg\max_{\mathcal S\in\mathcal P}H_{\mathcal S}$, with $\widehat\lambda=\lambda_{\widehat{\mathcal S}}$ and $h_{\mathrm{best}}=H_{\widehat{\mathcal S}}$.
Algorithm~\ref{alg:APS} summarizes the complete outer procedure. A call to SBMB solves~\eqref{eq:flambda} exactly and returns an optimal product line at the evaluated multiplier; [Pool-BP($\mathcal S$)] then optimizes that line exactly over the multiplier.

\begin{algorithm}[H]
\caption{Adaptive Parametric Search (APS)}\label{alg:APS}
{\small
\begin{algorithmic}[1]
\Require Instance $(\boldsymbol{r}, \widehat{\boldsymbol{u}}, N, \mathcal A, p)$, ambiguity radius $\epsilon \in [0, \epsdag)$, and outer tolerance $\delta \in (0, 1)$
\Ensure Incumbent pair $(\widehat{\mathcal S},\widehat\lambda)$, its value $h_{\mathrm{best}}=h_{\widehat{\mathcal S}}(\widehat\lambda)$, and absolute certificate gap $\Gamma$ satisfying
$R^{\mathrm{FC}}_{N,\epsilon}(\widehat{\mathcal S})=h_{\mathrm{best}}\ge(1-\delta)V^{\star}_{N,\epsilon}$
\State Compute $L^- \gets \epsilon$, \ $L^+ \gets \frac{1}{N}\sum_{n \in \mathcal{N}} \Delta_n - \epsilon$, and $\lambda^u$ by the ratio scan over~\eqref{eq:D_n}--\eqref{eq:Phi}
\State $(\mathcal S_0,f(0))\gets\mathrm{SBMB}(0)$ and $(\mathcal S_u,f(\lambda^u))\gets\mathrm{SBMB}(\lambda^u)$
\State Record $h(0)\gets f(0)$ and $h(\lambda^u)\gets-\lambda^u\epsilon+f(\lambda^u)$
\State $\mathcal P\gets\operatorname{distinct}\{\mathcal S_0,\mathcal S_u\}$; evaluate [Pool-BP($\mathcal S$)] once for each $\mathcal S\in\mathcal P$
\State Choose $\widehat{\mathcal S}\in\arg\max_{\mathcal S\in\mathcal P}H_{\mathcal S}$ and set $(\widehat\lambda,h_{\mathrm{best}})\gets(\lambda_{\widehat{\mathcal S}},H_{\widehat{\mathcal S}})$
\State Set $\IAPS\gets\{[0,\lambda^u]\}$
\While{$\IAPS\neq\emptyset$}
    \State $[a,b]\gets\arg\max_{[c,d]\in\IAPS}\mathrm{UB}^*([c,d])$
    \State $\mathrm{ub}_{\max}\gets\mathrm{UB}^*([a,b])$ and $\Gamma\gets\mathrm{ub}_{\max}-h_{\mathrm{best}}$
    \If{$\Gamma\le\delta\,\mathrm{ub}_{\max}$}
        \State \Return $(\widehat{\mathcal S},\widehat\lambda,h_{\mathrm{best}},\Gamma)$
    \EndIf
    \State $\lambda_{\mathrm{new}}\gets(a+b)/2$ and $(\mathcal S_{\mathrm{new}},f(\lambda_{\mathrm{new}}))\gets\mathrm{SBMB}(\lambda_{\mathrm{new}})$
    \State Record $h(\lambda_{\mathrm{new}})\gets-\lambda_{\mathrm{new}}\epsilon+f(\lambda_{\mathrm{new}})$
    \If{$\mathcal S_{\mathrm{new}}\notin\mathcal P$}
        \State Evaluate [Pool-BP($\mathcal S_{\mathrm{new}}$)] and set $\mathcal P\gets\mathcal P\cup\{\mathcal S_{\mathrm{new}}\}$
        \If{$H_{\mathcal S_{\mathrm{new}}}>h_{\mathrm{best}}$}
            \State $(\widehat{\mathcal S},\widehat\lambda,h_{\mathrm{best}})\gets(\mathcal S_{\mathrm{new}},\lambda_{\mathcal S_{\mathrm{new}}},H_{\mathcal S_{\mathrm{new}}})$
        \EndIf
    \EndIf
    \State $\IAPS\gets(\IAPS\setminus\{[a,b]\})\cup\{[a,\lambda_{\mathrm{new}}],[\lambda_{\mathrm{new}},b]\}$
    \State $\IAPS\gets\{[c,d]\in\IAPS:\mathrm{UB}^*([c,d])>h_{\mathrm{best}}\}$
\EndWhile
\State \Return $(\widehat{\mathcal S},\widehat\lambda,h_{\mathrm{best}},0)$ \Comment{all remaining intervals have been pruned}
\end{algorithmic}}
\end{algorithm}

Each product line enters $\mathcal P$ at most once, so its breakpoint set is evaluated once. The resulting value is the exact robust revenue of that line and hence a valid lower bound on $V^{\star}_{N,\epsilon}$, but it need not equal $h(\lambda_{\mathcal S})$, which also optimizes over product lines outside the pool. These fixed-line values update $h_{\mathrm{best}}$ and strengthen pruning; the exact SBMB values enter~\eqref{eq:ub_star}.

The following result establishes the end-to-end approximation guarantee and bounds the number of SBMB calls.

\begin{theorem}[Optimality Guarantee]\label{prop:delta_optimal}
Suppose $0<\epsilon<\epsdag$ and $\delta\in(0,1)$. Algorithm~\ref{alg:APS} returns a feasible pair
$(\widehat{\mathcal S},\widehat\lambda)$ satisfying
\[
    R^{\mathrm{FC}}_{N,\epsilon}(\widehat{\mathcal{S}})
    \;=\;h_{\mathrm{best}}
    \;\ge\;(1-\delta)V^{\star}_{N,\epsilon}.
\]
Let $K=L^+L^-/(L^++L^-)$ and let $B$ denote the total number of SBMB calls, including the two endpoint calls. Then
\[
B<2\max\left\{1,\frac{K\lambda^u}{\delta V^{\star}_{N,\epsilon}}\right\}+1.
\]
\end{theorem}
The maximum covers the possibility that the initial interval already satisfies the stopping condition, in which case $B=2$. For fixed instance-dependent quantities, the bound grows at rate $1/\delta$. It depends on the unknown value $V^{\star}_{N,\epsilon}$, but the algorithm itself does not: APS stops using the observable certificate gap $\Gamma$. Section~\ref{sec:experiments} reports the observed $B$ for $\delta\in\{0.01,0.05\}$.

The next corollary combines this computational certificate with the finite-sample guarantee in Theorem~\ref{thm:loss}.

\begin{corollary}[Population Performance of APS]
\label{cor:aps_population}
Under the sampling conditions of Theorem~\ref{thm:loss}, fix
$\epsilon>0$ and $\delta\in(0,1)$. If
$\epsilon<\epsilon^\dagger(\widehat Q_N)$, let
$\widehat{\mathcal S}$ be the product line returned by
Algorithm~\ref{alg:APS}; otherwise, set
$\widehat{\mathcal S}=\varnothing$. Let
$V_\epsilon^\star(Q)=\max_{\mathcal S\in\mathcal L}
R_\epsilon(\mathcal S;Q)$. Then
\begin{equation}
\label{eq:aps_population_expectation}
\mathbb E\!\left[R_\epsilon(\widehat{\mathcal S};Q)\right]
\ge
(1-\delta)V_\epsilon^\star(Q)
-
(2-\delta)\frac{\bar r}{\sqrt{2N}}
\left(
\sqrt{\pi}+\frac{1}{2}+\sqrt{1+\log|\mathcal L|}
\right),
\end{equation}
where the expectation is over the $N$ sampled utility vectors.
\end{corollary}
Thus, $\delta$ controls optimization error on the empirical problem, while the expected loss in~\eqref{eq:expected-loss} controls the error from replacing the population utility distribution by the empirical sample.

\begin{remark}[Finite-Sample Revenue Resolution]
Consider two product lines whose customer-level worst-case revenue profiles differ only for one customer, whose selected product changes from $i$ to $j$. Their average-revenue difference is $|r_i-r_j|/N$. Hence the smallest revenue change of this type is $N^{-1}\min_{i\ne j}|r_i-r_j|$,
which is $O(1/N)$. This finite-sample resolution is not used in the APS stopping rule; the guarantee in Theorem~\ref{prop:delta_optimal} is based on the certificate gap $\Gamma$.
\end{remark}

\section{Computational Experiments}\label{sec:experiments}

Section~\ref{sec:delta_optimal} develops APS and its finite-termination certificate. We now examine the computational cost of that certificate and compare the resulting solutions with those of the direct formulation [DRO-MILP]. The experiments vary the candidate set, customer sample, cardinality bound, and ambiguity radius to identify where the two approaches begin to separate.

\subsection{Experimental Design}\label{sec:exp_setup}

Each product $i$ has a distinct profile $\boldsymbol z_i$, sampled without replacement from $\{0,1\}^{10}\setminus\{\boldsymbol 0\}$, and a price $r_i\sim\mathrm{Uniform}(1,3)$ that also serves as unit revenue. We first draw a common partworth vector $\bar{\boldsymbol\theta}$ with independent $\mathrm{Uniform}(-0.75,0.75)$ coordinates. Customer $n$ has $\boldsymbol\theta_n=\bar{\boldsymbol\theta}+\boldsymbol\xi_n$, where the coordinates of $\boldsymbol\xi_n$ are independently drawn from $\mathrm{Uniform}(-0.70,0.70)$. Normalizing the magnitude of the price coefficient to one gives $\hat u_{ni}=\boldsymbol\theta_n^\top\boldsymbol z_i-r_i$.
The outside option has zero utility and revenue. The common partworths induce correlation across product utilities, while the continuous draws make utility ties a zero-probability event.

For each of 10 seeds, we generate a master instance with $A=N=200$; smaller instances retain the first $A$ products and first $N$ customers. This nesting preserves paired comparisons as the problem grows. We vary $A,N\in\{50,100,200\}$, $p\in\{5,10\}$, and $\epsilon\in\{0.25\epsdag,0.75\epsdag\}$, where $\epsdag$ is the critical radius in Corollary~\ref{corollary:welfare_lb}. We refer to these radii as low~(L) and high~(H), respectively. The design contains 36 configurations and 360 instances.

We compare APS at $\delta\in\{0.01,0.05\}$ with [DRO-MILP]
from Section~\ref{appendix:milp}. We solve [DRO-MILP] with
Gurobi~11.0.3 using one thread and a relative gap tolerance of
$10^{-4}$. Each SBMB call uses the same relative gap tolerance.
APS uses the upper bounds returned by SBMB to construct its
interval bounds. The outer tolerance $\delta$ applies to the
resulting global certificate for the empirical problem.
Each end-to-end run has a one-hour limit.
All experiments use a Windows~10 workstation with an Intel Xeon W-2245 CPU (3.90~GHz) and 64~GB RAM. Reported times cover the complete optimization routine: [DRO-MILP] includes model construction, while APS includes computation of $\lambda^u$ and all SBMB calls. Instance generation, calibration of $\epsdag$, and one-time compilation are excluded.

\subsection{Computational Performance}\label{sec:exp_results}

Each row of Table~\ref{tab:main_results} summarizes the
10 replications of one configuration. Time and $B$, the total
number of SBMB calls, are reported as arithmetic means.
Solved counts the [DRO-MILP] runs that reach the prescribed
$10^{-4}$ relative gap tolerance within one hour. A time-limited
run contributes 3600 seconds to the mean time. We compare the
returned objective values using
$\Delta=(h_{\mathrm{APS}}-h_{\mathrm{MILP}})/|h_{\mathrm{MILP}}|
\times100\%$, where $h_{\mathrm{APS}}$ is the objective value
returned by APS and $h_{\mathrm{MILP}}$ is the final
[DRO-MILP] incumbent objective. A positive $\Delta$ favors APS.

\begin{table}[htbp]
\centering
\caption{End-to-end comparison of APS and [DRO-MILP].}\label{tab:main_results}
\scriptsize
\renewcommand{\arraystretch}{1.0}
\setlength{\tabcolsep}{4pt}
\begin{tabular*}{\textwidth}{@{\extracolsep{\fill}}cccc rr rrr rrr@{}}
\toprule
& & & & \multicolumn{2}{c}{[DRO-MILP]} & \multicolumn{3}{c}{APS ($\delta{=}0.01$)} & \multicolumn{3}{c}{APS ($\delta{=}0.05$)} \\
\cmidrule(lr){5-6} \cmidrule(lr){7-9} \cmidrule(lr){10-12}
$p$ & $A$ & $N$ & $\epsilon$ & Mean Time & Solved & Mean Time & Mean $B$ & Mean $\Delta$ (\%) & Mean Time & Mean $B$ & Mean $\Delta$ (\%) \\
\midrule
5 & 50 & 50  & H &    1.1 & 10/10 & 0.07 & 26.4 & 0.00 & 0.05 & 14.5 & 0.00 \\
  &    &     & L &    1.9 & 10/10 & 0.07 & 30.0 & 0.00 & 0.04 & 13.2 & 0.00 \\
  &    & 100 & H &    3.4 & 10/10 & 0.12 & 37.0 & 0.00 & 0.08 & 18.6 & 0.00 \\
  &    &     & L &    5.6 & 10/10 & 0.09 & 30.5 & 0.00 & 0.06 & 13.7 & 0.00 \\
  &    & 200 & H &   11.4 & 10/10 & 0.22 & 41.5 & 0.00 & 0.17 & 19.5 & 0.00 \\
  &    &     & L &   17.2 & 10/10 & 0.22 & 32.9 & 0.00 & 0.16 & 13.6 & 0.00 \\
\cmidrule(l){2-12}
  &100 & 50  & H &    4.8 & 10/10 & 0.20 & 40.0 & 0.00 & 0.13 & 18.6 & 0.00 \\
  &    &     & L &   10.0 & 10/10 & 0.18 & 29.2 & 0.00 & 0.11 & 12.8 & 0.00 \\
  &    & 100 & H &   14.2 & 10/10 & 0.30 & 39.4 & 0.00 & 0.23 & 19.6 & 0.00 \\
  &    &     & L &   38.1 & 10/10 & 0.35 & 29.5 & 0.00 & 0.25 & 12.8 & 0.00 \\
  &    & 200 & H &   34.3 & 10/10 & 0.63 & 46.0 & 0.00 & 0.51 & 20.7 & 0.00 \\
  &    &     & L &   98.3 & 10/10 & 0.79 & 30.3 & 0.00 & 0.63 & 12.9 & 0.00 \\
\cmidrule(l){2-12}
  &200 & 50  & H &   18.0 & 10/10 & 0.55 & 38.3 & 0.00 & 0.42 & 19.2 & 0.00 \\
  &    &     & L &  114.4 & 10/10 & 1.1  & 29.5 & 0.00 & 0.93 & 13.4 & 0.00 \\
  &    & 100 & H &   50.8 & 10/10 & 1.0  & 45.9 & 0.00 & 0.84 & 21.3 & 0.00 \\
  &    &     & L &  371.5 & 10/10 & 2.8  & 30.7 & 0.00 & 2.5  & 13.3 & 0.00 \\
  &    & 200 & H &  142.3 & 10/10 & 2.0  & 46.6 & 0.00 & 1.7  & 21.8 & 0.00 \\
  &    &     & L &  905.6 & 10/10 & 6.3  & 30.4 & 0.00 & 5.8  & 13.4 & 0.00 \\
\midrule
10& 50 & 50  & H &   0.81 & 10/10 & 0.09 & 25.2 & 0.00 & 0.06 & 13.8 & 0.00 \\
  &    &     & L &    1.6 & 10/10 & 0.08 & 29.5 & 0.00 & 0.04 & 13.0 & 0.00 \\
  &    & 100 & H &    2.7 & 10/10 & 0.17 & 34.2 & 0.00 & 0.12 & 16.1 & 0.00 \\
  &    &     & L &    4.6 & 10/10 & 0.18 & 31.7 & 0.00 & 0.12 & 13.5 & 0.00 \\
  &    & 200 & H &    9.1 & 10/10 & 0.45 & 36.7 & 0.00 & 0.39 & 16.9 & 0.00 \\
  &    &     & L &   13.1 & 10/10 & 0.53 & 31.7 & 0.00 & 0.44 & 13.2 & 0.00 \\
\cmidrule(l){2-12}
  &100 & 50  & H &    3.5 & 10/10 & 0.28 & 37.8 & 0.00 & 0.17 & 16.8 & 0.00 \\
  &    &     & L &    7.6 & 10/10 & 0.25 & 30.8 & 0.00 & 0.15 & 13.1 & 0.00 \\
  &    & 100 & H &   12.4 & 10/10 & 0.45 & 34.0 & 0.00 & 0.34 & 16.5 & 0.00 \\
  &    &     & L &   45.3 & 10/10 & 0.65 & 29.6 & 0.00 & 0.46 & 12.6 & 0.00 \\
  &    & 200 & H &   34.2 & 10/10 & 1.2  & 37.4 & 0.00 & 0.99 & 16.7 & 0.00 \\
  &    &     & L &  205.6 & 10/10 & 2.5  & 30.8 & 0.00 & 2.1  & 12.8 & 0.00 \\
\cmidrule(l){2-12}
  &200 & 50  & H &   14.2 & 10/10 & 0.72 & 31.1 & 0.00 & 0.54 & 15.8 & 0.00 \\
  &    &     & L &  162.9 & 10/10 & 1.8  & 29.2 & 0.00 & 1.3  & 12.8 & 0.00 \\
  &    & 100 & H &   56.8 & 10/10 & 3.2  & 38.1 & 0.00 & 2.9  & 17.6 & 0.00 \\
  &    &     & L & 1414.1 &  7/10 & 25.3 & 30.8 & 0.00 & 23.7 & 13.0 & 0.00 \\
  &    & 200 & H &  281.0 & 10/10 & 9.3  & 39.1 & 0.00 & 8.7  & 17.2 & 0.00 \\
  &    &     & L & 2372.8 &  5/10 &165.9 & 30.5 & 0.14 &160.6 & 12.3 & 0.14 \\
\bottomrule
\end{tabular*}
\end{table}
APS meets its prescribed relative tolerance on all 360 instances
at both $\delta=0.01$ and $\delta=0.05$, as certified by its
own global bounds. [DRO-MILP] reaches its $10^{-4}$ relative
gap tolerance on 352 instances within one hour. On those
instances, the objective values returned by APS and [DRO-MILP]
agree to numerical precision. Among the eight time-limited
[DRO-MILP] runs, APS matches the final incumbent in six and
improves on it by $0.18\%$ and $1.19\%$ in the other two.

APS is faster in every paired run, including those in which
[DRO-MILP] reaches its prescribed tolerance. The mean paired
speedups are $52.9\times$ for $\delta=0.01$ and $70.7\times$
for $\delta=0.05$; the corresponding maxima are $508.0\times$
and $662.4\times$. The advantage is therefore not driven by the
eight time-limited benchmark runs.

Scale separates the methods most sharply. With $p=10$, $N=200$,
and the low radius, the mean [DRO-MILP] times at $A=50,100,200$
are 13.1, 205.6, and 2372.8 seconds, compared with 0.44, 2.1,
and 160.6 seconds for APS at $\delta=0.05$. All eight [DRO-MILP]
time limits occur at $A=200$. On the largest configuration,
the mean paired objective improvement over the final [DRO-MILP]
incumbent is $0.14\%$.

Both tolerances return the same product line, multiplier, and objective value on every instance. Reducing $\delta$ from $0.05$ to $0.01$ raises the mean number of SBMB calls from 15.5 to 34.0 and increases the mean paired runtime by a factor of $1.38$, tightening the certificate without changing the solution.

\begin{remark}[Decomposing the Computational Gain]
Table~\ref{tab:main_results} reports the end-to-end performance of APS. Section~\ref{appendix:sbmb_benchmark} of the online supplement fixes $\lambda$ and compares SBMB with [MILP], thereby isolating the contribution of the inner solver.
\end{remark}

\section{Case Study on Conjoint Data}\label{sec:Case}

Section~\ref{sec:experiments} establishes the computational performance of the
proposed algorithms on synthetic instances. We now turn to two empirical
questions: how the proposed model performs with new customers, and how
robustness reshapes the product line and thereby produces those outcomes. The
application concerns a smart vending machine whose limited capacity requires
the seller to choose a small menu of meal boxes from a large collection of
feasible configurations.

\subsection{Data and Experimental Protocol}\label{sec:Case::Data}

The data come from a choice-based conjoint survey of $45$ students. Each
respondent completed $10$ tasks involving three meal profiles and a
no-purchase option. Table~\ref{tab:case_levels} lists the five meal attributes.
We estimate respondent-level partworths using the generalized robust conjoint
method of \citet{evgeniou2005generalized}. After normalizing the magnitude of
the price coefficient to one, a meal's utility equals the sum of its four
non-price partworths minus its price. The estimates form a $45\times16$ matrix
of non-price partworths; price enters utility directly. The survey design, raw
responses, and estimated partworths are available at
\url{https://github.com/mogoooojay/Robust-PLD.git}.

\begin{table}[htbp]
\centering
\caption{Meal attributes and level indices.}
\label{tab:case_levels}
\small
\begin{tabular}{p{2.5cm}p{12cm}}
\toprule
\textbf{Attribute} & \textbf{Levels (global index)} \\
\midrule
Main Dish     & Chicken Breast (1), Shrimp (2), Fish Fillet (3), Beef Steak (4) \\
Carbohydrate  & Brown Rice (5), Ciabatta (6), Sweet Potato (7), Corn (8) \\
Side Dish     & Steamed Broccoli (9), Boiled Vegetables (10), Roasted Mushrooms (11), Salad (12) \\
Flavor        & Mild (13), Spicy (14), Garlic Herb (15), Teriyaki (16) \\
Price         & \textyen 8 (17), \textyen 10 (18), \textyen 12 (19), \textyen 15 (20) \\
\bottomrule
\end{tabular}
\end{table}

We fit a kernel density estimate with bandwidth $0.5$ to the respondent-level
partworth vectors and use it to generate repeated instances. In each
replication, we draw $A+20$ distinct non-price profiles from the $4^4=256$
possible combinations and independently assign each a price drawn uniformly
from
$\{\text{\textyen}8,\text{\textyen}10,\text{\textyen}12,\text{\textyen}15\}$.
The first $A$ meals form the seller's candidate set. The remaining $20$ are
competing alternatives, and their maximum utility defines each customer's
outside option.

For each $(A,p,N)$ configuration, we conduct $100$ replications. A replication
uses $N$ training customers and $100$ independently drawn validation customers.
For testing, we draw an independent sample of $10{,}000$ partworth vectors once
and hold it fixed across replications. Within each replication, all methods face
the same products and customer samples. The validation sample is used only for
model selection, and the test sample only for evaluation, so all reported
method comparisons are paired.

We measure utility relative to the competing alternatives. For each customer,
we normalize utilities as $\widetilde u_{ni}=u_{ni}-u_{n0}$ and set
$\widetilde u_{n0}=0$. This normalization preserves within-customer utility
differences and first choices. All optimization and evaluation use these
normalized utilities. For any customer sample $\mathcal D$, buyers' welfare is
$W(\mathcal S;\mathcal D)=|\mathcal D|^{-1}\sum_{n\in\mathcal D}\max_{i\in\mathcal S\cup\{0\}}\widetilde u_{ni}$.
With unit price sensitivity, $W(\mathcal S;\mathcal D)$ has a money-metric interpretation.
Revenue is reported per customer, in the same monetary units as price.

The MNL-based benchmarks use a representative-consumer model constructed from
the same training sample. Its attraction parameters are
$\bar v_i=\exp(N^{-1}\sum_{n=1}^N\widetilde u_{ni})$, so $\bar v_0=1$. For a
product line $\mathcal S$, its nominal revenue and welfare are
\[
R_{\mathrm{MNL}}(\mathcal S)
=\frac{\sum_{i\in\mathcal S}r_i\bar v_i}
{1+\sum_{i\in\mathcal S}\bar v_i},
\qquad
U_{\mathrm{MNL}}(\mathcal S)
=\log\!\left(1+\sum_{i\in\mathcal S}\bar v_i\right).
\]

We compare five methods under the common cardinality constraint
$|\mathcal S|\le p$:
\begin{itemize}
\item \textbf{Nominal-FC} maximizes empirical first-choice revenue.

\item \textbf{R\&T 2012} is the robust MNL model of
\citet{rusmevichientong2012robust}. We use their rectangular attraction set
$\mathcal V=\prod_{i\in\Aplus}[(1-\kappa)\bar v_i,(1+\kappa)\bar v_i]$,
which includes the outside-option attraction.

\item \textbf{B\&M 2017} uses the robust product-line design framework of
\citet{bertsimas2016robust}, which covers parameter uncertainty within a choice
model as well as structural uncertainty across choice models. We adopt their
first-choice box specification from Section~4.2 because it is the closest
counterpart to our empirical first-choice model. It retains the sampled utility
vectors and choice rule but allows the population weights of the sampled
customer types to vary according to
$\max\{0,(1-\eta)/N\}\leq \pi^k\leq(1+\eta)/N$, with
$\sum_k\pi^k=1$.

\item \textbf{Sumida et al. 2021} follows the direct revenue--utility tradeoff
of \citet{sumida2021revenue}. Using the representative-consumer MNL model
above, its supported solutions maximize
$R_{\mathrm{MNL}}(\mathcal S)+\lambda U_{\mathrm{MNL}}(\mathcal S)$ over
different values of $\lambda$.

\item \textbf{Ours} is the proposed $\ell_1$-Wasserstein DRO model. It maximizes
worst-case first-choice revenue over a Wasserstein ball centered at the
empirical utility distribution.
\end{itemize}

Nominal-FC requires no calibration. The next subsection specifies how the
remaining methods are calibrated and then compares the resulting product lines
on a common test market.

\subsection{The ``Doing Well by Doing Good'' Mechanism}
\label{sec:case_v2:mechanism}

\begin{figure}[htbp]
    \centering
    \includegraphics[width=0.8\textwidth]{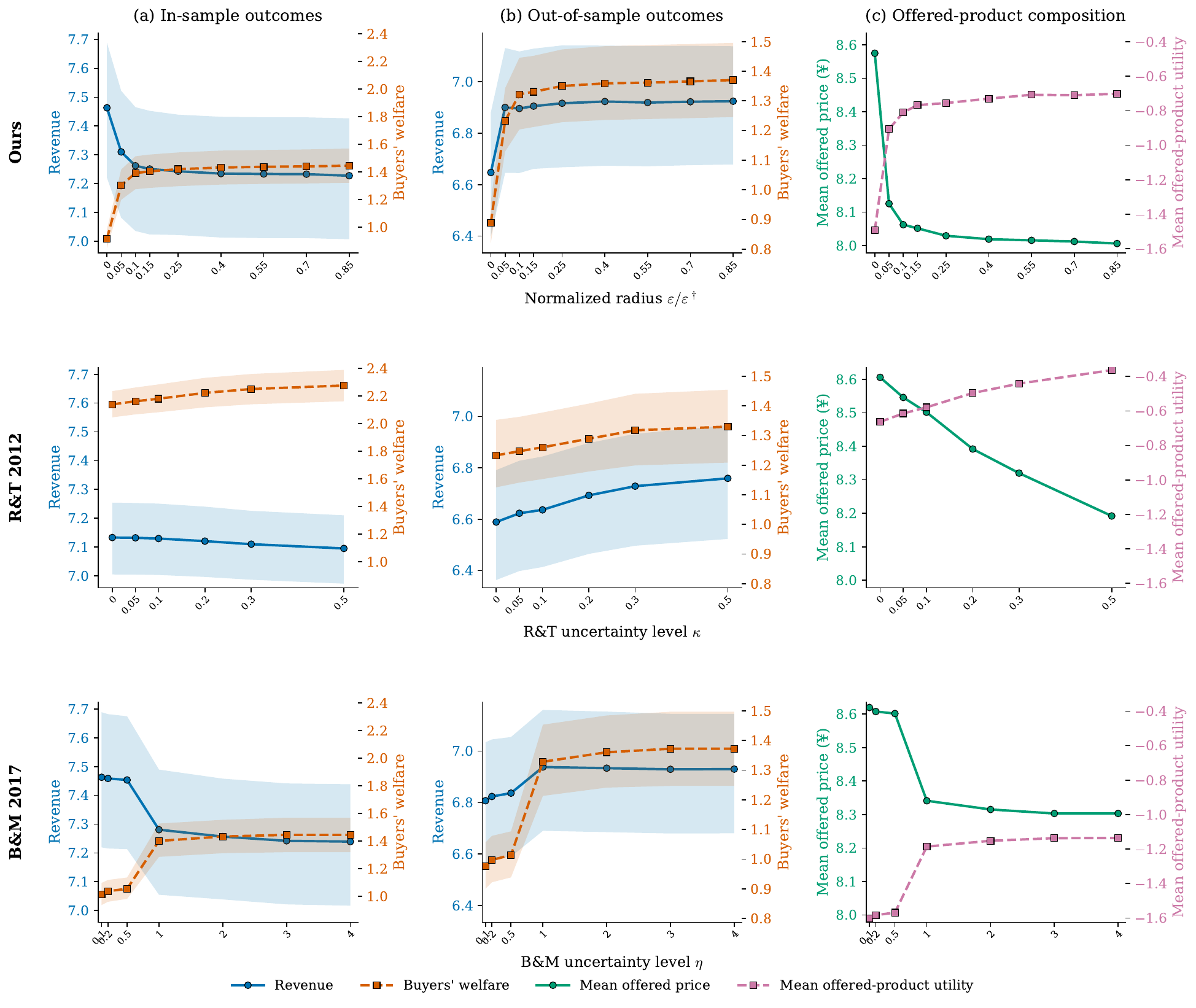}
    \caption{Revenue, welfare, and offered-product composition along the
    uncertainty paths.}
    \label{fig:case_v2_mechanism_paths}
\end{figure}

We first validate the ``doing good'' incentive in the robust product line design across all robust models, including R\&T, B\&M, and Ours, and then show that the resulting gain in buyers' welfare is repaid in out-of-sample revenue, so that the model also ``does well.''
We fix $(A,p,N)=(200,10,30)$ and use
Figure~\ref{fig:case_v2_mechanism_paths} to trace each robust model across its
uncertainty-parameter grid. The resulting paths connect changes in product-line
composition to revenue and buyers' welfare. Our model moves utility vectors,
R\&T bounds aggregate MNL attractions, and B\&M reweights sampled customer
types. Because $\epsilon/\epsdag$, $\kappa$, and $\eta$ measure different forms
of uncertainty, their magnitudes and slopes are meaningful within a row, not
across rows.

Panel~(a) evaluates each line under the nominal model from which it was
constructed: first-choice for our model and B\&M, and MNL for R\&T. These are
nominal in-sample outcomes, not robust objective values. Panel~(b) evaluates
every line on the common first-choice test population. Panel~(c) reports the
mean price of the offered products and mean utility across all offered
products and test customers, measured relative to the outside option. These
composition measures describe the line itself, not the products eventually
purchased. The zero endpoint is added for our model and R\&T. Each point is a
mean over $100$ replications; the bands in panels~(a) and~(b) are $95\%$
percentile-bootstrap confidence intervals.

Along all three paths, stronger protection shifts the selected product line
toward lower prices and higher utility relative to the outside option. The
broader appeal raises purchase incidence, and by the upper end of each grid the
additional purchases more than offset the decline in revenue per purchase:
both out-of-sample revenue and welfare exceed their low-protection levels. The
adjustment need not be monotone at every grid point, and its source differs
across models.

The top row traces the proposed model. As $\epsilon/\epsdag$ increases from $0$
to $0.85$, nominal in-sample revenue falls from $7.463$ to $7.227$, while
welfare rises from $0.914$ to $1.444$. Mean offered price falls from
\textyen$8.57$ to \textyen$8.01$, and mean utility rises from $-1.49$ to
$-0.70$. On the test population, purchase incidence increases from $0.809$ to
$0.865$, more than offsetting the decline in revenue per purchase from $8.214$
to $8.003$. Out-of-sample revenue consequently rises from $6.647$ to $6.925$,
alongside a welfare increase from $0.889$ to $1.370$. Most of the revenue gain
appears at the first positive radius, whereas welfare continues to improve
along the path. The dual representation in
Proposition~\ref{prop:dual_worst_case_revenue} explains why this movement arises
naturally under [DRO-FC]: robust revenue contains the term
$\lambda W_N(\mathcal S)$. Thus, buyers' welfare enters the seller's primary
robust-revenue objective rather than serving only as an ex post performance
measure.

For R\&T, the adjustment originates in aggregate-attraction protection. At the
relevant vertex of its uncertainty set, the
product-to-outside attraction ratio is scaled by
$(1-\kappa)/(1+\kappa)$. Increasing $\kappa$ therefore favors lines that retain
more aggregate attraction under the adverse model. From $\kappa=0$ to $0.5$,
mean offered price falls from \textyen$8.61$ to \textyen$8.19$, and mean
utility rises from $-0.66$ to $-0.36$. In the representative-consumer MNL
model, in-sample revenue decreases slightly from $7.133$ to $7.095$, while
welfare increases from $2.139$ to $2.275$. On the common first-choice test
population, purchase incidence rises from $0.816$ to $0.845$ and revenue per
purchase falls from $8.088$ to $8.004$. The increase in purchases dominates:
out-of-sample revenue rises from $6.589$ to $6.759$, together with a welfare
increase from $1.233$ to $1.330$. The capacity bound is binding throughout, so
the change comes from product replacement rather than line expansion;
Section~\ref{appendix:capacity_separation} discusses this distinction.

B\&M acts instead on customer-type weights. For a fixed line, its
adversary shifts probability toward sampled types that generate lower revenue,
penalizing lines whose performance relies on favorable population weights.
Within each B\&M solve, the in-sample tie rule specified in
Section~\ref{sec:case_v2:headline} invokes welfare only when robust revenue is
unchanged. From $\eta=0.1$ to $4$, nominal in-sample
revenue falls from $7.463$ to $7.239$, while welfare rises from $1.013$ to
$1.445$. Over the same path, mean offered price falls from \textyen$8.62$ to
\textyen$8.30$, and mean utility rises from $-1.60$ to $-1.14$. Test-market
purchase incidence increases from $0.831$ to $0.866$,
offsetting a decrease in revenue per purchase from $8.188$ to $8.002$.
Out-of-sample revenue reaches $6.938$ at $\eta=1$ and remains near $6.93$
thereafter; welfare rises from $0.976$ at $\eta=0.1$ to $1.372$ at $\eta=4$.
Most of the adjustment occurs between $\eta=0.5$ and $1$.

The three rows show that the ``doing well by doing good'' mechanism extends
beyond a single ambiguity set. In this application, utility perturbations,
aggregate-attraction protection, and customer-type reweighting all shift the
selected product line toward broader appeal, although through different
channels. [DRO-FC] is distinctive because its robust-revenue representation
contains the explicit term $\lambda W_N(\mathcal S)$, so buyers' welfare enters
the seller's primary objective. Proposition~\ref{prop:dual_worst_case_revenue}
identifies this link, and Corollary~\ref{corollary:welfare_lb} converts it into
a formal welfare guarantee.

\subsection{Out-of-Sample Comparison}\label{sec:case_v2:headline}

We then compare the out-of-sample performance of our model against the benchmarks across varying settings. R\&T, B\&M, and Sumida et al., which select their uncertainty levels on the independent
validation sample. Regardless of the model used to construct a product line,
we score it by first-choice validation revenue and welfare. For each method, we
solve the full parameter grid and retain the line with the highest validation
revenue; exact revenue ties are broken by higher validation welfare and then by
the larger uncertainty parameter.

We use the grids $\kappa\in\{0.05,0.10,0.20,0.30,0.50\}$, $\eta\in\{0.10,0.20,0.50,1,2,3,4\}$, and $\epsilon/\epsdag\in\{0.05,0.10,0.15,0.25,0.40,0.55,0.70,0.85\}$,
where $\epsdag$ is the critical radius in
Corollary~\ref{corollary:welfare_lb}. We exclude zero because it reproduces
Nominal-FC for our model and B\&M, and the nominal MNL endpoint for R\&T. The
B\&M grid follows the positive values in Section~4.2 of their paper. At each value of $\eta$, we
run the B\&M local search from $50$ random starts, add the exact Nominal-FC line
and the line from the preceding value of $\eta$ as warm starts, and retain the
solution with the highest robust objective. To accommodate the common
constraint $|\mathcal S|\le p$, its one-product neighborhood allows additions
and deletions as well as swaps. Exact robust-objective ties are
resolved by higher in-sample buyers' welfare and then by higher nominal
in-sample revenue; this comparison uses neither validation nor test outcomes.
Nominal-FC and R\&T are solved exactly. At each value of $\epsilon$, our model
is solved by APS
(Section~\ref{sec:delta_optimal}) with relative tolerance $\delta=0.01$.

The Sumida et al. benchmark is calibrated by a revenue-loss allowance rather
than validation. On the supported frontier generated by
$R_{\mathrm{MNL}}(\mathcal S)+\lambda U_{\mathrm{MNL}}(\mathcal S)$, we select
the highest-utility point satisfying
$R_{\mathrm{MNL}}(\mathcal S)\ge0.99R_{\mathrm{MNL}}^\star$, where
$R_{\mathrm{MNL}}^\star$ is the exact revenue optimum under that MNL model.
Thus, $1\%$ is the permitted revenue loss under the model, not a fixed value
of $\lambda$. The supported frontier and its revenue benchmark are computed
exactly.

After calibration, we evaluate every selected product line on the common
$10{,}000$-customer first-choice test sample. The test sample plays no role in
model selection.

\begin{figure}[htbp]
    \centering
    \includegraphics[width=0.8\textwidth]{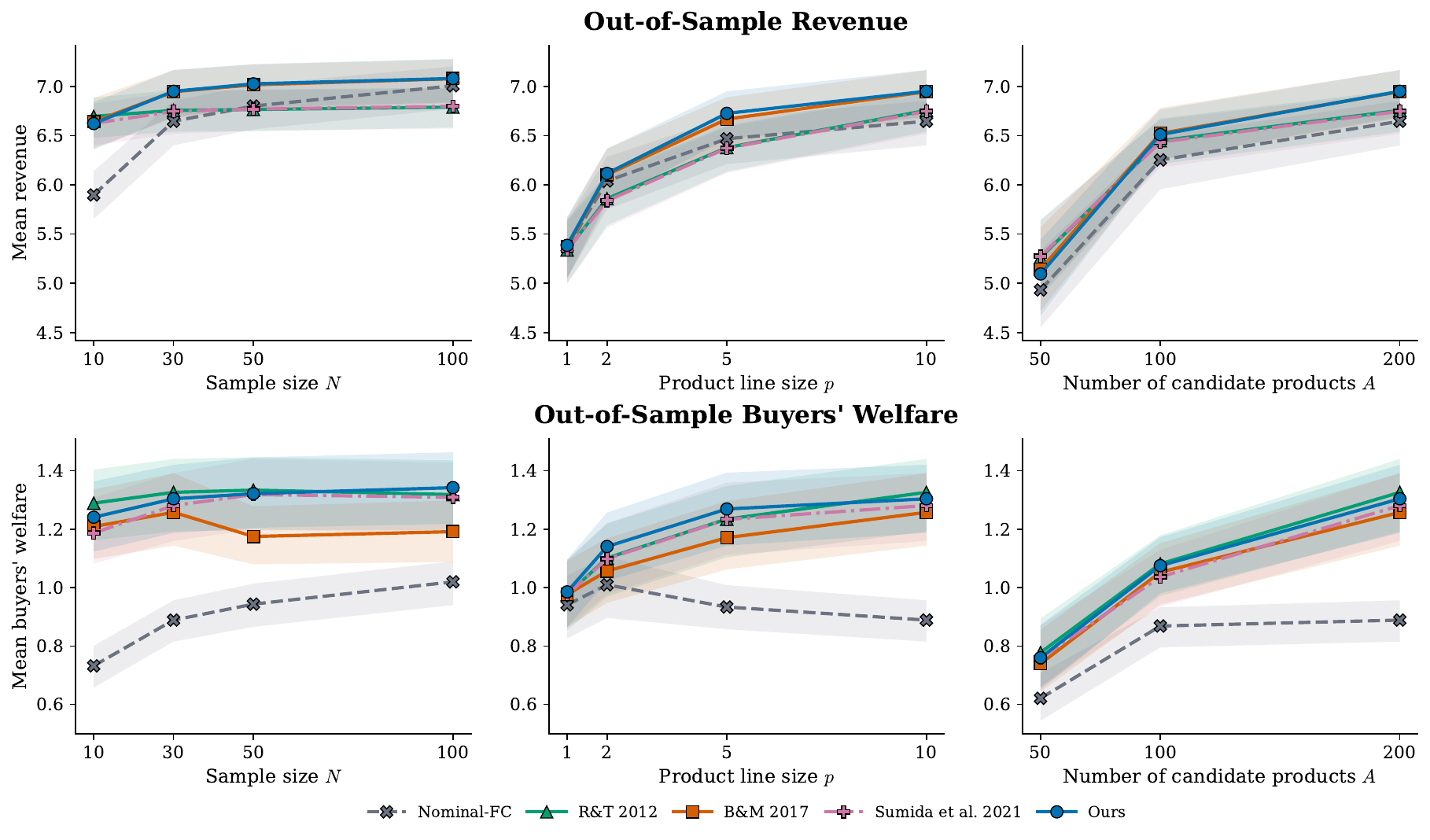}
    \caption{Out-of-sample revenue and buyers' welfare after model calibration.}
    \label{fig:case_v2_sensitivity}
\end{figure}

We vary $N$, $p$, and $A$ one at a time around the default setting
$(A,p,N)=(200,10,30)$. Thus, $(A,p)=(200,10)$ when $N$ varies,
$(A,N)=(200,30)$ when $p$ varies, and $(p,N)=(10,30)$ when $A$ varies.
Figure~\ref{fig:case_v2_sensitivity} reports outcome levels for all five methods.
Each marker is the mean over $100$ replications, and the shaded bands are
$95\%$ percentile-bootstrap confidence intervals for the mean.

Across the nine distinct configurations, our model's mean revenue and welfare
both exceed Nominal-FC in every case. It attains the highest mean revenue in
six configurations and ranks first or second in mean welfare in all nine. At
the default setting, its mean revenue is $6.951$, narrowly above B\&M at
$6.947$. Mean welfare is $1.304$, $0.022$ below R\&T and above the other three
methods.

The amount of training data changes the ranking. At $N=10$, all four
alternatives to Nominal-FC yield higher revenue. R\&T leads with a gain of
$0.799$, followed by B\&M at $0.745$; Sumida et al. and our model are close
behind, at $0.727$ and $0.725$. Our model moves into first place at $N=30$ and
remains there at $N=50$ and $100$, although B\&M is close at each of these
sample sizes. Its gain over Nominal-FC narrows from $0.725$ at $N=10$ to
$0.074$ at $N=100$ as the empirical distribution becomes more informative.
The two MNL-based benchmarks benefit most from the smallest sample, but their
revenue gaps turn negative at $N=50$ and $100$. Thus, the aggregate MNL
approach is most competitive with limited data, whereas the individual-level
utility information becomes more valuable as $N$ grows.

The apparent drop in B\&M welfare is a calibration effect. From $N=30$ to
$N=50$, validation selects $\eta=4$ less often ($79$ versus $63$ replications),
and mean welfare falls from $1.258$ to $1.175$ while mean revenue rises from
$6.947$ to $7.017$. At every fixed $\eta$, welfare instead increases.
Figure~\ref{fig:case_v2_sensitivity} therefore reports calibrated, rather than
fixed-parameter, performance.

Our model's advantage over Nominal-FC also grows with the size of the feasible
design space. Its revenue gain rises from $0.025$ at $p=1$ to $0.304$ at $p=10$,
and from $0.163$ at $A=50$ to $0.304$ at $A=200$. Sumida et al. leads at $A=50$,
B\&M leads narrowly at $A=100$, and our model leads at $A=200$. Larger values
of $p$ and $A$ give the seller more scope to fit a nominal line to
sample-specific preferences and, correspondingly, create more value from
distributional protection. The revenue gains are accompanied by strong welfare
performance: our model leads both measures at $N=100$ and at $p=1,2,5$, and
ranks first or second in welfare elsewhere.
Section~\ref{sec:case_v2:mechanism} next examines how the joint improvement
arises.

The outperformance of Sumida et al.\ over Nominal-FC provides further support
for the revenue--welfare alignment. Their model places greater weight on
buyers' welfare than Nominal-FC does when selecting the product line, which is
the ``doing good'' incentive. That emphasis improves both out-of-sample revenue
and buyers' welfare, which is the ``doing well'' payoff in market performance.

\section{Conclusion}\label{sec:Con}
This paper studies product line design when customer utilities are uncertain. A DRO model is formulated to maximize worst-case expected revenue over an $\ell_1$-Wasserstein ball around a nominal utility distribution. The analysis of the dual reveals a ``doing well by doing good'' mechanism: protecting revenue against utility uncertainty raises buyers' welfare, and that welfare in turn supports revenue. Under the MNL nominal with no cardinality constraint, the higher welfare is achieved by expanding the revenue-ordered assortment as the robustness radius grows. Under the empirical nominal, robustness shifts the line toward lower-priced products with greater utility relative to the outside option.
We also provide statistical guarantees on the revenue loss incurred when the model is applied to estimated or sampled utilities. The conjoint study shows the mechanism in the data: the robust line improves both out-of-sample revenue and buyers' welfare over the nominal benchmark in every configuration.
Our adaptive parametric search algorithm certifies that the
returned product line attains at least a $(1-\delta)$ fraction
of the empirical robust optimum. We separately bound its
expected robust-revenue loss relative to the population optimum,
accounting for both sampling error and the optimization tolerance.
On 360 instances it is $70.7\times$ faster on average and up to $662.4\times$ faster than the mixed-integer formulation.

Three directions remain open.  The first is a DRO formulation for choice models beyond random utility maximization. The second is attribute-level product line design, in which the decision is the attribute profile of each variant rather than a subset of a candidate list. The third is whether the ``doing well by doing good'' mechanism extends to other applications of robust models.

\par\bigskip
\appendix
\phantomsection
\addcontentsline{toc}{section}{Appendix}
\section*{Appendix}
\renewcommand{\thesection}{EC.\arabic{section}}
\renewcommand{\theHsection}{appendix.\arabic{section}}
\setcounter{equation}{0}
\setcounter{figure}{0}
\setcounter{table}{0}
\renewcommand{\theequation}{EC.\arabic{equation}}
\renewcommand{\thefigure}{EC.\arabic{figure}}
\renewcommand{\thetable}{EC.\arabic{table}}
\renewcommand{\theHequation}{appendix.\arabic{equation}}
\renewcommand{\theHfigure}{appendix.\arabic{figure}}
\renewcommand{\theHtable}{appendix.\arabic{table}}

\section{The [DRO-MILP] Reformulation}\label{appendix:milp}

We derive [DRO-MILP], an exact mixed-integer linear reformulation of [DRO-FC]. Let binary decision variables \(\boldsymbol{x} \in \mathbb{B}^{A+1}\) represent the product line, where \(x_i = 1\) if alternative \(i \in \Aplus\) is offered and \(x_i = 0\) otherwise. The feasible region of \(\boldsymbol{x}\), denoted by \(\mathcal{X}\), is in one-to-one correspondence with the feasible family \(\mathcal{L}\). Under the cardinality restriction used in Section~\ref{sec:delta_optimal}, \(\mathcal{L} = \{\mathcal{S} \subseteq \mathcal{A} : |\mathcal{S}| \le p\}\), and the corresponding decision space is
\[
\mathcal{X} = \left\{ \boldsymbol{x} \in \mathbb{B}^{A+1} : x_0 = 1,\ \sum_{i \in \Aplus} x_i \le p + 1 \right\}.
\]

With these definitions, the robust PLD problem can be reformulated as the following MIP:
\begin{align}
    \max_{\lambda \ge 0,\ \boldsymbol{x} \in \mathcal{X}} \quad & -\lambda\epsilon + \frac{1}{N} \left( \sum_{n \in \mathcal{N}} s_n + \lambda \sum_{n \in \mathcal{N}} \sum_{i \in \Aplus} \hat{u}_{ni} y_{ni} \right) \label{eq:obj} \\
    \text{s.t.} \quad
    & s_n \le r_i - \lambda \hat{u}_{ni} + \bar{C}_n(1 - x_i), && \forall n \in \mathcal{N},\ \forall i \in \Aplus, \label{eq:c1} \\
    & \sum_{i \in \Aplus} y_{ni} = 1, && \forall n \in \mathcal{N}, \label{eq:c2} \\
    & y_{ni} \le x_i, && \forall n \in \mathcal{N},\ \forall i \in \Aplus, \label{eq:c3} \\
    & \boldsymbol{y} \ge 0. \label{eq:c4}
\end{align}

Constraint~\eqref{eq:c1} enforces that \(s_n = \min_{i \in \bar{\mathcal{S}}} \{r_i - \lambda \hat{u}_{ni}\}\) at optimality, where \(\bar{C}_n\) is a constant large enough to deactivate the constraint for products that are not offered. The variables \(\boldsymbol{y} \in \mathbb{R}_+^{N \times (A+1)}\) represent customer choice across \(N\) samples. Constraints~\eqref{eq:c2} and \eqref{eq:c3} ensure that each customer allocates one unit of choice among the offered alternatives (including the outside option). Although \(\boldsymbol{y}\) is declared continuous, the formulation is exact in value: for fixed \((\boldsymbol{x},\lambda)\), the welfare term \(\lambda\sum_{i\in\Aplus}\hat{u}_{ni}y_{ni}\) is a linear function maximized over the polytope \(\{\boldsymbol{y}_n\ge0:\sum_{i\in\Aplus}y_{ni}=1,\ y_{ni}\le x_i\}\), whose vertices are integral; its maximum equals \(\lambda\max_{i\in\bar{\mathcal{S}}}\hat{u}_{ni}\) and is attained by placing all mass of \(\boldsymbol{y}_n\) on a utility-maximizing offered alternative. Utility ties leave the optimal value unchanged because tied alternatives carry the same coefficient \(\hat{u}_{ni}\). Averaging over customers, the optimal value of the MIP equals the [DRO-FC] value \(V^{\star}_{N,\epsilon}\).

Next, we apply a standard linearization technique to the bilinear term $\lambda y_{ni}$. By Lemma~\ref{prop:lambda_bounds}, [DRO-FC] admits an optimal multiplier in $[0, \lambda^u]$, which we use as the McCormick range. Introducing auxiliary variables \(\rho_{ni}\) in place of \(\lambda y_{ni}\) yields the following mixed-integer linear program:

\begin{align}
    \text{[DRO-MILP]} \qquad \max_{\lambda \in [0,\ \lambda^u],\ \boldsymbol{x} \in \mathcal{X}} \quad & -\lambda\epsilon + \frac{1}{N} \sum_{n \in \mathcal{N}} \left(s_n + \sum_{i \in \Aplus} \hat{u}_{ni}\, \rho_{ni} \right) \label{eq:obj_milp} \\
    \text{s.t.} \quad & \eqref{eq:c1}-\eqref{eq:c3}, \nonumber \\
    & 0 \le \rho_{ni} \le \lambda^u y_{ni}, && \forall n \in \mathcal{N},\ \forall i \in \Aplus, \label{eq:mc1} \\
    & \lambda^u (y_{ni} - 1) + \lambda \le \rho_{ni} \le \lambda, && \forall n \in \mathcal{N},\ \forall i \in \Aplus, \label{eq:mc2} \\
    & \sum_{i \in \Aplus} \rho_{ni} = \lambda, && \forall n \in \mathcal{N}. \label{eq:vi}
\end{align}

Constraints~\eqref{eq:mc1}--\eqref{eq:mc2} are the McCormick envelope of \(\rho_{ni}=\lambda y_{ni}\) over \(\lambda\in[0,\lambda^u]\), \(y_{ni}\in[0,1]\). Constraint~\eqref{eq:vi} is obtained by multiplying the assignment constraint~\eqref{eq:c2} by \(\lambda\); it is not implied by the McCormick system alone and serves two purposes. First, it strengthens the LP relaxation by linking the auxiliary variables \(\rho_{ni}\) across products for each customer. Second, it makes the linearization exact in value: because continuous \(\boldsymbol{y}\) admits fractional points at which the McCormick system alone permits \(\sum_{i\in\Aplus}\hat u_{ni}\rho_{ni}\) to exceed the bilinear value, [DRO-MILP] without~\eqref{eq:vi} is a strict relaxation in general.

We now verify that [DRO-MILP] is exact. Fix any feasible solution of [DRO-MILP] with product line \(\mathcal{S}=\{i\in\mathcal{A}:x_i=1\}\). Constraints~\eqref{eq:mc1} and~\eqref{eq:c3} give \(\rho_{ni}\le\lambda^u y_{ni}\le\lambda^u x_i\), so \(\rho_{ni}=0\) whenever \(x_i=0\). If \(\lambda>0\), constraint~\eqref{eq:vi} makes \(\{\rho_{ni}/\lambda\}_{i}\) a convex-combination weight vector supported on \(\bar{\mathcal{S}}\), whence
\[
\sum_{i \in \Aplus} \hat{u}_{ni}\,\rho_{ni} \;\le\; \lambda \max_{i \in \bar{\mathcal{S}}} \hat{u}_{ni};
\]
if \(\lambda=0\), \eqref{eq:vi} forces \(\boldsymbol{\rho}_n=\boldsymbol{0}\) and the same bound holds trivially. Combined with \(s_n \le \min_{i\in\bar{\mathcal{S}}}(r_i-\lambda\hat u_{ni})\) from~\eqref{eq:c1}, every feasible objective value is at most \(\max_{\mathcal{S}\in\mathcal{L},\,\lambda\in[0,\lambda^u]}\{-\lambda\epsilon+\lambda W_N(\mathcal{S})+C(\mathcal{S},\lambda)\}=V^{\star}_{N,\epsilon}\), using Lemma~\ref{prop:lambda_bounds} for the restriction of \(\lambda\) to \([0,\lambda^u]\). Conversely, the value \(V^{\star}_{N,\epsilon}\) is attained by taking an optimal \((\mathcal{S}^{\star},\lambda^{\star})\), setting binary \(y_{ni}\) to select a utility-maximizing alternative in \(\mathcal{S}^{\star}\cup\{0\}\) for each \(n\), and setting \(\rho_{ni}=\lambda^{\star}y_{ni}\), which satisfies \eqref{eq:mc1}--\eqref{eq:vi}. Hence the optimal value of [DRO-MILP] equals \(V^{\star}_{N,\epsilon}\).

With the $\lambda$ range from Lemma~\ref{prop:lambda_bounds}, the Big-M constants can be set as
\[
\bar{C}_n = \max_{0\le\lambda\le\lambda^u}\left\{\max_{i\in\Aplus}\{r_i - \lambda\hat{u}_{ni}\} -\min_{j\in\Aplus}\{r_j - \lambda\hat{u}_{nj}\}\right\}
\]
for each $n\in\mathcal{N}$.

The benchmark implementation used to generate Table~\ref{tab:main_results} solves an algebraically equivalent variant that folds the welfare term into the epigraph constraints, replacing~\eqref{eq:c1} by \(s_n \le r_i + \sum_{j\in\Aplus}\hat u_{nj}\rho_{nj} - \lambda\hat u_{ni} + \bar{C}_n(1-x_i)\) so that the objective carries only \(\frac{1}{N}\sum_n s_n\). The epigraph variables remain free. Apart from this algebraic rearrangement, the implementation uses the continuous assignment variables, the McCormick system~\eqref{eq:mc1}--\eqref{eq:mc2}, and the aggregation constraint~\eqref{eq:vi} exactly as displayed above. No first-choice ordering inequalities or nonempty-line restriction are imposed.

\section{SBMB Inner Solver Performance}\label{appendix:sbmb_benchmark}

To isolate the contribution of SBMB to the speedup in Table~\ref{tab:main_results}, we fix $\lambda$ and compare SBMB with a direct mixed-integer linear formulation of the same subproblem. We use the product-line variables $\boldsymbol{x}\in\mathcal X$ from Section~\ref{appendix:milp}, continuous assignment variables $\boldsymbol y$, and free epigraph variables $\boldsymbol s$:
\begin{align}
    \text{[MILP]}\qquad
    \max_{\boldsymbol{x}\in\mathcal X,\,\boldsymbol y\ge0,\,\boldsymbol s\in\mathbb R^N}\quad
    & \frac{1}{N}\sum_{n\in\mathcal N}
    \left(s_n+\lambda\sum_{i\in\Aplus}\hat u_{ni}y_{ni}\right) \label{eq:fixed_milp_obj}\\
    \text{s.t.}\quad
    & s_n\le r_i-\lambda\hat u_{ni}+\bar C_n(\lambda)(1-x_i),
    && \forall n\in\mathcal N,\ i\in\Aplus, \label{eq:fixed_milp_epi}\\
    & \sum_{i\in\Aplus}y_{ni}=1,
    && \forall n\in\mathcal N, \label{eq:fixed_milp_assign}\\
    & y_{ni}\le x_i,
    && \forall n\in\mathcal N,\ i\in\Aplus. \label{eq:fixed_milp_link}
\end{align}
Here, $\bar C_n(\lambda)$ deactivates~\eqref{eq:fixed_milp_epi} for alternatives outside the product line; the implementation uses
\[
    \bar C_n(\lambda)
    =1+\max_{i\in\Aplus}\{r_i-\lambda\hat u_{ni}\}
       -\min_{i\in\Aplus}\{r_i-\lambda\hat u_{ni}\}.
\]
At optimality, \eqref{eq:fixed_milp_epi} recovers the customer-level minimum in $C(\mathcal S,\lambda)$ and \eqref{eq:fixed_milp_assign}--\eqref{eq:fixed_milp_link} recover the maximum in $W_N(\mathcal S)$. Thus [MILP] and SBMB both compute $f(\lambda)=\max_{\mathcal S\in\mathcal L}G(\mathcal S,\lambda)$. Fixing $\lambda$ makes $-\lambda\epsilon$ constant and eliminates the need for McCormick variables.

We use the complete grid and all 10 seeds from Section~\ref{sec:exp_setup}: $A,N\in\{50,100,200\}$, $p\in\{5,10\}$, and $\epsilon\in\{0.25\epsdag,0.75\epsdag\}$. We solve each instance at $\lambda/\lambda^u\in\{0.25,0.50,1.00\}$, for 1,080 paired problems.

[MILP] uses Gurobi~11.0.3 with one thread. Both methods use a relative optimality tolerance of $10^{-4}$ and no time limit. The timing includes [MILP] model construction, which would recur at every call within a parametric search, but excludes shared preprocessing and one-time compilation. Table~\ref{tab:sbmb_vs_milp} reports mean times over the 10 seeds.

\begin{table}[htbp]
\centering
\caption{\hspace{0.5em}Mean solution times (seconds) for the fixed-$\lambda$ subproblem.}\label{tab:sbmb_vs_milp}
\scriptsize
\renewcommand{\arraystretch}{0.90}
\setlength{\tabcolsep}{3.5pt}
\begin{tabular*}{\textwidth}{@{\extracolsep{\fill}}cccc rr rr rr@{}}
\toprule
& & & & \multicolumn{2}{c}{$\lambda/\lambda^u=0.25$} & \multicolumn{2}{c}{$\lambda/\lambda^u=0.50$} & \multicolumn{2}{c}{$\lambda/\lambda^u=1.00$} \\
\cmidrule(lr){5-6} \cmidrule(lr){7-8} \cmidrule(lr){9-10}
$p$ & $A$ & $N$ & $\epsilon$ & MILP & SBMB & MILP & SBMB & MILP & SBMB \\
\midrule
5 & 50 & 50 & H & 0.17 & $<0.01$ & 0.19 & $<0.01$ & 0.19 & $<0.01$ \\
 &  &  & L & 0.17 & $<0.01$ & 0.20 & $<0.01$ & 0.21 & $<0.01$ \\
 &  & 100 & H & 0.33 & $<0.01$ & 0.35 & $<0.01$ & 0.38 & $<0.01$ \\
 &  &  & L & 0.34 & $<0.01$ & 0.37 & $<0.01$ & 0.41 & $<0.01$ \\
 &  & 200 & H & 0.72 & $<0.01$ & 0.69 & $<0.01$ & 0.77 & $<0.01$ \\
 &  &  & L & 0.77 & $<0.01$ & 0.78 & $<0.01$ & 0.85 & 0.01 \\
\cmidrule(l){2-10}
 & 100 & 50 & H & 0.32 & $<0.01$ & 0.34 & 0.01 & 0.39 & 0.02 \\
 &  &  & L & 0.36 & $<0.01$ & 0.44 & 0.03 & 0.58 & 0.03 \\
 &  & 100 & H & 0.69 & $<0.01$ & 0.68 & 0.02 & 0.78 & 0.03 \\
 &  &  & L & 0.79 & 0.02 & 0.85 & 0.05 & 1.19 & 0.08 \\
 &  & 200 & H & 1.35 & $<0.01$ & 1.47 & 0.03 & 1.73 & 0.05 \\
 &  &  & L & 1.59 & 0.03 & 2.00 & 0.11 & 3.20 & 0.22 \\
\cmidrule(l){2-10}
 & 200 & 50 & H & 0.71 & 0.01 & 0.68 & 0.02 & 0.86 & 0.07 \\
 &  &  & L & 0.80 & 0.03 & 1.12 & 0.17 & 2.18 & 0.56 \\
 &  & 100 & H & 1.39 & 0.01 & 1.50 & 0.03 & 1.89 & 0.17 \\
 &  &  & L & 1.63 & 0.08 & 2.58 & 0.42 & 7.26 & 1.71 \\
 &  & 200 & H & 2.75 & 0.02 & 3.06 & 0.05 & 4.04 & 0.32 \\
 &  &  & L & 3.36 & 0.13 & 5.69 & 0.85 & 16.82 & 4.03 \\
\midrule
10 & 50 & 50 & H & 0.17 & $<0.01$ & 0.18 & $<0.01$ & 0.19 & $<0.01$ \\
 &  &  & L & 0.17 & $<0.01$ & 0.20 & $<0.01$ & 0.21 & $<0.01$ \\
 &  & 100 & H & 0.33 & $<0.01$ & 0.33 & $<0.01$ & 0.37 & $<0.01$ \\
 &  &  & L & 0.35 & $<0.01$ & 0.37 & $<0.01$ & 0.40 & $<0.01$ \\
 &  & 200 & H & 0.71 & $<0.01$ & 0.67 & $<0.01$ & 0.74 & 0.01 \\
 &  &  & L & 0.76 & $<0.01$ & 0.77 & 0.01 & 0.83 & 0.01 \\
\cmidrule(l){2-10}
 & 100 & 50 & H & 0.32 & $<0.01$ & 0.35 & 0.01 & 0.39 & 0.01 \\
 &  &  & L & 0.37 & 0.01 & 0.47 & 0.02 & 0.83 & 0.02 \\
 &  & 100 & H & 0.70 & $<0.01$ & 0.68 & 0.02 & 0.79 & 0.04 \\
 &  &  & L & 0.79 & 0.02 & 0.97 & 0.04 & 3.40 & 0.11 \\
 &  & 200 & H & 1.35 & 0.01 & 1.46 & 0.03 & 1.74 & 0.13 \\
 &  &  & L & 1.57 & 0.06 & 2.34 & 0.24 & 6.98 & 0.82 \\
\cmidrule(l){2-10}
 & 200 & 50 & H & 0.70 & 0.01 & 0.68 & 0.03 & 0.94 & 0.10 \\
 &  &  & L & 0.83 & 0.05 & 1.75 & 0.25 & 22.16 & 0.65 \\
 &  & 100 & H & 1.35 & 0.02 & 1.48 & 0.15 & 2.98 & 1.82 \\
 &  &  & L & 1.79 & 0.71 & 8.79 & 4.73 & 146.06 & 16.56 \\
 &  & 200 & H & 2.71 & 0.04 & 3.00 & 0.30 & 6.59 & 6.49 \\
 &  &  & L & 3.55 & 1.53 & 20.57 & 22.60 & 325.30 & 131.13 \\
\bottomrule
\end{tabular*}
\end{table}

[MILP] and SBMB reach the tolerance on all 1,080 problems with identical objective values. SBMB is faster in 1,065 paired runs and in 107 of the 108 comparisons of configuration means. The exception is $A=N=200$, $p=10$, $\epsilon=0.25\epsdag$, and $\lambda/\lambda^u=0.50$: 22.60 seconds for SBMB versus 20.57 seconds for [MILP]. The ratios of total [MILP] time to total SBMB time are $12.70$, $2.25$, and $3.42$ at $\lambda/\lambda^u=0.25$, $0.50$, and $1.00$, respectively.

\section{Nestedness and Welfare--Revenue Alignment under Capacity Constraints}\label{appendix:capacity_separation}

A cardinality constraint can overturn assortment nestedness, but it does not remove the welfare term from the robust revenue of a fixed product line in [DRO-FC].

\subsection{Nestedness Need Not Survive a Cardinality Constraint}\label{appendix:capacity_separation:parametric_collapse}

Corollary~3.3 of \citet{rusmevichientong2012robust} and Proposition~6 of \citet{desir2024robust} guarantee nested unconstrained optimizers for nested uncertainty sets. The following examples show why neither conclusion extends to $\mathcal L_p=\{\mathcal S\subseteq\mathcal A:|\mathcal S|\le p\}$.

\subsubsection{Robust MNL}\label{appendix:capacity_separation:rt_collapse}

In the model of \citet[\S3]{rusmevichientong2012robust}, $v_0=1$ is the no-purchase attraction and
\[
R(\mathcal S,v)=\frac{\sum_{i\in\mathcal S}r_iv_i}{v_0+\sum_{i\in\mathcal S}v_i}.
\]
For a box $V$, let $\mathcal S_p^*(V)$ maximize $\min_{v\in V}R(\mathcal S,v)$ over $\mathcal L_p$. We take $A=4$, $p=2$, $r=(10,8,6,5)$, and the two boxes in Table~\ref{tab:capacity_boxes}. The inclusion $V'\subset V$ is strict.

\begin{table}[htbp]
\centering
\caption{Nested rectangular uncertainty sets in the robust-MNL counterexample.}
\label{tab:capacity_boxes}
\small
\setlength{\tabcolsep}{10pt}
\renewcommand{\arraystretch}{1.1}
\begin{tabular}{@{}c c c l@{}}
\toprule
Coordinate & $V$ box & $V'$ box & Tightening \\
\midrule
$v_1$ & $[0.1,\,5.0]$ & $[0.1,\,0.3]$ & shrunk to low end \\
$v_2$ & $[0.1,\,5.0]$ & $[4.0,\,5.0]$ & shrunk to high end \\
$v_3$ & $[0.5,\,1.0]$ & $[0.5,\,1.0]$ & unchanged \\
$v_4$ & $[3.0,\,5.0]$ & $[3.0,\,5.0]$ & unchanged \\
\bottomrule
\end{tabular}
\end{table}

By Lemma~3.7 of \citet{rusmevichientong2012robust}, it suffices to evaluate the extreme points of each box. Table~\ref{tab:capacity_per_assortment} reports the resulting worst-case revenue for every nonempty line in $\mathcal L_2$.

\begin{table}[htbp]
\centering
\caption{Worst-case MNL revenue for every nonempty line with $|\mathcal S|\le2$; bold starred entries are optimal.}
\label{tab:capacity_per_assortment}
\small
\setlength{\tabcolsep}{12pt}
\renewcommand{\arraystretch}{1.1}
\begin{tabular}{@{}l r r @{\qquad} l r r@{}}
\toprule
$\mathcal{S}$ & over $V$ & over $V'$ & $\mathcal{S}$ & over $V$ & over $V'$ \\
\midrule
$\{1\}$   & $0.9091$ & $0.9091$ & $\{1,3\}$ & $2.5000$ & $2.5000$ \\
$\{2\}$   & $0.7273$ & $6.4000$ & $\{1,4\}$ & $3.9024$ & $3.9024$ \\
$\{3\}$   & $2.0000$ & $2.0000$ & $\{2,3\}$ & $2.3750$ & $6.3333$ \\
$\{4\}$   & $3.7500$ & $3.7500$ & $\{2,4\}$ & $3.8537$ & $5.7000$ \\
$\{1,2\}$ & $1.5000$ & $\mathbf{6.4706}^{*}$ & $\{3,4\}$ & $\mathbf{4.0000}^{*}$ & $4.0000$ \\
\bottomrule
\end{tabular}
\end{table}

The unique constrained optima $\mathcal S_p^*(V)=\{3,4\}$ and $\mathcal S_p^*(V')=\{1,2\}$ are disjoint. Without the capacity limit, the optima are $\{1,2,3,4\}$ under $V$ and $\{1,2\}$ under $V'$, consistent with Corollary~3.3.

\subsubsection{Robust Markov Chain}\label{appendix:capacity_separation:desir_collapse}

In the model of \citet[\S2.1]{desir2024robust}, $\alpha_i>0$ is the arrival probability at product $i$, and $\rho_{ij}$ is the probability of moving from an unavailable product $i$ to product $j$. We use $\alpha$ for the arrival vector, which they denote by $\lambda$, to distinguish it from our dual multiplier. For a fixed line $\mathcal S$,
\[
R^{\mathrm{MC}}(\mathcal S;\alpha,\rho)=\alpha^\top v(\rho),
\qquad
v_i(\rho)=
\begin{cases}
r_i, & i\in\mathcal S,\\
\sum_{j\in\mathcal A}\rho_{ij}v_j(\rho), & i\notin\mathcal S.
\end{cases}
\]
The outside transition is $\rho_{i0}=1-\sum_{j\in\mathcal A}\rho_{ij}$. We impose $\rho_{ii}=0$ and $\sum_j\rho_{ij}\le0.99$, so every admissible row satisfies the strict condition in Assumption~2 of \citet{desir2024robust}.

Again take $A=4$, $p=2$, and $r=(10,8,6,5)$. Table~\ref{tab:capacity_mc_boxes} defines two row-wise uncertainty sets. The arrival coordinates lie in the displayed intervals and satisfy $\sum_i\alpha_i\le1$. Transition rows vary independently; each off-diagonal coordinate has upper bound one, subject to the common row-mass bound $0.99$. We write $\boldsymbol l_i$ for the lower-bound vector of row $i$. The smaller set raises selected arrival and transition lower bounds, so $U^{\mathrm{MC}\prime}\subset U^{\mathrm{MC}}$.

\begin{table}[htbp]
\centering
\caption{Nested row-wise uncertainty sets in the robust Markov-chain counterexample. Transition lower bounds are listed in product order $(1,2,3,4)$.}
\label{tab:capacity_mc_boxes}
\small
\setlength{\tabcolsep}{8pt}
\renewcommand{\arraystretch}{1.1}
\begin{tabular}{@{}c c c@{}}
\toprule
Coordinate & $U^{\mathrm{MC}}$ & $U^{\mathrm{MC}\prime}$ \\
\midrule
$\alpha_1$ & $[0.05,\,0.20]$ & $[0.15,\,0.20]$ \\
$\alpha_2$ & $[0.05,\,0.20]$ & $[0.15,\,0.20]$ \\
$\alpha_3$ & $[0.30,\,0.40]$ & $[0.30,\,0.40]$ \\
$\alpha_4$ & $[0.30,\,0.40]$ & $[0.30,\,0.40]$ \\
$\boldsymbol l_{1},\boldsymbol l_{2}$ & $\boldsymbol0$ & $\boldsymbol0$ \\
$\boldsymbol l_{3}$ & $\boldsymbol0$ & $(0.524097,0.298207,0,0.027696)$ \\
$\boldsymbol l_{4}$ & $\boldsymbol0$ & $(0.496549,0.281499,0.071952,0)$ \\
\bottomrule
\end{tabular}
\end{table}

\begin{aplemma}[Lower-Bound Worst Case]\label{lem:mc_inner_min}
For either uncertainty set in Table~\ref{tab:capacity_mc_boxes} and any fixed $\mathcal S$, the balance equations have a unique nonnegative solution, and
\[
\min_{(\alpha,\rho)}R^{\mathrm{MC}}(\mathcal S;\alpha,\rho)
= (\boldsymbol\alpha^l)^\top v(\boldsymbol l),
\]
where $\boldsymbol\alpha^l$ and $\boldsymbol l$ are the corresponding arrival and transition lower bounds.
\end{aplemma}

\textit{Proof.}
For fixed $\rho$, let $T_\rho$ denote the mapping on the right-hand side of the balance equations. Its contraction modulus under the sup norm is at most $0.99$, so $v(\rho)$ is unique and nonnegative. The mapping is coordinatewise nondecreasing in both its argument and $\rho$. Iteration from $\boldsymbol0$ therefore gives $v(\boldsymbol l)\le v(\rho)$. Since $\alpha\ge\boldsymbol\alpha^l$ and $v(\rho)\ge\boldsymbol0$,
\[
\alpha^\top v(\rho)\ge(\boldsymbol\alpha^l)^\top v(\rho)
\ge(\boldsymbol\alpha^l)^\top v(\boldsymbol l).
\]
The lower bounds are feasible, so equality is attained at $(\boldsymbol\alpha^l,\boldsymbol l)$. \hfill$\square$

Table~\ref{tab:capacity_mc_per_assortment} reports these exact worst-case values for every nonempty line in $\mathcal L_2$.

\begin{table}[htbp]
\centering
\caption{Worst-case Markov-chain revenue for every nonempty line with $|\mathcal S|\le2$; bold starred entries are optimal.}
\label{tab:capacity_mc_per_assortment}
\small
\setlength{\tabcolsep}{12pt}
\renewcommand{\arraystretch}{1.1}
\begin{tabular}{@{}l r r @{\qquad} l r r@{}}
\toprule
$\mathcal{S}$ & over $U^{\mathrm{MC}}$ & over $U^{\mathrm{MC}\prime}$ & $\mathcal{S}$ & over $U^{\mathrm{MC}}$ & over $U^{\mathrm{MC}\prime}$ \\
\midrule
$\{1\}$   & $0.5000$ & $4.7227$ & $\{1,3\}$ & $2.3000$ & $4.9192$ \\
$\{2\}$   & $0.4000$ & $2.6644$ & $\{1,4\}$ & $2.0000$ & $4.6138$ \\
$\{3\}$   & $1.8000$ & $1.9295$ & $\{2,3\}$ & $2.2000$ & $3.8051$ \\
$\{4\}$   & $1.5000$ & $1.5415$ & $\{2,4\}$ & $1.9000$ & $3.4572$ \\
$\{1,2\}$ & $0.9000$ & $\mathbf{7.3872}^{*}$ & $\{3,4\}$ & $\mathbf{3.3000}^{*}$ & $3.3000$ \\
\bottomrule
\end{tabular}
\end{table}

The unique constrained optima are $\{3,4\}$ under $U^{\mathrm{MC}}$ and $\{1,2\}$ under $U^{\mathrm{MC}\prime}$. Under $U^{\mathrm{MC}}$, zero product-to-product transitions are feasible, so the two lines earn $0.30(6)+0.30(5)=3.30$ and $0.05(10)+0.05(8)=0.90$, respectively. Under $U^{\mathrm{MC}\prime}$, higher arrival floors at products $1$ and $2$ and positive transition floors from products $3$ and $4$ raise the value of $\{1,2\}$ to $7.3872$; $\{3,4\}$ still earns $3.30$. The optima are thus disjoint, whereas the unconstrained optima $\{1,2,3,4\}$ and $\{1,2\}$ are nested, consistent with Proposition~6 of \citet{desir2024robust}.

In both examples, capacity forces a change in the optimal line to take the form of replacement rather than expansion. A change in the uncertainty set can reverse the ranking of the feasible pairs.

\subsection{The Welfare Term under Capacity Constraints}\label{appendix:capacity_separation:our_invariance}

In [DRO-FC], welfare enters the robust objective before product lines are compared. Proposition~\ref{prop:dual_worst_case_revenue}, specialized to the empirical distribution, gives, for every fixed $\mathcal S$,
\[
R^{\mathrm{FC}}_{N,\epsilon}(\mathcal{S}) = \max_{\lambda \ge 0}\ \Bigl\{-\lambda\,\epsilon\ +\ \lambda\,W_N(\mathcal{S})\ +\ \frac{1}{N}\sum_{n\in\mathcal{N}}\min_{i\in\bar{\mathcal{S}}}\bigl(r_i - \lambda\,\hat{u}_{ni}\bigr)\Bigr\}.
\]
This identity is established before optimizing over the feasible family. Replacing an unconstrained family by $\mathcal L_p$, or by any other operationally constrained family, changes which product lines are compared but not the term $+\lambda W_N(\mathcal S)$ for any retained line. The critical radius $\epsilon^\dagger(\widehat Q_N)$ is then computed over that family, and Corollary~\ref{corollary:welfare_lb} continues to provide the corresponding welfare guarantee.

\begin{figure}[htbp]
\centering
\includegraphics[width=\textwidth]{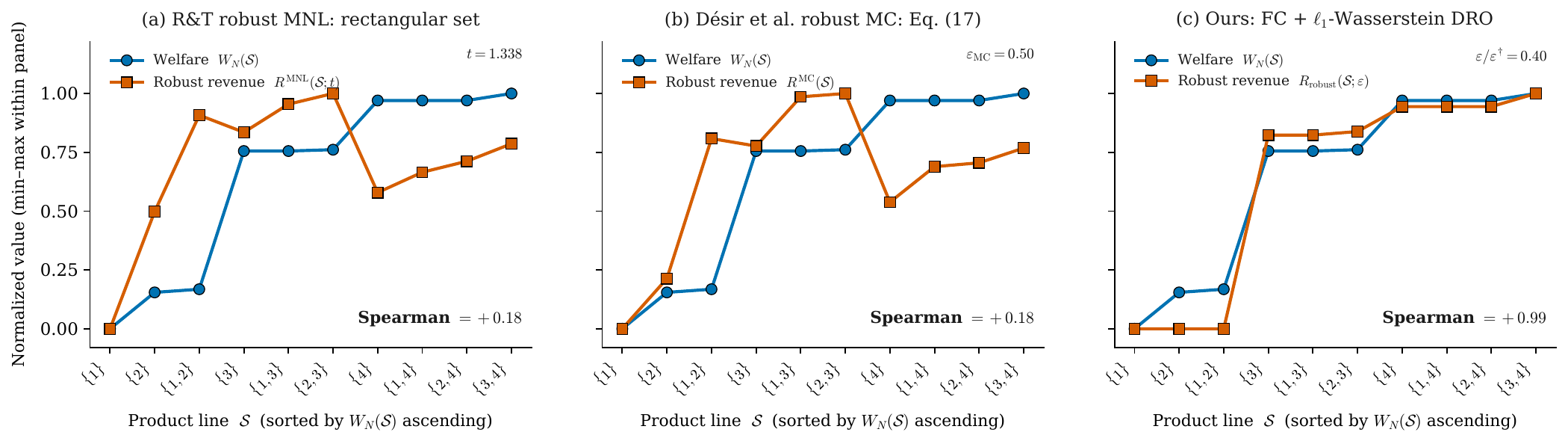}
\caption{Welfare and robust revenue across the ten nonempty lines with $|\mathcal S|\le2$, ordered by $W_N(\mathcal S)$. Values are min--max normalized within each panel; panel annotations report the ambiguity level and Spearman rank correlation.}
\label{fig:alignment_scatter}
\end{figure}

The last term in the display also varies with $\mathcal S$, so the strength of the revenue--welfare association may differ across instances. Figure~\ref{fig:alignment_scatter} shows this association on a common four-product instance with $r=(10,8,6,5)$, $N=15$, and $|\mathcal S|\le2$. Product and outside-option utilities are drawn independently with standard deviation $0.5$ around $(1.0,1.5,3.0,3.5)$ and $1.0$, respectively, using seed $20250501$. We evaluate all ten nonempty feasible lines and order them by $W_N(\mathcal S)$.

For the robust-MNL benchmark, $\bar v_i=\exp(N^{-1}\sum_n\hat u_{ni})$, $\bar v_0$ is fixed, and $v_i\in[\bar v_i-t,\bar v_i+t]$ with $t=1.338$. For the robust Markov-chain benchmark, modal arrival and transition probabilities are obtained from the same mean utilities, including the outside option in each normalization; Equation~(17) of \citet{desir2024robust} is then applied with $\epsilon_{\mathrm{MC}}=0.50$. We evaluate [DRO-FC] at $\epsilon=0.996=0.40\epsdag$.

At these ambiguity levels, robust-MNL and robust Markov-chain revenues are nonmonotone in the welfare ordering, with Spearman correlations of $0.18$. [DRO-FC] tracks the ordering closely (correlation $0.99$) and selects the welfare-maximizing line $\{3,4\}$.

\section{Proofs}\label{appendix:proofs}

    \subsection{Proof of Proposition~\ref{prop:dual_worst_case_revenue}}
    
    For a fixed product line $\mathcal S$, let $g_{\mathcal S}(\boldsymbol\xi)=\min\{r_i:i\in\arg\max_{j\in\bar{\mathcal S}}U_j\}$ be its realized revenue under the lower-revenue tie rule in Section~\ref{sec:MF}. This bounded, lower-semicontinuous function satisfies the conditions for Wasserstein strong duality~\citep{esfahani2018data,Blanchet2019Quantifying}:
    \begin{equation}\label{eq:general_wasserstein_step}
    R_\epsilon(\mathcal S;Q)
    =\sup_{\lambda\ge0}\left\{
    -\lambda\epsilon
    +\mathbb E_Q\!\left[
    \inf_{\boldsymbol\xi\in\mathbb R^{A+1}}
    \bigl\{g_{\mathcal S}(\boldsymbol\xi)
    +\lambda\|\boldsymbol\xi-\boldsymbol\xi'\|_1\bigr\}
    \right]
    \right\},
    \end{equation}
    where $\boldsymbol\xi'$ is distributed according to $Q$.
    
    For $i\in\bar{\mathcal S}$, the closure of its choice region is $\mathcal C_i(\mathcal S)=\{\boldsymbol\xi\in\mathbb R^{A+1}:U_i\ge U_j\ \forall j\in\bar{\mathcal S}\}$. These cells cover $\mathbb R^{A+1}$; on overlaps the tie rule takes the lowest revenue. Let $m=\max_{j\in\bar{\mathcal S}}U'_j$. Raising coordinate $i$ to $m$ reaches $\mathcal C_i(\mathcal S)$ at $\ell_1$ distance $m-U'_i$. For the reverse bound, fix $\boldsymbol\xi\in\mathcal C_i(\mathcal S)$ and choose $j^\star\in\arg\max_{j\in\bar{\mathcal S}} U'_j$. If $j^\star=i$, the bound is immediate. Otherwise, $U_i\ge U_{j^\star}$ implies
    \[
    m-U'_i
    =U'_{j^\star}-U'_i
    \le |U'_{j^\star}-U_{j^\star}|+|U_i-U'_i|
    \le\|\boldsymbol\xi-\boldsymbol\xi'\|_1.
    \]
    Taking the infimum over $\boldsymbol\xi\in\mathcal C_i(\mathcal S)$ gives the distance. Therefore
    \begin{align*}
    &\inf_{\boldsymbol\xi\in\mathbb R^{A+1}}
    \{g_{\mathcal S}(\boldsymbol\xi)
    +\lambda\|\boldsymbol\xi-\boldsymbol\xi'\|_1\}\\
    &\qquad=\min_{i\in\bar{\mathcal S}}
    \left\{r_i+\lambda\left(
    \max_{j\in\bar{\mathcal S}}U'_j-U'_i
    \right)\right\}\\
    &\qquad=\lambda\max_{j\in\bar{\mathcal S}}U'_j
    +\min_{i\in\bar{\mathcal S}}\{r_i-\lambda U'_i\}.
    \end{align*}
    Substituting this identity into~\eqref{eq:general_wasserstein_step} and using
    $\min_i\{r_i-\lambda U_i\}=-\lambda\max_i\{U_i-r_i/\lambda\}$ for $\lambda>0$ gives~\eqref{e4}. The expression at $\lambda=0$ equals zero because $r_0=0$ and is obtained as the continuous boundary value. The combined integrand is continuous in $\lambda$ and lies in $[0,\bar r]$: choosing a utility-maximizing alternative in the inner minimum gives the upper bound. Dominated convergence therefore gives continuity of the dual objective on $[0,\infty)$, while the term $-\lambda\epsilon$ forces it to $-\infty$ as $\lambda\to\infty$ for $\epsilon>0$. Hence, the supremum is attained. \qed
    
    \subsection{Proof of Corollary~\ref{corollary:welfare_lb}}
    
    For fixed $\mathcal S$, write the objective in~\eqref{e4} as
    \[
    \psi_{\mathcal S}(\lambda;\epsilon,Q)
    =-\lambda\epsilon+\lambda W(\mathcal S;Q)
    +\mathbb E_Q\!\left[
    \min_{i\in\bar{\mathcal S}}\{r_i-\lambda U_i\}
    \right].
    \]
    This objective is concave and equals zero at $\lambda=0$. Since the outside option uniquely minimizes $r_i$, dominated convergence gives $\partial_+\psi_{\mathcal S}(0;\epsilon,Q)=W(\mathcal S;Q)-\mathbb E_Q[U_0]-\epsilon$. If the derivative is nonpositive, concavity implies that the maximum is zero; if it is positive, the objective increases to the right of zero and the maximum is positive.
    
    Because $\epsilon<\epsilon^\dagger(Q)$, a welfare-maximizing product line has incremental welfare strictly greater than $\epsilon$ and hence positive robust revenue. Every robust-optimal product line must therefore have positive robust revenue, which requires
    $W(\mathcal S^*_{\epsilon}(Q);Q)-\mathbb E_Q[U_0]>\epsilon$. The weak inequality stated in the corollary follows. \qed
    
    \subsection{Proof of Proposition~\ref{prop:revenue_ordered_mnl}}
    
    Fix $\lambda>0$. The fixed-$\lambda$ expression for [DRO-MNL] is maximized by maximizing
    \[
    F_\lambda(\mathcal S)
    =\frac{v_0+\sum_{i\in\mathcal S}v_i}
    {v_0+\sum_{i\in\mathcal S}v_i e^{-r_i/\lambda}}.
    \]
    Let $\mathcal S_\lambda$ be optimal and denote its numerator and denominator by $A_\lambda$ and $B_\lambda$. Adding an excluded product $i$ cannot improve the ratio, which implies
    $e^{-r_i/\lambda}\ge B_\lambda/A_\lambda$. Removing an included product cannot improve it, which implies
    $e^{-r_i/\lambda}\le B_\lambda/A_\lambda$. Since $e^{-r/\lambda}$ decreases in $r$, an optimal product line can be chosen as a revenue-ordered prefix for every fixed $\lambda>0$. At $\lambda=0$, the continuous boundary value is zero for every product line, so a prefix can again be selected. Maximizing jointly over the finite set of prefixes and $\lambda$ proves the result. \qed
    
    \subsection{Proof of Proposition~\ref{prop:monotone_assortment_epsilon}}
    
    Fix $\lambda>0$ and write $t=1/\lambda$. For the prefix $\mathcal S_k=\{1,\ldots,k\}$, define
    \[
    \rho_i(t)=1-e^{-tr_i},\qquad
    A_k=v_0+\sum_{i=1}^k v_i,\qquad
    \Psi_k(t)=\frac{\sum_{i=1}^k v_i\rho_i(t)}{A_k},
    \]
    with $\Psi_0(t)=0$. Maximizing the fixed-$\lambda$ [DRO-MNL] objective over prefixes is equivalent to maximizing $\Psi_k(t)$.
    
    The sequence $\{\Psi_k(t)\}_{k=0}^A$ is single-peaked in $k$. Indeed,
    \begin{equation}\label{eq:mnl_adjacent_addition}
    \Psi_{k+1}(t)\ge\Psi_k(t)
    \quad\Longleftrightarrow\quad
    \rho_{k+1}(t)\ge\Psi_k(t).
    \end{equation}
    If the inequality fails at $k$, then $\Psi_{k+1}(t)$ is a strict weighted average of $\Psi_k(t)$ and $\rho_{k+1}(t)$. Since $r_{k+2}\le r_{k+1}$, we have
    $\rho_{k+2}(t)\le\rho_{k+1}(t)<\Psi_{k+1}(t)$, so the next addition also fails; the same argument applies inductively. The largest maximizing prefix is therefore the last index whose adjacent addition is weakly beneficial.
    
    Each adjacent addition also has a single-crossing property in $t$. Rearranging~\eqref{eq:mnl_adjacent_addition} gives
    \[
    H_k(t):=
    \frac{v_0e^{tr_{k+1}}+
    \sum_{i=1}^k v_i e^{-t(r_i-r_{k+1})}}{A_k}
    \ge1.
    \]
    The function $H_k$ is strictly convex, $H_k(0)=1$, and $H_k(t)\to\infty$ as $t\to\infty$. It follows that $\{t>0:H_k(t)\ge1\}$ is an upper interval. Combined with single-peakedness, the largest fixed-$\lambda$ optimal prefix is nondecreasing in $t=1/\lambda$; the conclusion extends to $\lambda=0$ by selecting the largest prefix.
    
    It remains to establish the monotonicity of the optimal multiplier. Let
    \[
    J_k(\lambda)
    =\lambda\log
    \frac{v_0+\sum_{i=1}^k v_i}
    {v_0+\sum_{i=1}^k v_i e^{-r_i/\lambda}},
    \]
    with its continuous boundary value at $\lambda=0$, so the objective associated with pair $(k,\lambda)$ is $J_k(\lambda)-\epsilon\lambda$. Fix $\epsilon_1<\epsilon_2$ and write $(k_q,\lambda_q)=(k^*(\epsilon_q),\lambda^*_{\epsilon_q})$ for $q=1,2$. Optimality at the two radii gives
    \begin{align*}
    J_{k_1}(\lambda_1)-\epsilon_1\lambda_1
    &\ge J_{k_2}(\lambda_2)-\epsilon_1\lambda_2,\\
    J_{k_2}(\lambda_2)-\epsilon_2\lambda_2
    &\ge J_{k_1}(\lambda_1)-\epsilon_2\lambda_1.
    \end{align*}
    Adding these inequalities yields
    $(\epsilon_2-\epsilon_1)(\lambda_1-\lambda_2)\ge0$, and hence $\lambda_1\ge\lambda_2$. Thus, $\lambda^*_{\epsilon}$ is nonincreasing in $\epsilon$, so $1/\lambda^*_{\epsilon}$ is nondecreasing, with $1/0$ interpreted as $+\infty$. Applying the fixed-multiplier result and the stated largest-prefix tie rule gives
    $k^*(\epsilon_1)\le k^*(\epsilon_2)$ whenever $\epsilon_1\le\epsilon_2$. \qed
    
    \subsection{Proof of Proposition~\ref{prop:mle}}
    
    With $v_0=1$, the multinomial log-likelihood is $\ell_T(v)=\sum_{i=1}^A\hat T_i\log v_i-T\log(1+\sum_{i=1}^Av_i)$ for $v\in(0,\infty)^A$. It extends continuously to $v_i=0$ when $\hat T_i=0$. If $\hat T_0>0$, an interior critical point satisfies $\hat T_k/v_k=T/(1+\sum_{i=1}^Av_i)$ for each $k=1,\ldots,A$.
    Summing $\hat{T}_k=v_k T/(1+\sum_i v_i)$ over $k\ge 1$ and using $\sum_{k\ge 1}\hat{T}_k=T-\hat{T}_0$ gives $\sum_{k\ge 1}v_k=(T-\hat{T}_0)/\hat{T}_0$. Substituting back yields $v_k=\hat{T}_k/\hat{T}_0=\hat{p}_k/\hat{p}_0$. The Hessian of $\ell_T$ is negative definite on $(0,\infty)^A$ whenever $\hat{T}_0>0$, so this critical point is the unique maximizer. Coordinates with $\hat{T}_k=0$ are maximized at $v_k=0$, which is the extended MLE on $[0,\infty)^A$. Every coordinate is strictly positive if and only if $\hat{T}_i>0$ for all $i=0,1,\ldots,A$, and only then is the maximizer interior.
    
    If $\hat T_0=0$, write $v=tu$ with $t>0$ and $\sum_{i=1}^A u_i=1$. Since $\sum_{i=1}^A\hat T_i=T$, we have $\ell_T(tu)=\sum_{i=1}^A\hat T_i\log u_i+T\log(t/(1+t))$. The last term is strictly increasing in $t$ and converges to zero, so the finite supremum is not attained. \qed
    
    \subsection{Proof of Lemma~\ref{lem:assortment-l1}}
    
    The multinomial counts give $\mathbb E[\|\hat{\boldsymbol p}-p^*\|_2^2]=T^{-1}(1-\sum_{i=0}^A(p_i^*)^2)$. Since $\sum_{i=0}^A(p_i^*)^2\ge 1/(A+1)$, this is at most $A/((A+1)T)$. Cauchy--Schwarz therefore gives
    \[
    \mathbb{E}\bigl[\|\hat{\boldsymbol{p}}-p^*\|_1\bigr]
    \le
    \sqrt{A+1}\,
    \sqrt{\mathbb{E}\bigl[\|\hat{\boldsymbol{p}}-p^*\|_2^2\bigr]}
    \le
    \sqrt{\frac{A}{T}}.
    \]
    On $\{\hat{p}_0>0\}$, fix $\mathcal{S}\subseteq\mathcal{A}$ and write $a_i=p_i^*$, $\hat{a}_i=\hat{p}_i$ for $i\in\mathcal{S}\cup\{0\}$, with sums $A_{\mathcal{S}}$ and $\hat{A}_{\mathcal{S}}$. Then
    \[
    \sum_{i\in\mathcal{S}\cup\{0\}}
    \bigl|q_i(\mathcal{S};\widehat{Q})-q_i(\mathcal{S};Q_0)\bigr|
    \le
    \frac{2\|\hat{a}-a\|_1}{A_{\mathcal{S}}}
    \le
    \frac{2}{\kappa}\|\hat{\boldsymbol{p}}-p^*\|_1,
    \]
    because $A_{\mathcal{S}}\ge p_0^*\ge\kappa$. The same bound holds for the supremum over $\mathcal{S}$, so
    \[
    \mathbb{E}\!\left[
    \sup_{\mathcal{S}\subseteq\mathcal{A}}
    \sum_{i\in\mathcal{S}\cup\{0\}}
    \bigl|q_i(\mathcal{S};\widehat{Q})-q_i(\mathcal{S};Q_0)\bigr|
    \cdot\mathbf{1}\{\hat{p}_0>0\}
    \right]
    \le
    \frac{2}{\kappa}\sqrt{\frac{A}{T}}.
    \]
    On $\{\hat{p}_0=0\}$ the error is $2$, and $\mathbb{P}(\hat{p}_0=0)=(1-p_0^*)^T\le\exp(-\kappa T)$. Adding these contributions yields \eqref{eq:assortment-l1-expectation}. \qed

    \subsection{Proof of Proposition~\ref{prop:mnl_nominal_revenue}}
    We first bound the effects of estimation and robustness separately.
    \begin{aplemma}[Nominal Revenue and Choice Probabilities]\label{lem:revenue_from_choice}
        Let $Q$ and $Q'$ be two utility distributions. For every product line $\mathcal{S}$,
        \[
        \bigl|R(\mathcal{S};Q)-R(\mathcal{S};Q')\bigr|
        \le
        \frac{\bar r}{2}
        \sum_{i\in\mathcal{S}\cup\{0\}}
        \bigl|q_i(\mathcal{S};Q)-q_i(\mathcal{S};Q')\bigr|.
        \]
    \end{aplemma}
    \begin{proof}[Proof of Lemma~\ref{lem:revenue_from_choice}]
        Let $d_i=q_i(\mathcal S;Q)-q_i(\mathcal S;Q')$ for $i\in\bar{\mathcal S}$. The outside option can be included in the revenue sum because $r_0=0$. Since both choice vectors sum to one, $\sum_{i\in\bar{\mathcal S}}d_i=0$. Subtracting $\bar r/2$ from each revenue coefficient therefore leaves the difference unchanged:
        \[
        \bigl|R(\mathcal S;Q)-R(\mathcal S;Q')\bigr|
        =\left|\sum_{i\in\bar{\mathcal S}}(r_i-\bar r/2)d_i\right|
        \le \frac{\bar r}{2}\sum_{i\in\bar{\mathcal S}}|d_i|,
        \]
        where the inequality uses $0\le r_i\le\bar r$.
    \end{proof}
    \begin{aplemma}[Robustness Penalty for an MNL Nominal]\label{lem:mnl_robustness_penalty}
        Let $Q$ be any MNL distribution. Then, for every $\epsilon\ge 0$,
        \begin{equation}
        \label{eq:mnl-robustness-penalty}
        V_0^\star(Q)-V_\epsilon^\star(Q)
        \le
        \bar r\sqrt{\frac{\epsilon}{2}}.
        \end{equation}
        \end{aplemma}
        
        \begin{proof}[Proof of Lemma~\ref{lem:mnl_robustness_penalty}]
        The claim is immediate when $\epsilon=0$. Fix $\epsilon>0$ and a product line $\mathcal{S}$. Let $X$ be the nominal revenue of the alternative chosen from $\mathcal{S}\cup\{0\}$, so $X\in[0,\bar r]$ and $\mathbb{E}[X]=R(\mathcal{S};Q)$. The Gumbel calculation in \eqref{eq:rho-v} gives
        \[
        R_\epsilon(\mathcal{S};Q)
        =
        \sup_{\lambda>0}
        \left\{
        -\lambda\epsilon-\lambda\log\mathbb{E}\!\left[e^{-X/\lambda}\right]
        \right\}.
        \]
        Hoeffding's lemma gives $\log\mathbb E[e^{-X/\lambda}]\le -R(\mathcal S;Q)/\lambda+\bar r^2/(8\lambda^2)$ for every $\lambda>0$. Multiplying by $-\lambda$ and substituting yields
        \[
        R_\epsilon(\mathcal{S};Q)
        \ge
        R(\mathcal{S};Q)
        -
        \inf_{\lambda>0}
        \left\{
        \lambda\epsilon+\frac{\bar r^2}{8\lambda}
        \right\}
        =
        R(\mathcal{S};Q)-\bar r\sqrt{\frac{\epsilon}{2}}.
        \]
        The infimum equals $\bar r\sqrt{\epsilon/2}$ at $\lambda=\bar r/\sqrt{8\epsilon}$. Evaluating the inequality at a nominal optimizer gives \eqref{eq:mnl-robustness-penalty}.
        \end{proof}
    
    Let $\widehat{\mathcal{S}}=\widehat{\mathcal{S}}_{T,\epsilon}$. On $\{\hat{p}_0=0\}$ the maximum-likelihood estimator is undefined and the loss is at most $\bar r$. On $\{\hat{p}_0>0\}$, decompose the loss as
    \[
    V_0^\star-R(\widehat{\mathcal{S}};Q_0)
    =
    \bigl(V_0^\star(Q_0)-V_0^\star(\widehat{Q})\bigr)
    +\bigl(V_0^\star(\widehat{Q})-V_\epsilon^\star(\widehat{Q})\bigr)
    +\bigl(V_\epsilon^\star(\widehat{Q})-R(\widehat{\mathcal{S}};Q_0)\bigr).
    \]
    Lemma~\ref{lem:mnl_robustness_penalty} bounds the second term by $\bar r\sqrt{\epsilon/2}$.
    
    For the first term, let $\mathcal{S}^\star$ be optimal for $R(\cdot;Q_0)$. Lemma~\ref{lem:revenue_from_choice} gives
    \[
    V_0^\star(Q_0)-V_0^\star(\widehat{Q})
    \le
    R(\mathcal{S}^\star;Q_0)-R(\mathcal{S}^\star;\widehat{Q})
    \le
    \frac{\bar r}{2}
    \sum_{i\in\mathcal{S}^\star\cup\{0\}}
    \bigl|q_i(\mathcal{S}^\star;\widehat{Q})-q_i(\mathcal{S}^\star;Q_0)\bigr|.
    \]
    For the third term, optimality of $\widehat{\mathcal{S}}$ for the robust objective and the inequality $R_\epsilon\le R$ give $V_\epsilon^\star(\widehat{Q})=R_\epsilon(\widehat{\mathcal{S}};\widehat{Q})\le R(\widehat{\mathcal{S}};\widehat{Q})$. Lemma~\ref{lem:revenue_from_choice} therefore yields
    \[
    V_\epsilon^\star(\widehat{Q})-R(\widehat{\mathcal{S}};Q_0)
    \le
    \frac{\bar r}{2}
    \sum_{i\in\widehat{\mathcal{S}}\cup\{0\}}
    \bigl|q_i(\widehat{\mathcal{S}};\widehat{Q})-q_i(\widehat{\mathcal{S}};Q_0)\bigr|.
    \]
    The two sums are each at most the supremum in Lemma~\ref{lem:assortment-l1}, so together they contribute at most $\bar r$ times that choice-probability error. On $\{\hat{p}_0>0\}$,
    \[
    V_0^\star-R(\widehat{\mathcal{S}};Q_0)
    \le
    \bar r\sqrt{\frac{\epsilon}{2}}
    +\bar r
    \sup_{\mathcal{S}\subseteq\mathcal{A}}
    \sum_{i\in\mathcal{S}\cup\{0\}}
    \bigl|q_i(\mathcal{S};\widehat{Q})-q_i(\mathcal{S};Q_0)\bigr|.
    \]
    On $\{\hat{p}_0=0\}$ the loss is at most $\bar r$, and $\mathbb{P}(\hat{p}_0=0)\le\exp(-\kappa T)$. Therefore
    \[
    \mathbb{E}\!\left[V_0^\star-R(\widehat{\mathcal{S}};Q_0)\right]
    \le
    \bar r\sqrt{\frac{\epsilon}{2}}
    +\bar r\,\mathbb{E}\!\left[
    \sup_{\mathcal{S}\subseteq\mathcal{A}}
    \sum_{i\in\mathcal{S}\cup\{0\}}
    \bigl|q_i(\mathcal{S};\widehat{Q})-q_i(\mathcal{S};Q_0)\bigr|
    \cdot\mathbf{1}\{\hat{p}_0>0\}
    \right]
    +\bar r\exp(-\kappa T).
    \]
    The proof of Lemma~\ref{lem:assortment-l1} bounds the expectation on $\{\hat{p}_0>0\}$ by $(2/\kappa)\sqrt{A/T}$. Substituting that bound gives \eqref{eq:mnl_expected_loss}. The displayed expectation in the lemma is larger by the additional $2\exp(-\kappa T)$ coming from the convention that the error equals $2$ on $\{\hat{p}_0=0\}$; the argument above charges that event only $\bar r$. \qed
    
    \subsection{Proof of Corollary~\ref{cor:mnl_sample_size}}
    The threshold $16A\bar r^2/(\kappa^2\eta^2)$ makes $(2\bar r/\kappa)\sqrt{A/T}\le\eta/2$. The threshold $\kappa^{-1}\log(2\bar r/\eta)$ makes $\bar r\exp(-\kappa T)\le\eta/2$ when $2\bar r>\eta$. When $\eta\ge 2\bar r$, one has $\log(2\bar r/\eta)\le 0$, so \eqref{eq:mnl_sample_size} imposes no further restriction from this term, and $\bar r\exp(-\kappa T)\le\bar r\le\eta/2$ for every $T\ge 0$. Adding these estimation terms to the robustness penalty gives the claim. \qed
    
    \subsection{Proof of Proposition~\ref{lemma:2}}
    
    The minimizing property of $k_n$ gives
    \begin{equation}\label{eq:kn_minimality}
    r_{k_n}-\lambda\hat u_{nk_n}
    \le r_i-\lambda\hat u_{ni},
    \qquad i\in\mathcal S\cup\{0\}.
    \end{equation}
    If $r_i<r_{k_n}$ and $\hat u_{ni}\ge\hat u_{nk_n}$, the right-hand side of~\eqref{eq:kn_minimality} is strictly smaller than the left-hand side, a contradiction. This proves dominance exclusion.
    
    If $r_i>r_{k_n}$ and $\hat u_{ni}>\hat u_{nk_n}$, rearranging~\eqref{eq:kn_minimality} gives the upper-margin condition $(r_i-r_{k_n})/(\hat u_{ni}-\hat u_{nk_n})\ge\lambda$. If $r_i<r_{k_n}$ and $\hat u_{ni}<\hat u_{nk_n}$, it gives the lower-margin condition $(r_{k_n}-r_i)/(\hat u_{nk_n}-\hat u_{ni})\le\lambda$. \qed
    
    \subsection{Proof of Proposition~\ref{pro:mon_lam}}
    
    Let
    \[
    F(\mathcal S,\lambda)
    =\lambda W_N(\mathcal S)
    +\frac{1}{N}\sum_{n\in\mathcal N}
    \min_{i\in\bar{\mathcal S}}\{r_i-\lambda\hat u_{ni}\}
    \]
    denote the part of the [DRO-FC] objective that does not contain $\epsilon$. For $\epsilon_1<\epsilon_2$, let $(\mathcal S_1,\lambda_1)$ and $(\mathcal S_2,\lambda_2)$ be joint optimizers. Their optimality gives
    \begin{align*}
    F(\mathcal S_1,\lambda_1)-\epsilon_1\lambda_1
    &\ge F(\mathcal S_2,\lambda_2)-\epsilon_1\lambda_2,\\
    F(\mathcal S_2,\lambda_2)-\epsilon_2\lambda_2
    &\ge F(\mathcal S_1,\lambda_1)-\epsilon_2\lambda_1.
    \end{align*}
    Adding the inequalities yields
    $(\epsilon_2-\epsilon_1)(\lambda_1-\lambda_2)\ge0$, and therefore $\lambda_1\ge\lambda_2$. \qed
    
    \subsection{Proof of Lemma~\ref{lem:envelope}}
    For convenience let
    $\bar {\mathcal S}={\mathcal S}\cup \{0\}$.
    For every $i\in\bar {\mathcal S}$,
    $\max_{j\in\bar {\mathcal S}}U_j-U_i\ge0$, so the second expression in
    \eqref{eq:phi} is nonnegative.  Choosing an index $i$ that maximizes utility
    gives $\phi_{{\mathcal S},\lambda}(U)\le r_i\le\bar r$.  Consequently, for
    $\lambda>\bar r/\epsilon$ the dual criterion is at most
    $-\lambda\epsilon+\bar r<0$.  Since $\lambda=0$ attains the nonnegative
    continuous-extension value, no optimizer lies beyond
    $\bar r/\epsilon$. \qed

    \subsection{Proof of Theorem~\ref{thm:loss}}

    For $\mathcal S\in\mathcal L$, let $Z_{N,\mathcal S}=\sup_{\lambda\ge0}|\hat\mu_{\mathcal S}(\lambda)-\mu_{\mathcal S}(\lambda)|$, so $Z_N=\max_{\mathcal S\in\mathcal L}Z_{N,\mathcal S}$. The inequality
    $\bigl|\sup_\lambda a(\lambda)-\sup_\lambda b(\lambda)\bigr|\le\sup_\lambda|a(\lambda)-b(\lambda)|$
    gives, for every $\epsilon>0$ and every $\mathcal{S}\in\mathcal{L}$,
    \[
    \bigl|R_\epsilon(\mathcal{S};Q)-\hat R_\epsilon(\mathcal{S})\bigr|
    \le Z_{N,\mathcal{S}}\le Z_N.
    \]
    The penalty $-\lambda\epsilon$ does not appear in $Z_N$, so this comparison is uniform in $\epsilon$. Let $\mathcal{S}^*$ maximize $R_\epsilon(\cdot;Q)$ and recall that $\Shat$ maximizes $\hat R_\epsilon$ over $\mathcal{L}$. Then
    \begin{align*}
    L_N
    &=
    \bigl[R_\epsilon(\mathcal{S}^*;Q)-\hat R_\epsilon(\mathcal{S}^*)\bigr]
    +\bigl[\hat R_\epsilon(\mathcal{S}^*)-\hat R_\epsilon(\Shat)\bigr]
    +\bigl[\hat R_\epsilon(\Shat)-R_\epsilon(\Shat;Q)\bigr]\\
    &\le Z_N+0+Z_N=2Z_N.
    \end{align*}
    The same comparison of maxima gives $\bigl|V^{\star}_{N,\epsilon}-V_\epsilon^\star(Q)\bigr|\le Z_N$. Thus \eqref{eq:value_statistical_error} and \eqref{eq:decision_statistical_error} follow from \eqref{eq:uniform_statistical_error}.

    It remains to control $Z_N$. Let $\widehat Q_N$ denote the empirical distribution of the $N$ utility vectors. Fix $\mathcal{S}$ and set $f_\lambda=\phi_{\mathcal{S},\lambda}/\bar r$. By Lemma~\ref{lem:envelope}, $f_\lambda$ takes values in $[0,1]$. Each term in the minimum \eqref{eq:phi} is affine in $\lambda$ with nonnegative slope, so $f_\lambda$ is continuous and nondecreasing in $\lambda$. For $t\in[0,1]$ define $\tau_t(\boldsymbol{U})=\inf\{\lambda\ge0:f_\lambda(\boldsymbol{U})>t\}$, with the infimum of the empty set equal to $+\infty$. The level sets
    \[
    E_{\lambda,t}=\{\boldsymbol{U}:f_\lambda(\boldsymbol{U})>t\}
    =\{\boldsymbol{U}:\tau_t(\boldsymbol{U})<\lambda\}
    \]
    are thresholds of the scalar random variable $\tau_t$. The Dvoretzky--Kiefer--Wolfowitz inequality with Massart's constant \citep{massart1990tight} therefore yields, for every $x>0$,
    \[
    \mathbb{P}\!\left\{
    \sup_{\lambda\ge0}
    \bigl|\widehat Q_N(E_{\lambda,t})-Q(E_{\lambda,t})\bigr|>x
    \right\}
    \le 2\exp(-2Nx^2).
    \]
    The layer-cake representation gives
    \[
    \frac{Z_{N,\mathcal{S}}}{\bar r}
    \le
    \int_0^1
    \sup_{\lambda\ge0}
    \bigl|\widehat Q_N(E_{\lambda,t})-Q(E_{\lambda,t})\bigr|\,dt.
    \]
    Taking expectations and applying Tonelli's theorem,
    \[
    \frac{\mathbb{E}[Z_{N,\mathcal{S}}]}{\bar r}
    \le
    \int_0^1\int_0^\infty 2e^{-2Nx^2}\,dx\,dt
    =
    \sqrt{\frac{\pi}{2N}}.
    \]
    Replacing one sampled utility vector changes $\hat\mu_{\mathcal{S}}(\lambda)$ by at most $\bar r/N$ for every $\lambda$, and hence changes $Z_{N,\mathcal{S}}$ by at most $\bar r/N$. Since $\mathbb{E}[Z_{N,\mathcal{S}}]\le\bar r\sqrt{\pi/(2N)}$, the bounded-difference inequality \citep{mcdiarmid1989method} gives
    \[
    \mathbb{P}\!\left\{
    Z_{N,\mathcal{S}}
    >
    \bar r\sqrt{\frac{\pi}{2N}}+x
    \right\}
    \le
    \exp\!\left(-\frac{2Nx^2}{\bar r^2}\right).
    \]
    Write $M=|\mathcal L|$, $\mu=\bar r\sqrt{\pi/(2N)}$, and $a=2N/\bar r^2$. A union bound gives $\mathbb P(Z_N>\mu+x)\le M\exp(-ax^2)$ for $x>0$. Taking $x=\sqrt{\log(M/\alpha)/a}$ and using $Z_N\le\bar r$ proves \eqref{eq:uniform_statistical_error}.

    To bound the expectation, set $t_0=\sqrt{\log(eM)/a}$, so $M\exp(-at_0^2)=1/e$. Integrating the tail gives
    \begin{align*}
    \mathbb{E}[Z_N]
    &\le
    \mu+t_0+\int_{t_0}^\infty M\exp(-ax^2)\,dx\\
    &\le
    \mu+t_0+\frac{1}{2at_0}.
    \end{align*}
    The integral bound uses $\int_{t_0}^\infty\exp(-ax^2)\,dx\le\exp(-at_0^2)/(2at_0)$. Since $\log(eM)\ge1$, one has $1/(2at_0)\le\bar r/\sqrt{8N}$. Substituting $t_0=\bar r\sqrt{\log(eM)/(2N)}$ therefore yields
    \[
    \mathbb{E}[Z_N]
    \le
    \frac{\bar r}{\sqrt{2N}}
    \left(
    \sqrt{\pi}+\frac12+\sqrt{1+\log|\mathcal{L}|}
    \right).
    \]
    Multiplying by two and using $L_N\le 2Z_N$ gives \eqref{eq:expected-loss}. \qed

    \subsection{Proof of Lemma~\ref{prop:lambda_bounds}}
    
    Since $\lambda \ge 0$ in the dual reformulation of Proposition~\ref{prop:dual_worst_case_revenue}, only the upper bound requires proof. We establish it via the exact right-derivative envelope of $h$.
    
    \emph{Right-derivative envelope.}
    Recall that $h(\lambda) = -\lambda\epsilon + f(\lambda) = -\lambda\epsilon + \max_{\mathcal{S}\in\mathcal{L}} G(\mathcal{S},\lambda)$. For each fixed $\mathcal{S}\in\mathcal{L}$, $G(\mathcal{S},\cdot)$ is concave and piecewise linear, being an average of pointwise minima of affine functions, and the customer-wise worst-case selection $k_n(\lambda) \in \arg\min_{i\in\bar{\mathcal{S}}}(r_i - \lambda\hat{u}_{ni})$ is well defined (breaking ties toward the utility-maximizer). Consequently
    \[
        h'_\mathcal{S}(\lambda^+) \;=\; -\epsilon + \tfrac{1}{N}\sum_{n\in\mathcal{N}} \bigl[\hat{u}_{n,i_n^\star(\mathcal{S})} - \hat{u}_{n,k_n(\mathcal{S},\lambda)}\bigr],
    \]
    where $i_n^\star(\mathcal{S}) \in \arg\max_{i\in\bar{\mathcal{S}}}\hat{u}_{ni}$ realizes $W_N(\mathcal{S})$.
    
    \emph{Customer-separable bound on the envelope.}
    Fix $\mathcal{S}$ and $n$. If $i_n^\star(\mathcal{S}) = 0$, then $\hat{u}_{n0} \ge \hat{u}_{ni}$ for all $i \in \mathcal{S}$ and $r_0 = 0 \le r_i$, so $r_0 - \lambda\hat{u}_{n0} \le r_i - \lambda\hat{u}_{ni}$ for all $i\in\mathcal{S}$; hence the outside option is itself a worst-case minimizer and $\hat{u}_{n,i_n^\star(\mathcal{S})} - \hat{u}_{n,k_n(\mathcal{S},\lambda)} = 0 \le D_n(\lambda)$. Otherwise, $i_n^\star(\mathcal{S}) \in \mathcal{A}$ and, by the minimality of $k_n(\mathcal{S},\lambda)$,
    \[
        r_{k_n} - \lambda\hat{u}_{n,k_n} \;\le\; r_{i_n^\star} - \lambda\hat{u}_{n,i_n^\star}.
    \]
    Rearranging gives $r_{i_n^\star} - r_{k_n} \ge \lambda\,(\hat{u}_{n,i_n^\star} - \hat{u}_{n,k_n})$. The utility gap $\hat{u}_{n,i_n^\star} - \hat{u}_{n,k_n}$ is non-negative by construction; if it vanishes, the bound $0 \le D_n(\lambda)$ is trivial. If it is strictly positive, then $r_{i_n^\star} > r_{k_n}$ (immediate from the display for $\lambda>0$; at $\lambda=0$ the minimizer is the outside option, so $r_{k_n}=0<r_{i_n^\star}$). Moreover, equality in the ratio would make $i_n^\star$ and $k_n$ simultaneous minimizers, in which case the right-derivative tie rule selects the higher-utility alternative $i_n^\star$ and the utility gap would be zero. Thus
    $\hat{u}_{n,i_n^\star} > \hat{u}_{n,k_n}$ and
    $(r_{i_n^\star}-r_{k_n})/(\hat{u}_{n,i_n^\star}-\hat{u}_{n,k_n}) > \lambda$, so the pair $(i_n^\star, k_n)$ is feasible for the maximization defining $D_n(\lambda)$ and therefore
    \[
        \hat{u}_{n,i_n^\star(\mathcal{S})} - \hat{u}_{n,k_n(\mathcal{S},\lambda)} \;\le\; D_n(\lambda).
    \]
    Averaging the bound $\hat{u}_{n,i_n^\star} - \hat{u}_{n,k_n} \le D_n(\lambda)$ over $n$ and taking the supremum over $\mathcal{S}\in\mathcal{L}$ yields, uniformly in $\mathcal{S}$,
    \begin{equation}\label{eq:envelope_bound}
        \sup_{\mathcal{S}\in\mathcal{L}} h'_\mathcal{S}(\lambda^+) \;\le\; \Phi(\lambda) - \epsilon.
    \end{equation}
    
    \emph{Monotonicity of $\Phi$.} As $\lambda$ increases, the feasibility set $\{(i,j,n) : (r_i-r_j)/(\hat{u}_{ni}-\hat{u}_{nj}) > \lambda\}$ underlying~\eqref{eq:D_n} shrinks; hence each $D_n(\cdot)$ is non-increasing and so is $\Phi$. At a ratio breakpoint, the corresponding pairs are excluded at the breakpoint itself. Thus $\Phi$ is right-continuous and piecewise constant.
    
    \emph{Upper bound on $\lambda^*$.}
    Let $\lambda^u = \min\{\lambda \ge 0 : \Phi(\lambda) \le \epsilon\}$, which is well defined because $\Phi$ is right-continuous, piecewise constant, and eventually zero. By monotonicity, $\Phi(\lambda) \le \epsilon$ for every $\lambda\ge\lambda^u$. Equation~\eqref{eq:envelope_bound} therefore gives $h'_\mathcal{S}(\lambda^+)\le0$ for every $\mathcal S$ and every $\lambda\ge\lambda^u$. Hence each $h_\mathcal S$ is non-increasing on $[\lambda^u,\infty)$ and attains its maximum over $[0,\lambda^u]$. This proves~\eqref{eq:lambda_bounds} for every product line. Taking the maximum over $\mathcal S\in\mathcal L$ then yields a joint optimal multiplier in $[0,\lambda^u]$. \qed

    \subsection{Proof of Lemma~\ref{lem:modularity}}
    
    Monotonicity of $W_N(\cdot)$ (resp. reverse-monotonicity of $C(\cdot,\lambda)$) is immediate: $\max_{i\in\bar{\mathcal{S}}}\hat{u}_{ni}$ is non-decreasing and $\min_{i\in\bar{\mathcal{S}}}(r_i-\lambda\hat{u}_{ni})$ is non-increasing in $\mathcal{S}$; averaging preserves the direction.
    
    \emph{Submodularity of $W_N(\cdot)$.} Fix $\mathcal{S}_1\subseteq\mathcal{S}_2\subseteq\mathcal{A}$ and $j\notin\mathcal{S}_2$. For every customer~$n$,
    \[
        \max_{i\in\mathcal{S}_k\cup\{j,0\}}\hat{u}_{ni} - \max_{i\in\mathcal{S}_k\cup\{0\}}\hat{u}_{ni} \;=\; \bigl[\hat{u}_{nj} - \max_{i\in\mathcal{S}_k\cup\{0\}}\hat{u}_{ni}\bigr]_+ \qquad (k=1,2).
    \]
    Because $\max_{i\in\mathcal{S}_1\cup\{0\}}\hat{u}_{ni} \le \max_{i\in\mathcal{S}_2\cup\{0\}}\hat{u}_{ni}$, the positive-part expression is non-increasing in $\mathcal{S}$. Averaging over $n$ and taking the difference gives
    \[
        W_N(\mathcal{S}_1\cup\{j\}) - W_N(\mathcal{S}_1) \;\ge\; W_N(\mathcal{S}_2\cup\{j\}) - W_N(\mathcal{S}_2),
    \]
    which is the submodularity of $W_N(\cdot)$.
    
    \emph{Supermodularity of $C(\cdot,\lambda)$.} By the same argument applied to the minimum, for every~$n$,
    \[
        \min_{i\in\mathcal{S}_k\cup\{0\}}(r_i - \lambda\hat{u}_{ni}) - \min_{i\in\mathcal{S}_k\cup\{j,0\}}(r_i - \lambda\hat{u}_{ni}) \;=\; \bigl[\min_{i\in\mathcal{S}_k\cup\{0\}}(r_i - \lambda\hat{u}_{ni}) - (r_j - \lambda\hat{u}_{nj})\bigr]_+.
    \]
    Because $\min_{i\in\mathcal{S}_1\cup\{0\}}(r_i-\lambda\hat{u}_{ni}) \ge \min_{i\in\mathcal{S}_2\cup\{0\}}(r_i-\lambda\hat{u}_{ni})$, this positive part is again non-increasing in $\mathcal{S}$. Averaging and rearranging yields
    \[
        C(\mathcal{S}_1,\lambda) - C(\mathcal{S}_1\cup\{j\},\lambda) \;\ge\; C(\mathcal{S}_2,\lambda) - C(\mathcal{S}_2\cup\{j\},\lambda),
    \]
    i.e., $-C(\cdot,\lambda)$ is submodular, so $C(\cdot,\lambda)$ is supermodular. \qed

    \subsection{Proof of Proposition~\ref{prop:modular_ub}}
    Fix a completion \(\mathcal{T}\subseteq\mathcal{F}\) of size
    \(t\le k\). By submodularity of \(W_N(\cdot)\),
    \[
    W_N(\mathcal{I}\cup\mathcal{T})-W_N(\mathcal I)
    \le
    \sum_{j\in\mathcal{T}}
    \bigl(W_N(\mathcal{I}\cup\{j\})-W_N(\mathcal I)\bigr)
    =
    \sum_{j\in\mathcal{T}}m^W_j(\mathcal I)
    \le
    \sum_{q=1}^{t}m^W_{[q]}.
    \]
    
    For the $C(\cdot,\lambda)$ term, direct calculation gives $C(\mathcal I,\lambda)-C(\mathcal I\cup\{j\},\lambda)=m^C_j(\mathcal I,\lambda)$. If $t\ge1$, monotonicity then gives, for every $j\in\mathcal T$,
    \[
    C(\mathcal{I}\cup\mathcal{T},\lambda)
    \le C(\mathcal{I}\cup\{j\},\lambda)
    =C(\mathcal I,\lambda)-m^C_j(\mathcal I,\lambda).
    \]
    Thus
    \[
    C(\mathcal{I}\cup\mathcal{T},\lambda)
    \le C(\mathcal I,\lambda)-\max_{j\in\mathcal{T}}m^C_j(\mathcal I,\lambda)
    \le C(\mathcal I,\lambda)-m^C_{(t)}.
    \]
    The last inequality follows because a set of \(t\) free products must contain a product whose
    decrease in the $C(\cdot,\lambda)$ term is at least the \(t\)-th smallest free-product decrease. For \(t=0\), the same expression
    holds under the convention \(m^C_{(0)}=0\).
    
    Combining the welfare and cost bounds yields
    \[
    G(\mathcal{I}\cup\mathcal{T},\lambda)
    \le\mathrm{UB}(t\mid\mathcal I,\mathcal F)
    \le\mathrm{UB}_{\mathrm{mod}}(\mathcal I,\mathcal F).
    \]
    Since every admissible completion has some size
    \(t\in\{0,\ldots,\min(k,|\mathcal{F}|)\}\), the modular component is valid for the at-most-$p$ model.
    
    It remains to verify the separable component. For customer $n$, write
    \[
    d^u_{nj}=[\hat u_{nj}-\bar u_n(\mathcal I)]_+,
    \qquad
    d^c_{nj}=[\bar c_n(\mathcal I,\lambda)-(r_j-\lambda\hat u_{nj})]_+.
    \]
    For a nonempty completion $\mathcal T$, take $j_n$ maximizing $d^u_{nj}$ over $\mathcal T$. The change in customer $n$'s contribution satisfies
    \[
    \lambda\max_{j\in\mathcal T}d^u_{nj}
    -\max_{j\in\mathcal T}d^c_{nj}
    \le \lambda d^u_{nj_n}-d^c_{nj_n}
    \le \max_{j\in\mathcal T}\eta_{nj}(\mathcal I,\lambda)
    \le \left[\max_{j\in\mathcal F}\eta_{nj}(\mathcal I,\lambda)\right]_+.
    \]
    The last positive part also covers the empty completion. Summing over customers and dividing by
    $N$ proves that $\mathrm{UB}_{\mathrm{sep}}(\mathcal I,\mathcal F)$ is valid. The minimum of
    two valid upper bounds is valid, establishing~\eqref{eq:node_sep_ub} and the proposition.
    
    Finally, SBMB evaluates \(G(\mathcal{I},\lambda)\) at every visited node with
    \(\mathcal I\ne\emptyset\), so nonempty feasible product lines of all cardinalities can update the
    incumbent. The empty root need not be evaluated because $G(\emptyset,\lambda)=0$ and the
    best-singleton initialization is nonnegative. Every child preserves \(|\mathcal{I}|\le p\), and the
    include/exclude branching tree is finite. A node is discarded only when its valid upper bound does
    not exceed the incumbent. Because SBMB processes nodes in decreasing order of their bounds, it
    terminates only when the queue is empty or its largest bound does not exceed the incumbent. No
    unexplored descendant can then improve the incumbent, so the returned value equals $f(\lambda)$.
    \qed

    \subsection{Proof of Lemma~\ref{lem:bilipschitz}}
    
    Consider any feasible $\mathcal{S} \in \mathcal{L}$. For each $n\in\mathcal{N}$, $\min_{i\in\bar{\mathcal{S}}}(r_i - \lambda \hat{u}_{ni})$ is the pointwise minimum of finitely many affine functions of $\lambda$ and is therefore concave and piecewise linear, so $G(\mathcal{S},\cdot)$ is concave in $\lambda$. Except on a finite set of breakpoints where the argmin is set-valued, the right derivative equals
    \[
        G'(\mathcal{S}, \lambda^+) \;=\; W_N(\mathcal{S}) \; - \; \tfrac{1}{N}\sum_{n \in \mathcal{N}} \hat{u}_{n, k_n(\lambda)},
    \]
    where $k_n(\lambda) \in \arg\min_{i \in \bar{\mathcal{S}}} (r_i - \lambda \hat{u}_{ni})$ is a worst-case product for customer~$n$. Since $k_n(\lambda) \in \bar{\mathcal{S}}$, we have $\hat{u}_{n, k_n(\lambda)} \le \max_{i\in\bar{\mathcal{S}}}\hat{u}_{ni}$ and hence $G'(\mathcal S,\lambda^+)\ge0$.
    
    Comparing the minimizing product with the outside option gives $r_{k_n(\lambda)}-\lambda\hat u_{n,k_n(\lambda)}\le -\lambda\hat u_{n0}$. If $k_n(\lambda)\ne0$, then $r_{k_n(\lambda)}>0$, which requires $\lambda>0$ and $\hat u_{n,k_n(\lambda)}>\hat u_{n0}$. The utilities are equal if $k_n(\lambda)=0$. Together with $W_N(\mathcal S)\le N^{-1}\sum_n\max_{i\in\Aplus}\hat u_{ni}$, this gives
    \[
        0 \;\le\; G'(\mathcal{S},\lambda^+) \;\le\; \tfrac{1}{N}\sum_n\bigl(\max_{i\in\Aplus}\hat{u}_{ni} - \hat u_{n0}\bigr) \;=\; \tfrac{1}{N}\sum_n\Delta_n.
    \]
    For $\lambda_1<\lambda_2$ we therefore have $0\le G(\mathcal{S},\lambda_2)-G(\mathcal{S},\lambda_1) \le (\tfrac{1}{N}\sum_n\Delta_n)(\lambda_2-\lambda_1)$. Taking the pointwise maximum over $\mathcal{S}\in\mathcal{L}$ preserves both the non-negative lower bound and the Lipschitz upper bound: for any $\mathcal{S}_1^\star\in\arg\max_\mathcal{S} G(\mathcal{S},\lambda_1)$,
    \[
        f(\lambda_2) - f(\lambda_1) \;\ge\; G(\mathcal{S}_1^\star,\lambda_2) - G(\mathcal{S}_1^\star,\lambda_1) \;\ge\; 0,
    \]
    and
    \[
        f(\lambda_2) - f(\lambda_1) \;\le\; \max_{\mathcal{S}\in\mathcal{L}}\bigl[G(\mathcal{S},\lambda_2)-G(\mathcal{S},\lambda_1)\bigr] \;\le\; \tfrac{1}{N}\sum_n\Delta_n\cdot(\lambda_2-\lambda_1).
    \]
    Hence $h'(\lambda^+)=-\epsilon+f'(\lambda^+)\in[-\epsilon,\, \tfrac{1}{N}\sum_n\Delta_n-\epsilon]$ wherever it exists; the same bounds extend to the finite difference quotient on the entire interval, which proves~\eqref{eq:bilipschitz}.
    
    For $\lambda\in[a,b]$, the slope bounds give $h(\lambda)\le h(a)+L^+(\lambda-a)$ and $h(\lambda)\le h(b)+L^-(b-\lambda)$. Thus $h$ lies below the pointwise minimum of an increasing and a decreasing affine function. Their
    intersection is $a+x_{ab}$; Equation~\eqref{eq:bilipschitz} places it in $[a,b]$. Substitution gives
    the closed form in~\eqref{eq:ub_star}. If $L^+\le0$, $h$ is non-increasing and
    $\lambda=0$ is optimal. \qed

    \subsection{Proof of Theorem~\ref{prop:delta_optimal}}
    Every product line entering the solution pool is feasible. The affine pieces of
    $h_{\mathcal S}(\cdot)$ can change only at the pairwise intersections collected in
    $\mathcal B_{\mathcal S}$, and a concave piecewise-linear function attains its maximum over a
    compact interval at an endpoint or a breakpoint. Lemma~\ref{prop:lambda_bounds} shows that
    $[0,\lambda^u]$ contains an optimal multiplier for every fixed product line. Consequently,
    [Pool-BP($\mathcal S$)] is exact and APS retains a line $\widehat{\mathcal S}$ satisfying
    \[
    R^{\mathrm{FC}}_{N,\epsilon}(\widehat{\mathcal{S}})
    =\max_{\lambda\ge0}h_{\widehat{\mathcal{S}}}(\lambda)
    \;=\;H_{\widehat{\mathcal S}}
    \;=\;h_{\mathrm{best}},
    \]
    so $h_{\mathrm{best}}$ is a valid lower bound on $V^{\star}_{N,\epsilon}$. Every SBMB call returns
    the exact value $h(\lambda)$, so Lemma~\ref{lem:bilipschitz} makes
    $\mathrm{UB}^{*}([a,b])$ a valid upper bound for every active interval.
    
    Let \(\lambda^{\star}\) maximize \(h\) on \([0,\lambda^u]\). If the active interval containing
    \(\lambda^{\star}\) has not been pruned, then
    \[
    V^{\star}_{N,\epsilon}=h(\lambda^{\star})
    \le \max_{[a,b]\in\IAPS}\mathrm{UB}^{*}([a,b])
    =:\mathrm{ub}_{\max}.
    \]
    If that interval was pruned earlier, its upper bound did not exceed the incumbent at that time, so
    the current incumbent already satisfies \(h_{\mathrm{best}}\ge V^{\star}_{N,\epsilon}\ge(1-\delta)V^{\star}_{N,\epsilon}\). In the unpruned case, the termination test $\mathrm{ub}_{\max}-h_{\mathrm{best}}\le\delta\,\mathrm{ub}_{\max}$ implies
    \[
    h_{\mathrm{best}}\ge(1-\delta)\mathrm{ub}_{\max}
    \ge(1-\delta)V^{\star}_{N,\epsilon}.
    \]
    
    It remains to prove the stated evaluation bound. Let
    $\mathcal{S}^{\mathrm W}\in\arg\max_{\mathcal{S}\in\mathcal{L}}W_N(\mathcal{S})$. Under the empirical distribution, Corollary~\ref{corollary:welfare_lb} gives
    \[
    \epsdag=W_N(\mathcal{S}^{\mathrm W})-\frac{1}{N}\sum_{n\in\mathcal{N}}\hat u_{n0}
    \le \frac{1}{N}\sum_{n\in\mathcal{N}}\Delta_n.
    \]
    Thus $\epsilon<\epsdag$ implies $L^+>0$, while the assumption $\epsilon>0$ gives
    $L^-=\epsilon>0$. First note that $V^{\star}_{N,\epsilon}>0$. Indeed,
    $h_{\mathcal S^{\mathrm W}}(0)=0$, and because
    $r_0=0<r_i$ for $i\in\mathcal A$, its right derivative at zero is
    \[
    h'_{\mathcal{S}^{\mathrm W}}(0^+)
    =-\epsilon+W_N(\mathcal{S}^{\mathrm W})-\frac{1}{N}\sum_{n\in\mathcal{N}}\hat u_{n0}
    =\epsdag-\epsilon>0.
    \]
    
    A direct calculation from~\eqref{eq:ub_star} gives, for every active interval $[a,b]$,
    \[
    \mathrm{UB}^{*}([a,b])-\max\{h(a),h(b)\}
    \le
    K(b-a),
    \qquad
    K=\frac{L^+L^-}{L^++L^-}.
    \]
    APS always selects an interval with the largest active upper bound, and both endpoints of that
    interval have already been evaluated. Hence, if the selected interval has width $b-a$,
    \[
    \mathrm{ub}_{\max}-h_{\mathrm{best}}
    \le
    \mathrm{UB}^{*}([a,b])-\max\{h(a),h(b)\}
    \le K(b-a).
    \]
    Define
    \[
    x_\delta=\frac{K\lambda^u}{\delta V^{\star}_{N,\epsilon}}\ge 0,
    \qquad
    d_\delta=\left\lceil\log_2\max\{1,x_\delta\}\right\rceil.
    \]
    An interval at depth $q$ in the midpoint-bisection tree has width $\lambda^u/2^q$. By the
    definition of $d_\delta$, every interval at depth $q\ge d_\delta$ satisfies
    \[
    K(b-a)\le \frac{K\lambda^u}{2^{d_\delta}}
    \le \delta V^{\star}_{N,\epsilon}.
    \]
    If such an interval were selected, either the interval containing an optimizer had already been
    pruned, in which case the incumbent is optimal, or $\mathrm{ub}_{\max}\ge
    V^{\star}_{N,\epsilon}$. In the latter case,
    \[
    \mathrm{ub}_{\max}-h_{\mathrm{best}}
    \le \delta V^{\star}_{N,\epsilon}
    \le \delta\,\mathrm{ub}_{\max},
    \]
    so APS would stop before bisecting that interval. Consequently, APS can bisect only nodes at
    depths $q<d_\delta$. The complete binary tree contains $2^{d_\delta}-1$ such nodes; this
    quantity is zero when $d_\delta=0$. Pruning can only reduce this count.
    Each bisection requires one new exact fixed-$\lambda$ evaluation, in addition to the two initial
    endpoint evaluations, yielding $B\le 2^{d_\delta}+1$. Since
    $d_\delta=\lceil\log_2\max\{1,x_\delta\}\rceil$, we have $2^{d_\delta}<2\max\{1,x_\delta\}$. Thus, on every empirical sample,
    \[
    V^{\star}_{N,\epsilon}-h_{\mathrm{best}}
    \le
    \delta V^{\star}_{N,\epsilon}
    ,\qquad
    B<2\max\left\{1,\frac{K\lambda^u}{\delta V^{\star}_{N,\epsilon}}\right\}+1.
    \]
    \qed
    
    \subsection{Proof of Corollary~\ref{cor:aps_population}}

    If $\epsilon<\epsilon^\dagger(\widehat Q_N)$,
    Theorem~\ref{prop:delta_optimal} gives
    $R^{\mathrm{FC}}_{N,\epsilon}(\widehat{\mathcal S})
    \ge (1-\delta)V^\star_{N,\epsilon}$.
    Otherwise, taking the outside option in each minimum in
    $C(\mathcal S,\lambda)$ gives
    \[
    h_{\mathcal S}(\lambda)
    \le \lambda\left(
    W_N(\mathcal S)-N^{-1}\sum_{n\in\mathcal N}\hat u_{n0}
    -\epsilon
    \right)
    \le 0
    \qquad(\mathcal S\in\mathcal L,\ \lambda\ge0),
    \]
    where the last inequality uses the definition of
    $\epsilon^\dagger(\widehat Q_N)$. Since $h_{\mathcal S}(0)=0$,
    the empty line is optimal in this case, so the same bound
    holds for every sample.

    By the uniform comparison following~\eqref{eq:ZN},
    $R_\epsilon(\widehat{\mathcal S};Q)\ge
    R^{\mathrm{FC}}_{N,\epsilon}(\widehat{\mathcal S})-Z_N$
    and $V^\star_{N,\epsilon}\ge V_\epsilon^\star(Q)-Z_N$.
    Combining these inequalities gives
    \[
    R_\epsilon(\widehat{\mathcal S};Q)
    \ge (1-\delta)V^\star_{N,\epsilon}-Z_N
    \ge (1-\delta)V_\epsilon^\star(Q)-(2-\delta)Z_N.
    \]
    Taking expectations and using the bound on $\mathbb E[Z_N]$
    from the proof of Theorem~\ref{thm:loss} yields
    \eqref{eq:aps_population_expectation}. \qed


\begin{thebibliography}{28}
\providecommand{\natexlab}[1]{#1}
\providecommand{\url}[1]{\texttt{#1}}
\expandafter\ifx\csname urlstyle\endcsname\relax
  \providecommand{\doi}[1]{doi: #1}\else
  \providecommand{\doi}{doi: \begingroup \urlstyle{rm}\Url}\fi

\bibitem[Bertsimas and Mi{\v{s}}i{\'c}(2017)]{bertsimas2016robust}
Dimitris Bertsimas and Velibor~V Mi{\v{s}}i{\'c}.
\newblock Robust product line design.
\newblock \emph{Operations Research}, 65\penalty0 (1):\penalty0 19--37, 2017.

\bibitem[Bertsimas and Mi{\v{s}}i{\'c}(2019)]{bertsimas2019exact}
Dimitris Bertsimas and Velibor~V Mi{\v{s}}i{\'c}.
\newblock Exact first-choice product line optimization.
\newblock \emph{Operations Research}, 67\penalty0 (3):\penalty0 651--670, 2019.

\bibitem[Blanchet and Murthy(2019)]{Blanchet2019Quantifying}
Jose Blanchet and Karthyek Murthy.
\newblock Quantifying distributional model risk via optimal transport.
\newblock \emph{Mathematics of Operations Research}, 44\penalty0 (2):\penalty0 565--600, 2019.

\bibitem[Chen and Hausman(2000)]{Chen2000Mathematical}
Kyle~D. Chen and Warren~H. Hausman.
\newblock Mathematical properties of the optimal product line selection problem using choice-based conjoint analysis.
\newblock \emph{Management Science}, 46\penalty0 (2):\penalty0 327--332, 2000.

\bibitem[D{\'e}sir et~al.(2021)D{\'e}sir, Goyal, Jiang, Xie, and Zhang]{desir2021robust}
Antoine D{\'e}sir, Vineet Goyal, Bo~Jiang, Tian Xie, and Jiawei Zhang.
\newblock Robust assortment optimization under the markov chain model.
\newblock \emph{Available at SSRN 3938924}, 2021.

\bibitem[D{\'e}sir et~al.(2024)D{\'e}sir, Goyal, Jiang, Xie, and Zhang]{desir2024robust}
Antoine D{\'e}sir, Vineet Goyal, Bo~Jiang, Tian Xie, and Jiawei Zhang.
\newblock Robust assortment optimization under the {M}arkov chain choice model.
\newblock \emph{Operations Research}, 72\penalty0 (4):\penalty0 1595--1614, 2024.

\bibitem[Esfahani and Kuhn(2018)]{esfahani2018data}
Peyman~Mohajerin Esfahani and Daniel Kuhn.
\newblock Data-driven distributionally robust optimization using the {W}asserstein metric: Performance guarantees and tractable reformulations.
\newblock \emph{Mathematical Programming}, 171\penalty0 (1-2):\penalty0 115--166, 2018.

\bibitem[Evgeniou et~al.(2005)Evgeniou, Boussios, and Zacharia]{evgeniou2005generalized}
Theodoros Evgeniou, Constantinos Boussios, and Giorgos Zacharia.
\newblock Generalized robust conjoint estimation.
\newblock \emph{Marketing Science}, 24\penalty0 (3):\penalty0 415--429, 2005.

\bibitem[Feng et~al.(2017)Feng, Li, and Wang]{feng2017relation}
Guiyun Feng, Xiaobo Li, and Zizhuo Wang.
\newblock On the relation between several discrete choice models.
\newblock \emph{Operations Research}, 65\penalty0 (6):\penalty0 1516--1525, 2017.

\bibitem[Green and Krieger(1985)]{10.2307/183706}
Paul~E. Green and Abba~M. Krieger.
\newblock Models and heuristics for product line selection.
\newblock \emph{Marketing Science}, 4\penalty0 (1):\penalty0 1--19, 1985.

\bibitem[Gupta and Kohli(1990)]{Gupta1990Designing}
Sunil Gupta and Rajeev Kohli.
\newblock Designing products and services for consumer welfare: Theoretical and empirical issues.
\newblock \emph{Marketing Science}, 9\penalty0 (3):\penalty0 230--246, 1990.

\bibitem[Jin et~al.(2025)Jin, Long, Sun, and Hu]{jin2022distributionally}
Qingwei Jin, Daniel~Zhuoyu Long, Yu~Sun, and Bin Hu.
\newblock Choice-based assortment optimization with distributional ambiguity.
\newblock \emph{Available at SSRN 4045001}, 2025.

\bibitem[Kohli and Sukumar(1990)]{KohliHeuristics1990}
Rajeev Kohli and R.~Sukumar.
\newblock Heuristics for product-line design using conjoint analysis.
\newblock \emph{Management Science}, 36\penalty0 (12):\penalty0 1464--1478, 1990.

\bibitem[Liu et~al.(2023{\natexlab{a}})Liu, Zhang, Gong, and Liu]{liu2023explainable}
Maoqi Liu, Xun Zhang, Hailei Gong, and Changchun Liu.
\newblock Explainable data-driven share-of-choice product line design optimization.
\newblock \emph{arXiv preprint arXiv:2310.07098}, 2023{\natexlab{a}}.

\bibitem[Liu et~al.(2023{\natexlab{b}})Liu, Zheng, Liu, and Zhang]{liu2023from}
Maoqi Liu, Li~Zheng, Changchun Liu, and Zhi-Hai Zhang.
\newblock From share of choice to buyers' welfare maximization: Bridging the gap through distributionally robust optimization.
\newblock \emph{Production and Operations Management}, 32\penalty0 (4):\penalty0 1205--1222, 2023{\natexlab{b}}.
\newblock \doi{10.1111/poms.13921}.

\bibitem[Massart(1990)]{massart1990tight}
Pascal Massart.
\newblock The tight constant in the {D}voretzky--{K}iefer--{W}olfowitz inequality.
\newblock \emph{The Annals of Probability}, 18\penalty0 (3):\penalty0 1269--1283, 1990.
\newblock \doi{10.1214/aop/1176990746}.

\bibitem[McDiarmid(1989)]{mcdiarmid1989method}
Colin McDiarmid.
\newblock On the method of bounded differences.
\newblock In J.~Siemons, editor, \emph{Surveys in Combinatorics, 1989}, volume 141 of \emph{London Mathematical Society Lecture Note Series}, pages 148--188. Cambridge University Press, 1989.
\newblock \doi{10.1017/CBO9781107359949.008}.

\bibitem[Mishra et~al.(2014)Mishra, Natarajan, Padmanabhan, Teo, and Li]{mishra2014theoretical}
Vinit~Kumar Mishra, Karthik Natarajan, Dhanesh Padmanabhan, Chung-Piaw Teo, and Xiaobo Li.
\newblock On theoretical and empirical aspects of marginal distribution choice models.
\newblock \emph{Management Science}, 60\penalty0 (6):\penalty0 1511--1531, 2014.

\bibitem[Mussa and Rosen(1978)]{mussa1978monopoly}
Michael Mussa and Sherwin Rosen.
\newblock Monopoly and product quality.
\newblock \emph{Journal of Economic Theory}, 18\penalty0 (2):\penalty0 301--317, 1978.

\bibitem[Natarajan et~al.(2009)Natarajan, Song, and Teo]{natarajan2009persistency}
Karthik Natarajan, Miao Song, and Chung-Piaw Teo.
\newblock Persistency model and its applications in choice modeling.
\newblock \emph{Management Science}, 55\penalty0 (3):\penalty0 453--469, 2009.

\bibitem[Rusmevichientong and Topaloglu(2012)]{rusmevichientong2012robust}
Paat Rusmevichientong and Huseyin Topaloglu.
\newblock Robust assortment optimization in revenue management under the multinomial logit choice model.
\newblock \emph{Operations Research}, 60\penalty0 (4):\penalty0 865--882, 2012.

\bibitem[Sturt(2021)]{sturt2021value}
Bradley Sturt.
\newblock The value of robust assortment optimization under ranking-based choice models.
\newblock \emph{arXiv preprint arXiv:2112.05010}, 2021.

\bibitem[Sumida et~al.(2021)Sumida, Gallego, Rusmevichientong, Topaloglu, and Davis]{sumida2021revenue}
Mika Sumida, Guillermo Gallego, Paat Rusmevichientong, Huseyin Topaloglu, and James Davis.
\newblock Revenue-utility tradeoff in assortment optimization under the multinomial logit model with totally unimodular constraints.
\newblock \emph{Management Science}, 67\penalty0 (5):\penalty0 2845--2869, 2021.

\bibitem[Toubia et~al.(2003)Toubia, Simester, Hauser, and Dahan]{toubia2003fast}
Olivier Toubia, Duncan~I Simester, John~R Hauser, and Ely Dahan.
\newblock Fast polyhedral adaptive conjoint estimation.
\newblock \emph{Marketing Science}, 22\penalty0 (3):\penalty0 273--303, 2003.

\bibitem[Toubia et~al.(2007)Toubia, Hauser, and Garcia]{Toubia2007Probabilistic}
Olivier Toubia, John Hauser, and Rosanna Garcia.
\newblock Probabilistic polyhedral methods for adaptive choice-based conjoint analysis: Theory and application.
\newblock \emph{Marketing Science}, 26\penalty0 (5):\penalty0 596--610, 2007.

\bibitem[Wang and Gutierrez(2022)]{Wang2021Robust}
Tan Wang and Genaro Gutierrez.
\newblock Robust product line design by protecting the downside while minding the upside.
\newblock \emph{Production and Operations Management}, 31\penalty0 (1):\penalty0 194--217, 2022.
\newblock \doi{10.1111/poms.13495}.
\newblock URL \url{https://onlinelibrary.wiley.com/doi/abs/10.1111/poms.13495}.

\bibitem[Wang and Curry(2012)]{Wang2012A}
Xinfang Wang and David~J. Curry.
\newblock A robust approach to the share-of-choice product design problem.
\newblock \emph{Omega}, 40\penalty0 (6):\penalty0 818--826, 2012.

\bibitem[Zufryden(1977)]{zufryden1977conjoint}
Fred~S Zufryden.
\newblock A conjoint measurement-based approach for optimal new product design and market segmentation.
\newblock In \emph{Analytic Approaches to Product and Market Planning}, pages 100--114. Marketing Science Institute, 1977.

\end{thebibliography}
\end{document}